\documentclass[11pt]{article}

\usepackage{booktabs} 
\usepackage{amsmath, amsthm, amssymb}
\usepackage{enumitem}
\usepackage{pdflscape}
\usepackage{caption}
\usepackage{bm}
\usepackage{ifpdf}
	\ifpdf
\usepackage[pdftex]{graphicx}
	\else
\usepackage[dvips]{graphicx}
	\fi
\usepackage[all]{xy}
\usepackage{tocvsec2}
\usepackage{bbm}
	\input xy
	\xyoption{all}
\usepackage[pdftex,plainpages=false,hypertexnames=false,pdfpagelabels]{hyperref}
	\newcommand{\arxiv}[1]{\href{http://arxiv.org/abs/#1}{\tt arXiv:\nolinkurl{#1}}}
	\newcommand{\arXiv}[1]{\href{http://arxiv.org/abs/#1}{\tt arXiv:\nolinkurl{#1}}}
	\newcommand{\doi}[1]{\href{http://dx.doi.org/#1}{{\tt DOI:#1}}}

	\newcommand{\googlebooks}[1]{(preview at \href{http://books.google.com/books?id=#1}{google books})}
	
\usepackage{xcolor}
	\definecolor{dark-red}{rgb}{0.7,0.25,0.25}
	\definecolor{dark-blue}{rgb}{0.15,0.15,0.55}
	\definecolor{medium-blue}{rgb}{0,0,.8}
	\definecolor{DarkGreen}{RGB}{0,150,0}
	\definecolor{rho}{named}{red}
	\hypersetup{
	   colorlinks, linkcolor={purple},
	   citecolor={medium-blue}, urlcolor={medium-blue}
	}
\usepackage{longtable}
\usepackage{fullpage}

\usepackage{mathbbol}

\theoremstyle{plain}
\newtheorem{thm}{Theorem}[section]
\newtheorem*{thm*}{Theorem}
\newtheorem{thmalpha}{Theorem}

\newtheorem{cor}[thm]{Corollary}
\newtheorem{coralpha}[thmalpha]{Corollary}
\newtheorem*{cor*}{Corollary}

\newtheorem*{conj*}{Conjecture}
\newtheorem{lem}[thm]{Lemma}
\newtheorem{prop}[thm]{Proposition}

\newtheorem*{quest*}{Question}
\newtheorem*{claim*}{Claim}

\theoremstyle{definition}
\newtheorem{defn}[thm]{Definition}

\newtheorem{construction}[thm]{Construction}

\newtheorem{exs}[thm]{Examples}
\newtheorem{ex}[thm]{Example}
\newtheorem*{ex*}{Example}
\newtheorem{sub-ex}[thm]{Sub-Example}
\newtheorem{counter-ex}[thm]{Counter-Example}

\newtheorem*{rem*}{Remark}
\newtheorem{remark}[thm]{Remark}

\DeclareMathOperator{\Aut}{Aut}
\DeclareMathOperator{\coev}{coev}

\DeclareMathOperator{\End}{End}
\DeclareMathOperator{\ev}{ev}
\DeclareMathOperator{\Hom}{Hom}

\DeclareMathOperator{\id}{id}

\DeclareMathOperator{\Spec}{\mathsf{spec}}

\DeclareMathOperator{\Tr}{Tr}
\DeclareMathOperator{\tr}{tr}

\newcommand{\comment}[1]{}

\newcommand{\noshow}[1]{}
\newcommand{\MR}[1]{}

\newcommand{\Rep}{{\sf Rep}}
\newcommand{\QAut}{{\sf QAut}}

\newcommand{\Hilb}{{\sf Hilb}}
\newcommand{\fdHilb}{{\sf Hilb_{fd}}}

\def\semicolon{;}
\def\applytolist#1{
    \expandafter\def\csname multi#1\endcsname##1{
        \def\multiack{##1}\ifx\multiack\semicolon
            \def\next{\relax}
        \else
            \csname #1\endcsname{##1}
            \def\next{\csname multi#1\endcsname}
        \fi
        \next}
    \csname multi#1\endcsname}

\def\calc#1{\expandafter\def\csname c#1\endcsname{{\mathcal #1}}}
\applytolist{calc}QWERTYUIOPLKJHGFDSAZXCVBNM;
\def\bbc#1{\expandafter\def\csname bb#1\endcsname{{\mathbb #1}}}
\applytolist{bbc}QWERTYUIOPLKJHGFDSAZXCVBNM;
\def\bfc#1{\expandafter\def\csname bf#1\endcsname{{\mathbf #1}}}
\applytolist{bfc}QWERTYUIOPLKJHGFDSAZXCVBNM;
\def\sfc#1{\expandafter\def\csname s#1\endcsname{{\sf #1}}}
\applytolist{sfc}QWERTYUIOPLKJHGFDSAZXCVBNM;
\def\fc#1{\expandafter\def\csname f#1\endcsname{{\mathfrak #1}}}
\applytolist{fc}QWERTYUIOPLKJHGFDSAZXCVBNM;

\usepackage{tikz}
\usepackage{tikz-cd}
\usetikzlibrary{arrows,backgrounds,patterns.meta}
\usetikzlibrary{positioning,shadings,cd}
\usetikzlibrary{shapes}
\usetikzlibrary{backgrounds}
\usetikzlibrary{decorations,decorations.pathreplacing,decorations.markings}
\usetikzlibrary{fit,calc,through}
\usetikzlibrary{external}
\usetikzlibrary{arrows}
\tikzset{vertex/.style = {shape=circle,draw,fill=black,inner sep=0pt,minimum size=5pt}}
\tikzset{edge/.style = {->,> = latex', bend right}}
\tikzset{
	super thick/.style={line width=3pt}
}
\tikzset{
    quadruple/.style args={[#1] in [#2] in [#3] in [#4]}{
        #1,preaction={preaction={preaction={draw,#4},draw,#3}, draw,#2}
    }
}
\tikzstyle{shaded}=[fill=red!10!blue!20!gray!30!white]
\tikzstyle{unshaded}=[fill=white]
\tikzstyle{empty box}=[circle, draw, thick, fill=white, opaque, inner sep=2mm]
\tikzstyle{annular}=[scale=.7, inner sep=1mm, baseline]
\tikzstyle{rectangular}=[scale=.75, inner sep=1mm, baseline=-.1cm]
\tikzstyle{mid>}=[decoration={markings, mark=at position 0.5 with {\arrow{>}}}, postaction={decorate}]
\tikzstyle{mid<}=[decoration={markings, mark=at position 0.5 with {\arrow{<}}}, postaction={decorate}]
\tikzstyle{over}=[double, draw=white, super thick, double=]

\tikzdeclarepattern{
  name=primeddots,
  type=uncolored,
  bounding box={(-.6pt,-.6pt) and (.6pt,.6pt)},
  tile size={(5pt,5pt)},
  tile transformation={rotate=60},
  parameters={none}, 
  code={
    \fill(0pt,0pt) circle (.35pt);
  }
}

\tikzdeclarepattern{
  name=primeddots2,
  type=uncolored,
  bounding box={(-.3pt,-.3pt) and (.3pt,.3pt)},
  tile size={(2.5pt,2.5pt)},
  tile transformation={rotate=60},
  parameters={none}, 
  code={
    \fill(0pt,0pt) circle (.35pt);
  }
}

\tikzstyle{primedregion}[none]=[
	preaction={fill=#1},
	pattern=primeddots,
  draw=#1,
]

\tikzstyle{primedregion2}[none]=[
	preaction={fill=#1},
	pattern=primeddots2,
  draw=#1,
]

\newcommand{\roundNbox}[6]{
	\draw[rounded corners=5pt, very thick, #1] ($#2+(-#3,-#3)+(-#4,0)$) rectangle ($#2+(#3,#3)+(#5,0)$);
	\coordinate (ZZa) at ($#2+(-#4,0)$);
	\coordinate (ZZb) at ($#2+(#5,0)$);
	\node at ($1/2*(ZZa)+1/2*(ZZb)$) {#6};
}

\newcommand{\tikzmath}[2][]
     {\vcenter{\hbox{\begin{tikzpicture}[#1]#2
                     \end{tikzpicture}}}
     }

\newcommand{\AColor}{gray!30}

\tikzcdset{scale cd/.style={every label/.append style={scale=#1},
    cells={nodes={scale=#1}}}}

\let\OLDthebibliography\thebibliography
\renewcommand\thebibliography[1]{
  \OLDthebibliography{#1}
  \setlength{\parskip}{0pt}
  \setlength{\itemsep}{0pt plus 0.3ex}
}

\begin{document}
\title{Quantum graphs and spin models}
\date{\today}
\author{
N\'estor Bravo Hern\'andez\footnote{CIMAT, Guanajuato, M\'exico \hfill \url{nestor.bravo@cimat.mx}}\hspace{.1cm}, Roberto Hern\'andez Palomares\footnote{Pure Mathematics, University of Waterloo \hfill \url{robertohp.math@gmail.com}}\hspace{.1cm} and Fabio Viales Sol\'is\footnote{Escuela de Matem\'atica, Universidad de Costa Rica \hfill \url{jose.vialessolis@ucr.ac.cr}}
}
\maketitle

\begin{abstract}
    We quantize the regularity properties of classical graphs that determine spin models for singly-generated Yang-Baxter planar algebras, including the Kauffman polynomial, and construct explicit examples. 
    A source of examples comes from deforming graphs using higher-idempotent splittings of quantum isomorphisms for which we prove that the relevant algebraic, combinatorial, and topological properties of the original graphs are preserved along with the quantum automorphism group. 
    We also obtain exotic examples of highly regular quantum graphs using the quantum Fourier transform and a method of iterated convolution. 
    Our examples include quantum versions of the strongly regular $9$-Paley, $16$-Clebsch and the Higman-Sims graphs, yielding new models for their regularity parameters. 
    As applications, we construct a compact quantum group that is monoidally equivalent to $SO_q(5)$ at the square of the golden ratio, whose dual is infinite with property (T), and exhibit a highly-regular quantum graph with no classical analogue. 
    Finally, we introduce quantum spin models, construct explicit examples and make contact with quantum Hadamard matrices. 
\end{abstract}

\section{Introduction}
Spin models are abstractions of physical systems with locally interacting components, appearing in statistical mechanics, and are known to yield knot and link invariants via state-sum constructions \cite{MR975085, MR990215, MR1247144}. 
For example, in \cite{MR899057}, Kauffman recovered the Jones polynomial \cite{MR766964}—originally defined using von Neumann algebras and traces of representations of braid groups.  
Shortly after, the Kauffman polynomial was introduced in \cite{MR958895}, generalizing the Jones polynomial. 
Jaeger later obtained new spin models for specializations of the Kauffman polynomial in \cite{MR1188082}, classifying them using certain association schemes and highly-regular \emph{classical} graphs such as the Higman-Sims graph associated to the simple sporadic group of order $44, 352, 000$ from \cite{MR227269}. 

This classification was further clarified by Edge in \cite{2019arXiv190208984E}, casting it in the framework of \emph{singly generated Yang-Baxter planar algebras}—a generalization of link diagrams, classified in \cite{MR1733737, MR1972635, MR3592517}—and providing a nice description of a graph's higher-regularity using skein relations. 
The conditions on a graph to determine a spin model for the Kauffman polynomial are extremely restrictive, and hence only a handful of examples are known: lattice graphs such as the square and the $9$-Paley graphs, as well as the pentagon, the $16$-Clebsch and the Higman-Sims graphs.
It is an open problem to decide whether there are any other such classical graphs, and herein we intended to quantize this data and produce new examples of highly-regular \emph{quantum} graphs. 

Quantum graphs are a vast non-commutative generalization of classical graphs and play important roles in many fields. 
They arise in quantum information as confusability graphs of noisy quantum channels in zero-error communication problems \cite[\S7]{MR3549479}, as models for completeness problems in quantum complexity theory \cite{Culf-Mehta}, and in the study of non-local games \cite{MR3849575,MR3907958, MR4091496, MR4505911, MR4507619, MR4752739, 2024arXiv240815444G}. 
Quantum graphs also appear as Cayley graphs of quantum groups \cite{MR2130588, MR4831258, 2025arXiv250311149B}, and were generalized by the second-named author and Brannan in \cite{doi:10.1142/S0129167X25500399} using subfactors and the diagrammatics of tensor categories. 
Certain well-studied quantum graphs known as biprojections—generalizing cliques—arise in the study of intermediate subfactors \cite{MR1262294, Bisch97}, and appear as well as examples of quantum relations in von Neumann algebras \cite{MR2908249, MR4302212}. 
In recent years, quantum graphs have also been studied in their own right by various authors from the perspective of dynamics, algebraic graph theory, quantum groups and operator systems \cite{MR4082994, MR4481115, MR4514486, MR4555986, MR4706978, Matsuda_2024, 2025arXiv250522519C}. 

\smallskip

The main theme of this manuscript is to develop an operator algebraic theory of quantum graphs, incorporating tools from quantum algebra, quantum topology, category theory, quantum information, quantum groups, subfactors, and graph theory. 
Herein, taking inspiration from the classification from \cite{MR1188082, 2019arXiv190208984E}, we introduce the relevant properties for simple quantum graphs such as its girth (in particular triangle-freeness), connected components, and crucially \emph{3-point regularity} \cite{MR1469634}.
We use the quantum graph Laplacian \cite{Matsuda_2024} to define the number of connected components of a quantum graph and to detect cycles, introducing the notion of a quantum forest. 

We hope that quantum graphs will become a new source of spin models and knot invariants, which classically are incredibly difficult to find.
To that effect we make special emphasis on constructing examples as explicitly as possible, aided by computer software developed by us, and investigating various methods of construction. 
In fact, we show that most known examples of classical graphs giving spin models linked to the Kauffman polynomial do not to exist in isolation.  
\begin{thmalpha}\label{thmalpha:qDeformations} {\bf [Theorem \ref{thm:qDeformations}]}
    For each of the classical graphs  $\cG\in\{9P,\ 16Cl,\ HS\}$, there exists a concrete quantum (non-classical) graph  $\cG_q\in\{9P_q,\ 16Cl_q,\ HS_q \}$ which is connected, $3$-point regular with the same parameters, Laplacian, minimal polynomial, and a concrete quantum graph isomorphism 
    $$
        \cG\cong_q\cG_q.
    $$
\end{thmalpha}

With the exception of the pentagon—which admits no deformation—all these graphs are realized as Cayley graphs (c.f. \cite{Heinze2001, MR1975767} and \cite[Lemma 16.3]{MR347649}) over groups admitting subgroups of \emph{central type}.
This observation together with tools arising from quantum information \cite{MR3907958}, and Theorem \ref{thmalpha:QuInvariants} which we shall explain below, afford us to deform their classical analogues into \emph{new} quantum graphs. 
That is, none of the $\cG_q$ are \emph{classically isomorphic} to their corresponding $\cG.$
This is a crucial point as all the classical graphs in Theorem \ref{thmalpha:qDeformations} are uniquely determined by their regularity parameters, which are invariants under quantum isomorphisms, as stated in Theorem \ref{thmalpha:QuInvariants}.
We shall in fact establish that our newly built quantum graphs do define \emph{quantum spin models} for the Kauffman polynomial in way we shall make rather explicit in Section \ref{sec:QuSpinModels}. 

The deformation method that we use from \cite{MR3907958}—herein referred to as \emph{bubbling}—applies to general finite quantum graphs and is realized entirely using finite-dimensional Hilbert spaces. 
Noteworthy, it is not known if other results involving quantum isomorphisms obtained from combinatorial tools such as \cite{MR4232075}, and \cite{MR4658469, MR4855314} can be implemented using finite-dimensional resources. 
In the context of non-local games, it is therefore plausible that the quantum isomorphisms from Theorem \ref{thmalpha:qDeformations} are implementable in a lab with finite resources, and so $HS$ can be realized in a rather \emph{tame} way.
In fact, that $\cG_q\cong_q\ \cG$ implies the \emph{Linking Algebra} between their quantum automorphism groups is nonzero, as this quantum isomorphism gives a finite-dimensional representation of this $*$-algebra \cite{MR4091496}.
In turn, addressing in the positive {\bf a question of Noah Snyder} asking whether there are (quantum) graphs with the same quantum automorphism group as $HS$.

Bubbling as a method to deform graphs stems from the Morita theory and higher idempotent splitting of the associated tensor categories (see \cite{MR3849575, MR3907958}, \cite{MR4419534} and references therein), and can be traced back to Jones' work on planar algebras \cite{MR4374438}. 
Given a quantum graph $\cG = (\ell^2(X), \hat{T})$, where $\ell^2(X)$ is the \emph{quantum set of vertices} (i.e. a $\delta$-form, see Equation (\ref{eqn:FrobAalgebra})) and $\hat{T}$ is its \emph{adjacency operator}, and a quantum automorphism $\bbX$ with the additional structure of a C* Frobenius algebra, one can \emph{split} $\bbX$ as  $\bbX\cong P\bullet P^*,$ where $P$ is some quantum isomorphism $P:\ell^2(X)\cong_q\ell^2(X_q)$ onto another quantum set $\ell^2(X_q),$ and $-\bullet-$ denotes the composition of quantum functions from \cite{MR3849575}.
One can use $P$ to bubble $\hat{T}$ into another adjacency operator $\hat{T}_q$ defined over $\ell^2(X_q)$ (Construction \ref{const:Bubbling}) yielding a new quantum graph $\cG_q= (\ell^2(X_q), \hat{T}_q)$ quantum isomorphic to the original $\cG.$
The point of bubbling is that it preserves desirable properties of $\cG$, including those from the classification of spin models for the Kauffman polynomial. 
\begin{thmalpha}\label{thmalpha:QuInvariants}{ \bf [Theorem \ref{thm:QuInvariants}]}
    Let $\cG = (\ell^2(X), \hat{T})$ be a quantum graph and $\mathbb{X}$ a C*-Frobenius algebra quantum isomorphism.
    That is, $\bbX \cong P\bullet P^*$ for some quantum isomorphism
    $P: H\otimes \ell^2(X_q) \to \ell^2(X)\otimes H,$ where $H$ is a finite-dimensional Hilbert space, and $\ell^2(X_q)$ is some quantum set.
    Then, the $P$-bubbling mapping  
    $\cG\mapsto \cG_q=(\ell^2(X_q), \hat{T}_q)$ preserves:
    \begin{enumerate}[label={\bf (B\arabic*)}]
        \item\label{alphaitem:AlgRels} {\bf Algebraic relations:} for any polynomial $f(z),$ we have 
        $f(\hat{T}_q)= (f(\hat{T}))_q.$
        \item\label{alphaitem:TopSpecFeatures} {\bf Topological and spectral features:} The Laplacian (Definition \ref{defn: laplacian}) and Degree Matrix (Definition \ref{defn: Degree matrix}) transform as 
        $L(\hat{T}_q) = L(\hat{T})_q,$ and $\mathsf{Deg}(\hat{T}_q) = \mathsf{Deg}(\hat{T})_q,$ preserving spectra and eigenspaces, respectively.  
        Consequently, the number of Edges (Equation (\ref{eqn:NumQuEdges})), number of connected components, and existence of cycles (Definition \ref{defn:pi0}) are invariants. 
        \item\label{alphaitem:HighReg} {\bf Higher regularity:} $t$-point regularity (Definition \ref{defn:3ptReg}) for $t\in\{1,2,3\}$, preserving the parameters. 
        Similarly, girth and triangle-freeness are quantum invariants.
        \item\label{alphaitem:BubblingQAUT} {\bf Quantum symmetries:} Quantum automorphism group up to monoidal equivalence.
    \end{enumerate}
\end{thmalpha}

We now present some {\bf applications concerning non-local games, quantum groups}, and make contact with {\bf algebraic graph theory}.
\begin{coralpha}
    Bubbling (quantum) graphs have the following consequences:
    \begin{enumerate}[label={\bf (C\arabic*)}]
        \item\label{item:Obstructions} The $2$-point and $3$-point regularity parameters of a quantum graph $\cG$ yield {\bf obstructions for $\cG$ to have classical representatives}. 
        That is, $\cG$ cannot be quantum isomorphic to a classical graph if any of these numbers is not a non-negative integer. Similarly, these properties with their parameters are {\bf obstructions to \emph{pseudotelepathy}} (c.f. \cite{MR3907958}).

        \item\label{item:(T)} Since $Aut^+(HS) \overset{\otimes}{\cong} SO_q(5)$, where $q$ is the square of the golden ratio \cite{MR1469634}, the monoidal equivalence  \ref{alphaitem:BubblingQAUT} yields a new {\bf infinite discrete quantum group with property $(T)$}. 
        Furthermore, the adjacency matrix of $HS_q$ gives a new fiber functor $\Rep(SO_q(5))\to \fdHilb.$

        \item\label{item:quTriangleFree} The $16Cl_q$ and $HS_q$ graphs are {\bf examples of connected triangle-free strongly regular quantum graphs} with the same parameters as their classical analogues (Example \ref{ex:params}). 

        \item\label{item:GirthAbelian} Quantum graphs with girth $\geq 5$ are not quantum isomorphic to Cayley graphs over abelian groups. 
    \end{enumerate}
\end{coralpha}
\noindent
Item \ref{item:quTriangleFree} provides {\bf an answer to the open problem in algebraic combinatorics} \cite[\S 2]{MR1312732} asking if there are only seven connected triangle-free strongly regular graphs. 
While the problem in the classical setting remains open to our knowledge, the answer is negative if one considers quantum graphs. 
Regarding Item \ref{item:Obstructions}, in contrast, Matsuda has announced that the degree $k$ of a $1$-point regular tracial quantum graph must be a natural number, which we observed in all of our examples. 
Item \ref{item:GirthAbelian} is a known obstruction in the classical setting and follows directly from Theorem \ref{thmalpha:QuInvariants}.
We briefly comment on Items \ref{item:(T)} and \ref{alphaitem:BubblingQAUT} in Remark \ref{remark:Prop(T)}.

{\bf We developed software} designed to verify if a given operator corresponds to a quantum graph on any given finite dimensional quantum set, compute the number of connected components, detect cycles, decide if it is $t$-point regular for $t=1,2,3$ and return the regularity parameters $(n,\ k,\ \lambda,\ \mu)$ for strong regularity (i.e. $2$-point regularity) and $q_0,\ q_1,\ q_2,\ q_3$ for $3$-point regularity. 
Our Software is available at the GitHub repository \href{https://github.com/nestorbravo/Regularity_for_quantum_graphs.git}{\texttt{Regularity\_for\_quantum\_graphs}}, where the interested reader can carry out most of our computations. 
Equipped with these tools, we constructed detailed examples in various dimensions, with {\bf applications to quantum topology and (quantum) association schemes}.
\begin{coralpha}\label{exalpha:M3TripleRegularity}
We constructed a number examples over various quantum spaces: 
    \begin{enumerate}[label={\bf (D\arabic*)}:]
        \item {\bf Examples over $M_2$:} The classification by Matsuda and Gromada of quantum graphs defined on $M_2$ from \cite{MR4481115, MR4514486} says that (up to complementation) the only graphs are the complete graph $\cK_{M_2}$ and a deformation of the square. 
        Aided with Theorem \ref{thmalpha:QuInvariants} and their classical representatives, we conclude these are connected and 3-point regular and obtain their parameters in Example \ref{ex:M2TripleRegularity}.
        \item {\bf Examples over $M_3$:} In example Example \ref{ex:M3TripleRegularity} we construct $4$ different $3$-point regular graphs over $M_3$ up to complementation $G_3, G_4, 9P_q$ and the irreflexive complete graph ${\cK}_{M_3}^{\mathsf{irref}},$ and obtained their parameters. 
        \item In Examples \ref{exs:OtherDeformations} we also deformed the {\bf Schl\"{a}fli Graph} and the {\bf Shrikhande graph}, modeled over non-commutative quantum sets of vertices.  
    \end{enumerate}
\end{coralpha}

The graph $G_3$ is constructed from the quantum Fourier transform \cite{Ocn90, Bisch97, MR4236188} —a notion closely related to the duality maps of association schemes (c.f. \cite{MR1414468})—of a rank $3$ projection, and $G_4$ is obtained as a limit of a geometric series of convolutions starting from appropriate initial conditions. 
We notice that all the graphs above satisfy condition $q_3-3q_2+3q_1-q_0\neq 0$ on the parameters of $3$-point regularity, arising from the classification of Edge \cite[Theorem 3]{2019arXiv190208984E} and the formal self-duality of the related association schemes investigated by Jaeger \cite{MR1358632}. 
Among these, the most interesting example is $G_4$ for several reasons.
First, as observed in Example \ref{ex:M3TripleRegularity}, both $G_4$ and its complement $G_6$ are connected \emph{genuinely quantum} graphs visibly from their parameters.  
The existence of $G_4$ then {\bf answers Gromada's open question} on whether there are any strongly regular quantum graphs with parameter $\mu\neq 0$ not quantum isomorphic to a classical graph \cite[page 17]{2024arXiv240406157G}. 
By inspection of its parameters, the graph $G_7,$ the complement of $G_3,$ appears to be quantum isomorphic to the complete tripartite graph $\cK_{3,3,3}.$ 

Motivated by having all these new examples of highly-regular quantum graphs, we introduce {\bf quantum spin models} in Section \ref{sec:QuSpinModels}. 
These are generalizations of their classical counterparts, where the set of possible states of a system's constituents is replaced by a quantum set. 
Using the association scheme tools developed by Jaeger in \cite{MR1188082}, following Construction \ref{const:BoltzmannWeights}, we write down explicit {\bf new examples of quantum association schemes} and their respective {\bf Boltzmann weights}—some with a classical counterpart—in Examples \ref{ex:SpinQSquare}, \ref{ex:SpinQPaley}, \ref{ex:SpinQClebsch}, \ref{ex:SpinQHS} and \ref{ex:SpimModG4}. 
In all these cases, we find these determine specializations of the Kauffman polynomial along with the relevant parameters.
In the \emph{topologically trivial} cases \ref{ex:SpinQClebsch} and \ref{ex:SpimModG4},  we obtain new examples of {\bf (real) quantum Hadamard matrices}. 
These claims can alternatively be checked computationally using our Software. 

\bigskip

{\bf Acknowledgements:} 
NBH and FVS were supported by {\bf Outsourcing Math Research}'s \emph{2024 Summer Research Program} and the \emph{Fall 2024 -  Fall 2025 Research and Learning Program}. \url{https://sites.google.com/view/outsourcingmathresearch/home}
NBH was also supported by CIMAT, \emph{Universidad de Guanajuato}, and the Workshop on Quantum Graphs at Saarland University. 
The authors benefited from fruitful conversations with Daniel Gromada, Michael Brannan, Junichiro Matsuda, Noah Snyder, Priyanga Ganesan, Ada Chan, Ofir Schnabel, and Kenny De Commer.
We are also grateful to the other participants of OMR 2024 Summer Research Program, Sergio Ching Heredia, Ivan Ortiz Arana, and Emmanuel Silva Amaya, who provided invaluable input during the initial stages of this project. 
RHP was partially supported by an NSERC Discovery Grant.

\tableofcontents

\section{Preliminaries}
    Let $A$ be a finite-dimensional $C^*$-algebra equipped with a faithful state $\psi$ and consider the GNS Construction with respect to the non-degenerate inner product 
    $$
    \langle a,b\rangle=\psi(a^*b) \qquad \forall a,b\in A.
    $$

    We will use the following diagrammatic representation of the multiplication and unit on $A$:
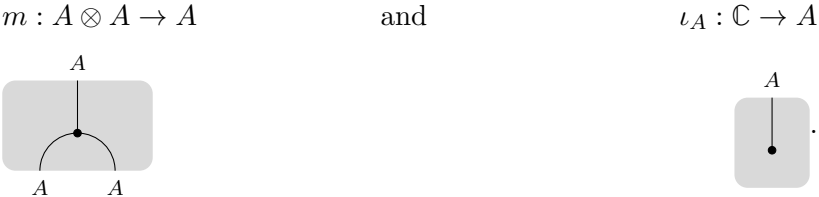
    \begin{align*}
    &m: A\otimes A\to A\qquad\qquad\qquad  \text{ and }\quad &\iota_A: \mathbb{C}\to A\\
    &\tikzmath{
        \begin{scope}
        \clip[rounded corners=5pt] (-1,0) rectangle (1,1.2);
        \fill[\AColor] (-1,0) rectangle (1,1.2);
        \end{scope}
        \draw (0,1.2) node[above]{$\scriptstyle A$};
        \draw (-.5,0) arc (180:0:.5cm);
        \draw (0,.5) -- (0,1.2);
        \filldraw (0,.5) circle (.05cm);
        \draw (-.5,0) node[below]{$\scriptstyle A$};
        \draw (0.5,0) node[below]{$\scriptstyle A$};
    }\qquad
    &\tikzmath{
        \begin{scope}
        \clip[rounded corners=5pt] (-.5,0) rectangle (.5,1.2);
        \fill[\AColor] (-.5,0) rectangle (.5,1.2);
        \end{scope}
        \draw (0,1.2) node[above]{$\scriptstyle A$};
        \draw (0,.5) -- (0,1.2);
        \filldraw (0,.5) circle (.05cm);
    }.
    \end{align*}
    
    Identifying $A$ with its GNS construction by a slight abuse of notation, we can take the adjoints of its multiplication and its unit map  obtaining bounded maps represented as 
    \begin{align*}
        &m^\dag:A\to A\otimes A \qquad\qquad\qquad \text{and} \qquad &\iota^\dagger: A\to \mathbb{C}\\
        &\tikzmath{
        \begin{scope}
        \clip[rounded corners=5pt] (-1,0) rectangle (1,-1.2);
        \fill[\AColor] (-1,0) rectangle (1,-1.2);
        \end{scope}
        \draw (0,-1.2) node[below]{$\scriptstyle A$};
        \draw (-.5,0) arc (180:360:.5cm);
        \draw (0,-.5) -- (0,-1.2);
        \filldraw (0,-.5) circle (.05cm);
        \draw (-.5,0) node[above]{$\scriptstyle A$};
        \draw (0.5,0) node[above]{$\scriptstyle A$};
    }\qquad
    &\tikzmath{
        \begin{scope}
        \clip[rounded corners=5pt] (-.5,0) rectangle (.5,-1.2);
        \fill[\AColor] (-.5,0) rectangle (.5,-1.2);
        \end{scope}
        \draw (0,-1.2) node[below]{$\scriptstyle A$};
        \draw (0,-.5) -- (0,-1.2);
        \filldraw (0,-.5) circle (.05cm);
    }.
    \end{align*}
    Then, the structure $X=(A, m, m^\dag, \iota, \iota^\dag)$ becomes a {\bf $\delta$-form}, where $\delta^2\in\mathbb{R}_{>0}.$ That is, $X$ is a {\bf Frobenius C*-algebra}; i.e. a (co)-associative (co)-unital (co)-algebra satisfying
    \begin{align}\label{eqn:FrobAalgebra}
        m\circ m^\dag = \delta^2\id_A \quad\text{ and }\quad m^\dag \circ m = (\id_A\otimes m)\circ (m^\dag\otimes \id_A) = (m\otimes \id_A)\circ(\id_A\otimes m^\dag). 
     \end{align}
     Diagrammatically, these identities are respectively expressed as
    \begin{align*}
    \tikzmath{
        \begin{scope}
        \clip[rounded corners=5pt] (-1.1,-1.2) rectangle (1.1,1.2);
        \fill[\AColor] (-1.1,-1.2) rectangle (1.1,1.2);
        \end{scope}
        \draw (-.5,0) arc (180:0:.5cm);
        \draw (0,.5) -- (0,1.2);
        \filldraw (0,.5) circle (.05cm);
        \draw (-.5,0) arc (180:360:.5cm);
        \draw (0,-.5) -- (0,-1.2);
        \filldraw (0,-.5) circle (.05cm);
    } 
    = \delta^2\cdot
    \tikzmath{
        \begin{scope}
        \clip[rounded corners=5pt] (-.5,-1.2) rectangle (.5,1.2);
        \fill[\AColor] (-1.1,-1.5) rectangle (1.1,1.5);
        \end{scope}
        \draw (0,-1.2) -- (0,1.2);
    },\quad\text{and }\quad
    \tikzmath{
        \begin{scope}
        \clip[rounded corners=5pt] (-1.1,-1.2) rectangle (1.1,1.2);
        \fill[\AColor] (-1.1,-1.2) rectangle (1.1,1.2);
        \end{scope}
        \draw (-.5,-1.2) arc (180:0:.5cm);
        \filldraw (0,.7) circle (.05cm);
        \draw (-.5,1.2) arc (180:360:.5cm);
        \draw (0,-.7) -- (0,.7);
        \filldraw (0,-.7) circle (.05cm);
    }
    =
    \tikzmath{
        \begin{scope}
        \clip[rounded corners=5pt] (-1.1,-1.2) rectangle (1.1,1.2);
        \fill[\AColor] (-1.1,-1.2) rectangle (1.1,1.2);
        \end{scope}
        \draw (-.6, -1.2) -- (-.6,.0) arc (180:0:.3) -- (-0,0) arc (180:360:.3) --(.6, 1.2);
        \draw (-.3, .3) -- (-.3, 1.2);
        \draw (.3, -.3) -- (.3, -1.2);
        \filldraw (.3, -.3) circle (.05cm);
        \filldraw (-.3, .3) circle (.05cm);
    }
    =
    \tikzmath{
        \begin{scope}
        \clip[rounded corners=5pt] (-1.1,-1.2) rectangle (1.1,1.2);
        \fill[\AColor] (-1.1,-1.2) rectangle (1.1,1.2);
        \end{scope}
        \draw (-.6, 1.2) -- (-.6, 0) arc (180:360:.3) -- (-0,0) arc (180:0:.3) --(.6, -1.2);
        \draw (-.3, -.3) -- (-.3, -1.2);
        \draw (.3, .3) -- (.3, 1.2);
        \filldraw (.3, .3) circle (.05cm);
        \filldraw (-.3, -.3) circle (.05cm);
    }.
    \end{align*}

    \begin{defn}\label{dfn:FinQuSet}
    A {\bf finite quantum set} $X$ is the $\delta$-form obtained from a finite dimensional C*-algebra equipped with a faithful state $(A, \psi)$.
    Whenever we want to refer to the underlying GNS Hilbert space from $(A, \psi)$ we denote it by $\ell^2(X)$. 
\end{defn}
    We will often obfuscate the C*-algebra $A$ and simply refer to it as $X$, emphasizing it's Frobenius algebra structure. Similarly, we write $\id_X$ instead of $\id_A,$ or sometimes simply $\id$ when no confusion may arise. 

\begin{ex}\label{ex:ClassicalDeltaForm}[{\bf Classical sets as $\delta$-forms}]
    We will interpret Definition \ref{dfn:FinQuSet} for classical sets. To this end, we start with $\mathbb{C}^N$ endowed with the normalized trace $\tr.$ Equipped with the orthonormal basis $\{u_a:= \sqrt{N}e_a\}_{a=1}^N$ we expand the comultiplication
    $$
    m^\dag(u_a) = \sqrt{N}u_a\otimes u_a,
    $$
    and so for any $a\in \{1,2,...,N\}$
    $$m\circ m^\dag(u_a) = \sqrt{N}^2u_a.
    $$
    And thus $\delta^2 = N.$
    Furthermore, $\iota^\dag (u_a) =\tr(u_a)=1/\sqrt{N}$,  where $[\iota]^a = 1/\sqrt{N}$ entrywise. 
\end{ex}

    From the relations above, it follows at once that any Frobenius C*-algebra $X$ is {\bf self-dual} with the following maps as {\bf solutions to the duality equations} (cf \cite[\S2]{2024arXiv240901951B})
    \begin{align*}
    \ev_X&:=
    \tikzmath{
        \begin{scope}
        \clip[rounded corners=5pt] (-1,0) rectangle (1,1.2);
        \fill[\AColor] (-1,0) rectangle (1,1.2);
        \end{scope}
        \draw (-.5,0) arc (180:0:.5cm);
    }
    :=
    \tikzmath{
        \begin{scope}
        \clip[rounded corners=5pt] (-1,0) rectangle (1,1.2);
        \fill[\AColor] (-1,0) rectangle (1,1.2);
        \end{scope}
        \draw (-.5,0) arc (180:0:.5cm);
        \draw (0,.5) -- (0,.8);
        \filldraw (0,.5) circle (.05cm);
        \filldraw (0,.8) circle (.05cm);
    }
    :X\otimes X\to \mathbb{C},\\ \\
    \coev_X&:=
    \tikzmath{
        \begin{scope}
        \clip[rounded corners=5pt] (-1,0) rectangle (1,-1.2);
        \fill[\AColor] (-1,0) rectangle (1,-1.2);
        \end{scope}
        \draw (-.5,0) arc (180:360:.5cm);
    }
    :=
    \tikzmath{
        \begin{scope}
        \clip[rounded corners=5pt] (-1,0) rectangle (1,-1.2);
        \fill[\AColor] (-1,0) rectangle (1,-1.2);
        \end{scope}
        \draw (-.5,0) arc (180:360:.5cm);
        \draw (0,-.5) -- (0,-.8);
        \filldraw (0,-.5) circle (.05cm);
        \filldraw (0,-.8) circle (.05cm);
    }
    : \mathbb{C}\to X\otimes X.
    \end{align*}
    Satisfying the {\bf Zig-Zag Equations}
    \begin{align*}
    (\ev_X\otimes \id)\circ(\id\otimes\coev_X)
    =
    \tikzmath{
        \begin{scope}
        \clip[rounded corners=5pt] (-1.1,-1.2) rectangle (1.1,1.2);
        \fill[\AColor] (-1.1,-1.2) rectangle (1.1,1.2);
        \end{scope}
        \draw (-.6, -1.2) -- (-.6,.0) arc (180:0:.3) -- (-0,0) arc (180:360:.3) --(.6, 1.2);
       }
    =
    \tikzmath{
        \begin{scope}
        \clip[rounded corners=5pt] (-.5,-1.2) rectangle (.5,1.2);
        \fill[\AColor] (-1.1,-1.5) rectangle (1.1,1.5);
        \end{scope}
        \draw (0,-1.2) -- (0,1.2);
    }
    =
    \tikzmath{
        \begin{scope}
        \clip[rounded corners=5pt] (-1.1,-1.2) rectangle (1.1,1.2);
        \fill[\AColor] (-1.1,-1.2) rectangle (1.1,1.2);
        \end{scope}
        \draw (-.6, 1.2) -- (-.6, 0) arc (180:360:.3) -- (-0,0) arc (180:0:.3) --(.6, -1.2);
    }
    =
    (\id\otimes\ev_X)\circ(\coev_X\otimes \id). 
    \end{align*}
    Notice that 
    $$
        \ev_X\circ \coev_X = \delta^2
    $$
    yields the {\bf quantum dimension of $X$}.

The endomorphisms of $\ell^2(X)$ have the structure of a C*-algebra with the usual composition of linear maps, the operator norm, and the usual adjoints. However, the structure of a Frobenius algebra/$\delta$-form endows $\End(\ell^2(X))$ with a second C*-algebra structure under a convolution product and a conjugate involution which we now introduce:
\begin{defn}
For a finite quantum set $X$, we define the {\bf Schur/convolution product} on $\End(\ell^2(X))\times \End(\ell^2(X))$ as follows:
\begin{equation}
S\star T:= \delta^{-2}\cdot m\circ(S\otimes T)\circ m^\dag=\  
\delta^{-2}\cdot\tikzmath{
\begin{scope}
\clip[rounded corners=5pt] (-1.1,-1.5) rectangle (1.1,1.5);
\fill[\AColor] (-1.1,-1.5) rectangle (1.1,1.5);
\end{scope}
\draw (-.5,.3) arc (180:0:.5cm);
\draw (0,.8) -- (0,1.5);
\filldraw (0,.8) circle (.05cm);
\roundNbox{fill=white}{(-.5,0)}{.3}{0}{0}{$S$}
\roundNbox{fill=white}{(.5,0)}{.3}{0}{0}{$T$}
\draw (-.5,-.3) arc (180:360:.5cm);
\draw (0,-.8) -- (0,-1.5);
\filldraw (0,-.8) circle (.05cm);
}.
\end{equation}

Combining the evaluation map with the Hilbert space adjoint $\dag,$ we define a new involution called the {\bf conjugation} on $\End(\ell^2(X))$ as follows:
\begin{equation}
T^*:=\  
\tikzmath{
\begin{scope}
\clip[rounded corners=5pt] (-1.1,-1.5) rectangle (1.1,1.5);
\fill[\AColor] (-1.1,-1.5) rectangle (1.1,1.5);
\end{scope}
\draw (-.7, -1.5) --(-.7,.3) arc (180:0:.3cm) ;
\roundNbox{fill=white}{(0,0)}{.3}{0}{0}{$T^\dag$}
\draw(0,-.3) arc (180:360:.3cm) -- (.6, 1.5);
}=
\tikzmath{
\begin{scope}
\clip[rounded corners=5pt] (-1.1,-1.5) rectangle (1.1,1.5);
\fill[\AColor] (-1.1,-1.5) rectangle (1.1,1.5);
\end{scope}
\draw (.7, -1.5) --(.7,.3) arc (0:180:.3cm) ;
\roundNbox{fill=white}{(0,0)}{.3}{0}{0}{$T^\dag$}
\draw (0,-.3) arc (360:180:.3cm) -- (-.6, 1.5);
}.
\end{equation}
\end{defn}
\noindent We remark that in the classical setting, $(T^*)^\dag = (T^\dag)^*$ recovers the usual matrix transpose. 

In both cases, when considering $(\End(\ell^2 X), \circ, \dag)$ and $(\End(\ell^2 X), \star, *)$ we get C*-algebras which are not necessarily isomorphic.

This second C*-algebra structure affords an elegant definition of a quantum graph (cf \cite[\S 3]{2024arXiv240901951B}):
\begin{defn}
    Let $X$ be a finite quantum set.
    A {\bf Schur idempotent} $\hat{T}:\ell^2(X)\to \ell^2(X)$ is a linear map satisfying 
    $$
    \hat{T}\star\hat{T} = \hat{T}.
    $$
    A {\bf quantum graph} on $X$ is a pair $\cG=(X,\ \hat{T})$, where $\hat{T}: \ell^2(X) \to \ell^2(X)$ is a  Schur/convolution idempotent. 
    We say that $\cG$ is {\bf reflexive} if there exists some $z\in Z(\End(X))^+\setminus\{0\}$, that is, a nonzero positive central endomorphism $z$, with
    $$
        \hat{T}\star \id_X = z\cdot\id_X = \id_X\star \hat{T}.
    $$
    In case 
    $$
        \hat{T}\star \id_X = 0 = \id_X\star \hat{T},
    $$
    we say $\cG$ is {\bf irreflexive}. Notice that if, say, $\cG$ is irreflexive, doing $\hat{T}+\id_X$ makes it reflexive. If $\cG$ is reflexive, then $\hat{T}-\id_X$ makes the graph irreflexive.

    We say $\cG $ {\bf is real} if 
    $$
    \hat{T}^* = \hat{T},
    $$ 
    and say $\cG$ is {\bf bi-directed/undirected} whenever 
    $$
    \hat{T}^\dag = \hat{T}.
    $$
\end{defn}

Notice that classical graphs corresponds to $X= (\bbC^N, \tr);$ ie a commutative finite dimensional C*-algebra equipped with its normalized trace, and the adjacency operators $\hat{T}$ are precisely the $\{0,1\}$-matrices with respect to the canonical basis.   

\begin{ex}
    Consider the identity map $\id_X: \ell^2(X)\to \ell^2(X)$. This trivially gives a quantum Schur idempotent since Equation (\ref{eqn:FrobAalgebra}) gives
    $$\id_X \star \id_X = \delta^{-2}\cdot m(\id_X\boxtimes \id_X)m^\dagger = \id_X.$$
    We call $\cT=(X,\ \id_X)$ the \textbf{trivial reflexive graph on $X$}. 

    Similarly, the map $\hat{J} = \delta^2\cdot \iota\circ \iota^\dagger: \ell^2(X)\to \ell^2(X)$ is also a Schur Idempotent, since 
    $$\hat{J} \star \hat{J} =  \delta^{-2}\cdot m((\delta^2\iota\circ \iota^\dagger)\boxtimes (\delta^2\iota\circ \iota^\dagger))m^\dagger = \delta^{2}\cdot\iota\circ \iota^\dagger = \hat{J}.$$
    The quantum graph $\cJ = (X, \hat{J})$ is named the \textbf{complete reflexive graph on $X$.} 
\end{ex}

\begin{sub-ex}
We shall interpret $\cT$ and $\cJ$ in the classical setting, following Example \ref{ex:ClassicalDeltaForm}. To this end, we take the $\delta$-form from the example and the respective orthonormal basis, then we have
\begin{align*}
    \hat{J}(u_a)&=N\iota\circ \iota^\dagger(u_a) \\
    &=\sqrt{N}\id_X \\ &=\sum_{a=1}^{n} u_a.
\end{align*}
\noindent This shows that, classically, $\hat{J}$ is the \emph{all ones matrix} and $\cJ$ represents the \emph{complete graph} of $X$. On the other hand, it's clear that $\id_X$ is the \emph{identity matrix}, and then $\cT$ represents the graph with only self loops.
\end{sub-ex}

Given a quantum graph $\cG=(X,\ \hat{T}),$ {\bf the complement of $\cG$} is the quantum graph
\begin{align}\label{eqn:GraphComplement}
    \cG^c:=(X,\ \hat{T}^c:=\hat{J} -  \delta^{-2}\cdot\id_X - \hat{T}).
\end{align}
The complement of an irreflexive graph as defined above is irreflexive, and similarly for reflexive graphs. Notice this convention is cosmetic as one can pass back and forth between reflexive and irreflexive by adding or substracting $\id_X$ to $\hat{T}^c$.

An equivalent way to describe a quantum graph $\cG = (X,\ \hat{T})$ is by its {\bf edge projector} $T$ which is given by 
\begin{align}\label{eqn:EdgeProj}
T:=
\tikzmath{
\begin{scope}
\clip[rounded corners=5pt] (-1.1,-1.5) rectangle (1.1,1.5);
\fill[\AColor] (-1.1,-1.5) rectangle (1.1,1.5);
\end{scope}
\draw (-.7, -1.5) --(-.7,.3) arc (180:0:.3cm);
\draw(-.4,.6) -- (-.4,1.5);
\filldraw (-.4,.6) circle (.05cm);
\draw(.3,-.6) -- (.3,-1.5);
\filldraw (.3,-.6) circle (.05cm);
\roundNbox{fill=white}{(0,0)}{.3}{0}{0}{$\hat T$}
\draw(0,-.3) arc (180:360:.3cm) -- (.6, 1.5);
}\in \End(X\otimes X).
\end{align}
The operator $T$ is readily seen to be a $\circ$-idempotent whenever $\hat{T}$ is a $\star$-idempotent. Whenever $\hat{T}$ is real (ie $\hat{T}= \hat{T}^*$) and undirected (ie $\hat{T}= \hat{T}^\dag$) we moreover have 
$$
T=
\tikzmath{
\begin{scope}
\clip[rounded corners=5pt] (-1.1,-1.5) rectangle (1.1,1.5);
\fill[\AColor] (-1.1,-1.5) rectangle (1.1,1.5);
\end{scope}
\draw (.7, -1.5) --(.7,.3) arc (0:180:.3cm) ;
\draw(.4,.6) -- (.4,1.5);
\filldraw (.4,.6) circle (.05cm);
\draw(-.3,-.6) -- (-.3,-1.5);
\filldraw (-.3,-.6) circle (.05cm);
\roundNbox{fill=white}{(0,0)}{.3}{0}{0}{$\hat{T}$}
\draw (0,-.3) arc (360:180:.3cm) -- (-.6, 1.5);
}.
$$
In this case, $T$ is readily seen to be a $\dag$-idempotent. 

The map $\hat{T}\mapsto T$ is closely related to the {\bf Quantum Fourier Transform}, whose role in the context of quantum graphs was explained in \cite{2024arXiv240901951B}. We will however not need this approach here. 

\bigskip

We will be concerned with classical and quantum isomorphisms of quantum graphs as introduced in \cite{MR3849575} and we recall here for the reader's convenience. 
We say that quantum graphs $\cG=(X,\hat{T})$ and $\cG'=(X',\hat{T}')$ are {\bf quantum-isomorphic} if and only if there exists a \emph{finite-dimensional} Hilbert space $H$ and a  linear map 
\begin{align}
    P:H\otimes X\to X'\otimes H
\end{align}
satisfying Equations (\ref{eqn:QuantumFunction}), (\ref{eqn:QuantumBijection}) and (\ref{eqn:QuantumGraphHom}):\\
Being a {\bf quantum function}
\begin{align}\label{eqn:QuantumFunction}
    \tikzmath{
    \begin{scope}
    \clip[rounded corners=5pt] (-.8,-1) rectangle (.6,1.8);
    \fill[\AColor] (-.8,-1) rectangle (.6,1.8);
    \end{scope}
    \filldraw (-.3,1.5) circle (.05cm);
    \draw (-.3,1.5) --(-.3,.3) .. controls ++(270:.2cm) and ++ (90:.2cm) .. (.3,-.3) --(.3,-1);
    \draw(-.6,1.8) arc(180:360:.3);
     \draw[thick] (-.3,-1) --(-.3,-.3) .. controls ++(90:.2cm) and ++ (270:.2cm) .. (.3,.3) --(.3,1.8);
    \roundNbox{fill=white}{(0,.2)}{.4}{0}{0}{$P$}
    \draw [->] (-.3,-.7);
    \draw [->] (.3,1.3);
    } 
    = 
    \tikzmath{
    \begin{scope}
    \clip[rounded corners=5pt] (-.6,-1.2) rectangle (1.3,1.8);
    \fill[\AColor] (-.6,-1.2) rectangle (1.3,1.8);
    \end{scope}
    \draw (-.3,1.8)  --(-.3,.3) .. controls ++(270:.2cm) and ++ (90:.2cm) .. (.3,-.3) --(.3,-.5);
    \draw(.3,-.5) arc(180:360:.3) -- (.9,1.2);
    \draw[thick] (.9, 1.2)-- (.9, 1.8);
    \draw(.6,-.8) -- (.6,-1.2);
    \draw[thick] (-.3,-1.2) --(-.3,-.3) .. controls ++(90:.2cm) and ++ (270:.2cm) .. (.3,.3) --(.4,.7);
    \draw (.4,.7)-- (.4, 1.8);
    \roundNbox{fill=white}{(0,0)}{.4}{0}{0}{$P$}
    \roundNbox{fill=white}{(0.7,1)}{.4}{0}{0}{$P$}
    \filldraw (.6,-.8) circle (.05cm);
    \draw [->] (-.3,-.7);
    \draw [->] (.35,.55);
    \draw [->] (.9,1.6);
    }, 
    \quad
    \tikzmath{
    \begin{scope}
    \clip[rounded corners=5pt] (-.6,-1) rectangle (.6,1.5);
    \fill[\AColor] (-.6,-1) rectangle (.6,1.5);
    \end{scope}
    \filldraw (-.3,1) circle (.05cm);
    \draw (-.3,1) --(-.3,.3) .. controls ++(270:.2cm) and ++ (90:.2cm) .. (.3,-.3) --(.3,-1);
    \draw[thick] (-.3,-1) --(-.3,-.3) .. controls ++(90:.2cm) and ++ (270:.2cm) .. (.3,.3) --(.3, 1.5);
    \roundNbox{fill=white}{(0,0)}{.4}{0}{0}{$P$}
    \draw [->] (-.3,-.7);
    \draw [->] (.3,1);
    }\ 
    =\ 
    \tikzmath{
    \begin{scope}
    \clip[rounded corners=5pt] (-.6,-1) rectangle (.8,1.5);
    \fill[\AColor] (-.6,-1) rectangle (.8,1.5);
    \end{scope}
    \draw[thick] (-.3,-1)  --(-.3,-.3) .. controls ++(90:.2cm) and ++ (270:.2cm) .. (.3,.3) --(.3,1.5);
    \draw (.2, -.3) -- (.2, -1);
    \filldraw (.2,-.3) circle (.05cm);
    \draw [->] (-.3,-.5);
    \draw [->] (.3,1);
    },
    \quad \text{and}\quad
    \tikzmath{
    \begin{scope}
    \clip[rounded corners=5pt] (-.6,-1) rectangle (.6,1.2);
    \fill[\AColor] (-.6,-1) rectangle (.6,1.2);
    \end{scope}
    \draw[thick] (-.3,1.2) --(-.3,.3) .. controls ++(270:.2cm) and ++ (90:.2cm) .. (.3,-.3) --(.3,-1);
    \draw (-.3,-1) --(-.3,-.3) .. controls ++(90:.2cm) and ++ (270:.2cm) .. (.3,.3) --(.3, 1.2);
    \roundNbox{fill=white}{(0,0)}{.4}{0}{0}{$P^{\dag}$}
    \draw [->] (.3,-.7);
    \draw [->] (-.3,1);
    } 
    = 
    \tikzmath{
    \begin{scope}
    \clip[rounded corners=5pt] (-1.2,-1.2) rectangle (1.2,1.2);
    \fill[\AColor] (-1.2,-1.2) rectangle (1.2,1.2);
    \end{scope}
    \draw (-.9,-1.2) -- (-.9,.5) arc(180:0:.3) -- (-.3,.5) --(-.3,.3) .. controls ++(270:.2cm) and ++ (90:.2cm) .. (.3,-.3) --(.3,-.5) arc(180:360:.3) -- (.9,1.2);
    \draw[thick] (-.3,-1.2) --(-.3,-.3) .. controls ++(90:.2cm) and ++ (270:.2cm) .. (.3,.3) --(.3, 1.2);
    \roundNbox{fill=white}{(0,0)}{.4}{0}{0}{$P$}
    \draw [->] (-.3,-.7);
    \draw [->] (.3,.7);
    }.
\end{align}
The thick wire corresponding to $H$ is marked with an arrow pointing upwards to differentiate it from the wires labeled by Frobenius C*-algebras. 
A quantum function can be though of a random function whose evaluation at an input obeys a probability distribution modeled over subspaces of $H$. 
Notice that if $H\cong \bbC$, the first two properties imply $P$ is a unital co-homomorphism in the classical sense.

A quantum function is moreover a {\bf quantum bijection} if it satisfies moreover 
\begin{align}\label{eqn:QuantumBijection}
    \tikzmath{
    \begin{scope}
    \clip[rounded corners=5pt] (-.6,-1.3) rectangle (.8,1.5);
    \fill[\AColor] (-.6,-1.3) rectangle (.8,1.5);
    \end{scope}
    \filldraw (.3,-1) circle (.05cm);
    \draw (-.3,1.5) --(-.3,.3) .. controls ++(270:.2cm) and ++ (90:.2cm) .. (.3,-.3) --(.3,-1);
    \draw(0,-1.3) arc(180:0:.3);
     \draw[thick] (-.3,-1.3) --(-.3,-.3) .. controls ++(90:.2cm) and ++ (270:.2cm) .. (.3,.3) --(.3,1.5);
    \roundNbox{fill=white}{(0,.2)}{.4}{0}{0}{$P$}
    \draw [->] (-.3,-.7);
    \draw [->] (.3,1);
    }
    =
    \tikzmath{
    \begin{scope}
    \clip[rounded corners=5pt] (-.7,-1) rectangle (1.4,2.2);
    \fill[\AColor] (-.7,-1) rectangle (1.4,2.2);
    \end{scope}
    \draw (-.3,1.5)  --(-.3,.3) .. controls ++(270:.2cm) and ++ (90:.2cm) .. (.3,-.3) --(.3,-1);
    \draw(-.3,1.5) arc(180:0:.35);
    \draw[thick] (.9, 1.2)-- (.9, 2.2);
    \draw(.95,.7) -- (.95,-1);
    \draw[thick] (-.3,-1) --(-.3,-.3) .. controls ++(90:.2cm) and ++ (270:.2cm) .. (.3,.3) --(.4,.7);
    \draw (.4,.7)-- (.4, 1.5);
    \roundNbox{fill=white}{(0,0)}{.4}{0}{0}{$P$}
    \roundNbox{fill=white}{(0.7,1)}{.4}{0}{0}{$P$}
    \filldraw (0.05,1.85) circle (.05cm);
    \draw(.05, 1.85) -- (.05, 2.2);
    \draw [->] (-.3,-.6);
    \draw [->] (.9,1.8);
    \draw [->] (.35,.55);
    }, 
    \qquad\text{and}\qquad
    \tikzmath{
    \begin{scope}
    \clip[rounded corners=5pt] (-.6,-1.5) rectangle (.6,1);
    \fill[\AColor] (-.6,-1.5) rectangle (.6,1);
    \end{scope}
    \filldraw (.3,-1) circle (.05cm);
    \draw (-.3,1) --(-.3,.3) .. controls ++(270:.2cm) and ++ (90:.2cm) .. (.3,-.3) --(.3,-1);
    \draw[thick] (-.3,-1.5) --(-.3,-.3) .. controls ++(90:.2cm) and ++ (270:.2cm) .. (.3,.3) --(.3, 1);
    \roundNbox{fill=white}{(0,0)}{.4}{0}{0}{$P$}
    \draw [->] (-.3,-.7);
    \draw [->] (.3,.7);
    }
    =
    \tikzmath{
    \begin{scope}
    \clip[rounded corners=5pt] (-.6,-1) rectangle (.6,1.5);
    \fill[\AColor] (-.6,-1) rectangle (.6,1.5);
    \end{scope}
    \draw[thick] (-.3,-1)  --(-.3,-.3) .. controls ++(90:.2cm) and ++ (270:.2cm) .. (.3,.3) --(.3,1.5);
    \draw (-.2, .3) -- (-.2, 1.5);
    \filldraw (-.2,.3) circle (.05cm);
    \draw [->] (-.3,-.7);
    \draw [->] (.3,.7);
    }.
\end{align}
It was shown in \cite[Theorem 4.8]{MR3849575} that a quantum function $P$ is a quantum bijection if and only if the underlying linear map is unitary.

If $\cG=(X, \hat{T})$ and $\cG' = (X', \hat{T'})$ are quantum graphs, the quantum function $P$ is a {\bf homomorphism of quantum graphs} if
\begin{align}\label{eqn:QuantumGraphHom}
    \tikzmath{
    \begin{scope}
    \clip[rounded corners=5pt] (-1,-2) rectangle (2,3.5);
    \fill[\AColor] (-1,-2) rectangle (2,3.5);
    \end{scope}
    \draw(-.6, 3.5) -- (-.6, 1.8) arc(180:360:.3);
    \draw(-.3,1.5) -- (-.3,.4);
    \filldraw (-.3,1.5) circle (.05cm);
    \filldraw (.5,-1.5) circle (.05cm);
    \roundNbox{fill=white}{(0.2,2.2)}{.4}{0}{0}{$\hat{T'}$}
    \draw (.4,2.6) -- (.4,2.8) arc (180:0:.3) --(1, 1.6);
    \filldraw (.7,3.1) circle (.05cm);
    \draw(.7, 3.1) -- (.7, 3.5);
    \filldraw (1.5,.1) circle (.05cm);
    \draw (1.5,.1) -- (1.5, .8);
    \draw[thick](1.5,1.2)--(1.5,3.5);
    \draw[->] (1.5, 3);
    \roundNbox{fill=white}{(0,0)}{.4}{0}{0}{$P$}
    \draw(.2,-.4) -- (.2,-1.2) arc(180:360:.3);
    \filldraw (.5,-1.5) circle (.05cm);
    \roundNbox{fill=white}{(1,-.8)}{.4}{0}{0}{$\hat{T}$}
    \draw (1.2,-.4) -- (1.2,-.2) arc (180:0:.3) --(1.8,-2);
    \filldraw (1.5,.1) circle (.05cm);
    \draw (1.5,.1) -- (1.5, .8);
    \roundNbox{fill=white}{(1.2,1.2)}{.4}{0}{0}{$P$}
    \draw[thick] (0.3,0.4) .. controls (.8,.5).. (.9,.8);
    \draw[thick] (-.3,-.4) -- (-.3,-2);
    \draw[->] (-.3, -1);
    \draw[->] (.85, .7);
    \draw (.5,-1.5) -- (.5,-2);
    }\ 
    =\ 
    \tikzmath{
    \begin{scope}
    \clip[rounded corners=5pt] (-.8,-2) rectangle (2,3.5);
    \fill[\AColor] (-.8,-2) rectangle (2,3.5);
    \end{scope}
    \roundNbox{fill=white}{(0,0)}{.4}{0}{0}{$P$}
    \draw(.2,-.4) -- (.2,-1.2) arc(180:360:.3);
    \filldraw (.5,-1.5) circle (.05cm);
    \roundNbox{fill=white}{(1,-.8)}{.4}{0}{0}{$\hat{T}$}
    \draw (1.2,-.4) -- (1.2,-.2) arc (180:0:.3) --(1.8,-2);
    \filldraw (1.5,.1) circle (.05cm);
    \draw (1.5,.1) -- (1.5, .8);
    \roundNbox{fill=white}{(1.2,1.2)}{.4}{0}{0}{$P$}
    \draw[ thick] (0.3,0.4) .. controls (.8,.5) .. (.9,.8);
    \draw[thick] (-.3,-.4) -- (-.3,-2);
    \draw (.5,-1.5) -- (.5,-2);
    \draw[thick](1.5,1.6)--(1.5,3.5);
    \draw[->] (1.5, 3);
    \draw(1,1.6)--(1,3.5);
    \draw(-.3,.4)--(-.3,3.5);
    \draw[->] (-.3, -1);
    \draw[->] (.85, .7);
    }.
\end{align}
In case $H\cong \bbC$, the above condition reduces to $\cG$ and $\cG'$ being {\bf classically isomorphic}. 

Quantum functions can be composed in a way that generalizes the usual composition of functions. 
Given quantum sets $X,  X', X'',$ finite dimensional Hilbert spaces $H, K$ and quantum functions $P:H\otimes X\to X'\otimes H$ and $Q:K\otimes X'\to X''\otimes K$ the {\bf composition law for quantum function} is defined as
\begin{align}\label{eqn:QuFunComp}
    Q\bullet P:=
    \tikzmath{
    \begin{scope}
    \clip[rounded corners=5pt] (-2,-2) rectangle (2,2);
    \fill[\AColor] (-2,-2) rectangle (2,2);
    \end{scope}
    \draw[thick] (-.8, -2) -- (-.8,1.1);
    \draw (-.8,1.1) -- (-.8,2);
    \draw (.8,-.7) -- (.8,-2);
    \draw[thick] (.8, -.7) -- (.8,2);
    \draw (.2, -.4) .. controls (-.2,0) .. (-.2,1.1);
    \draw[thick] (-.2,1.1) .. controls (.5,1.5) .. (.6,2);
    \draw[thick] (-.6,-2) .. controls (-.5,-1.5) .. (.2,-1.1);
    \draw [->] (-.8,-1);
    \draw [->] (.8,1);
    \draw [->] (.51,1.7);
    \draw [->] (-.51,-1.6);
    \roundNbox{fill=white}{(.5,-.7)}{.4}{0}{0}{$P$}
    \roundNbox{fill=white}{(-.5,.7)}{.4}{0}{0}{$Q$}
    \draw (-.8,-2) node[below]{${\scriptstyle K\otimes H}$};
    \draw (.8,2) node[above]{${\scriptstyle K\otimes H}$};
    \draw (-.8,2) node[above]{${\scriptstyle X''}$};
    \draw (-.28, 0) node[below]{${\scriptstyle X'}$};
    \draw (.8,-2) node[below]{${\scriptstyle X}$};
    }
\end{align}
which can be seen to define a quantum function $Q\bullet P: (K\otimes H)\otimes X\to X''\otimes(K\otimes H).$

Equivalently, a quantum graph isomorphism is a linear map $P: H\otimes X\to X'\otimes H$ satisfying Equations (\ref{eqn:QuantumFunction}), (\ref{eqn:QuantumBijection}) and also intertwining the adjacency matrices
$$
P\bullet \hat{T} = \hat{T'}\bullet P. 
$$

Given a quantum graph $\cG =(X, \hat{T})$ we denote {\bf the category of quantum graph automorphisms} by $\mathsf{QAut}(\cG)$ whose objects are the quantum automorphisms of $\cG$ and the morphisms are the \emph{intertwiners of quantum functions}. 
We say that a linear map $f: H \to K$ is an {\bf intertwiner of quantum functions} $P:H\otimes X \to X'\otimes H$ and $Q:K\otimes X\to X'\otimes K$ if and only if
\begin{align}\label{eqn:QuFunInter}
    Q\circ( f\otimes \id_X) = (\id_{X'}\otimes f)\circ P.
\end{align}
Here, we used $\circ$ to denote the usual composition of linear maps, which is in general different from the composition of quantum functions $\bullet$ defined in Equation (\ref{eqn:QuFunComp}). 
Composition of quantum functions $\bullet$ and the usual tensor product of intertwiners makes $\QAut(\cG)$ into a tensor category. 
The category $\mathsf{QAut}(\cG)$ was originally introduced in \cite{MR3849575} where the interested reader can find a more detailed exposition.

\smallskip

The category $\mathsf{QAut}(\cG)$ can be rephrased in terms of the quantum automorphism group of a graph, which we shall briefly describe. 
Recall that a {\bf compact quantum group} $\bbG$ consists of a pair $\bbG = (C(\bbG),\ \Delta),$ where $C(\bbG)$ is a unital C*-algebra, and $\Delta: C(\bbG)\to C(\bbG)\otimes C(\bbG)$ is a unital $*$-homomorphism satisfying the usual co-associativity (i.e. $(\id_{C(\bbG)}\otimes \Delta)\circ\Delta = (\Delta\otimes id_{C(\bbG)})\circ \Delta$), and cancellation properties (i.e. $(C(\bbG)\otimes 1)\Delta[C(\bbG)]$ and $(1\otimes C(\bbG))\Delta[C(\bbG)]$ are dense in $C(\bbG)\otimes C(\bbG)$). 
For an introductory text on compact quantum groups, we refer the reader to \cite{MR3204665}. 
 
Given a quantum graph $\cG = (X, \hat{T}),$ the {\bf quantum automorphism group of $\cG$} is a compact quantum group 
$$
\bbG = (C(\Aut^+(\cG)),\ \Delta)
,$$
where $C(\Aut^+(\cG))$ is the universal C*-algebra generated by the coefficients $u_i^k$ of a unitary $u = [u_i^k]\in M_n(C(\Aut^+(\cG)))$ making the operator 
$$
\rho: X\ni e_i\mapsto\sum_k e_k\otimes u_i^k\in X\otimes C(\Aut^+(\cG))
$$
a unital $*$-homomorphism satisfying $ \rho\circ \hat{T} = (\hat{T}\otimes \id_{C(\Aut^+(\cG))})\circ \rho.$
Here,  $\{e_i\}_1^n\subset \ell^2(X)$ denotes an arbitrary orthonormal basis.  
The comultiplication is given by 
$$
\Delta(u_i^k) = \sum_j u_j^k\otimes u_i^j.
$$
We shall not discuss the other ingredients making $\Aut^+(\cG)$ into a Hopf C*-algebra as we will not use them. 
Further details can be found in \cite[\S2]{MR4481115}, and \cite[\S3]{MR4091496}. 

We now recall the relationship between $\QAut(\cG)$ and $\Aut^+(\cG)$ (c.f. \cite[Proposition 5.19]{MR3849575} and \cite[Theorem 2.43]{MR4481115}). 
As tensor categories, we have 
\begin{align}\label{eqn:QAutRep}
    \mathsf{QAut}(\cG) \cong \Rep(C(\Aut^+(\cG))),
\end{align}
where we only consider the finite dimensional objects. 

\smallskip

The structure of $\QAut(\cG)$ is in fact that of a \emph{rigid} C*-tensor category. 
As shown in \cite[Theorem 5.16]{MR3849575}, a quantum graph homomorphism in $\QAut(\cG)$ is \emph{dualizable} if and only if it is a quantum graph isomorphism. 
By definition, a map $P: H\otimes X\to X\otimes H$ in $\QAut(\cG)$ is {\bf dualizable} if and only if, there exists a $Q: \overline{H}\otimes X\to X\otimes \overline{H}$ in $\QAut(\cG)$, where $H$ is some finite dimensional Hilbert space, such that 
\begin{align}\label{eqn:dualityQAut}
    \coev_H\in\QAut(\cG)(\id \to P\bullet Q)\quad \text{ and }\quad \ev_H\in\QAut(\cG)(Q\bullet P \to \id),
\end{align}
where $\coev_H:\bbC\to H\otimes\overline{H}$ and $\ev_H:\overline{H}\otimes H\to \bbC$ denote the usual coevaluation and evaluation maps given by $1\mapsto \sum_{i}e_i\otimes \overline{e_i}$, summing over an orthonormal basis of $H$, and $\overline{\eta}\otimes \xi\mapsto \langle \eta\ |\ \xi\rangle^H$ is the inner product in $H$. 
As shown in \cite[Theorem 4.8]{MR3849575}, $P$ is dualizable in $\QAut(\cG)$ if and only if $P$ is a unitary map. 
We then say that $Q$ is the dual of $P$, and can be identified with its conjugate $P^*,$ corresponding to the rotation of $P$ by $180$ degrees. 

\begin{ex}
    Given quantum graphs $\cG = (X,\ \hat{T})$ and $\cG' = (X',\ \hat{T}')$ and a quantum graph isomorphism $Q:K\otimes \ell^2(X)\to \ell^2(X')\otimes K$ we obtain the quantum graph automorphisms $$
     Q\bullet Q^*\ \ \ \text{   and   }\ \ \ Q^*\bullet Q,
    $$
    which on the nose have structures of {\bf C*-Frobenius algebras in  $\mathsf{QAut}(\cG).$ }
    By Equation (\ref{eqn:dualityQAut}) and the discussion around it, the solutions to the duality equations $\ev_{K}: \overline{K}\otimes K\to \bbC$ and $\coev_K:\bbC\to K\otimes \overline{K}$ automatically are intertwiners in $\mathsf{QAut}(\cG)$:
    \begin{align}\label{eqn:DualityQAutG}
        \ev_K\in\mathsf{QAut}(\cG)(Q^*\bullet Q\to\id_X)\quad\text{ and }\quad \coev_K \in\mathsf{QAut}(\cG)(\id_X\to Q\bullet Q^*),   
    \end{align}
    which in diagrams becomes: 
    \begin{align}\label{eqn:DiagDualityQAut}
    \tikzmath{
    \begin{scope}
    \clip[rounded corners=5pt] (-1,-1) rectangle (1,1);
    \fill[\AColor] (-1,-1) rectangle (1,1);
    \end{scope}
    \draw (.9,-1) .. controls (.6,-.2) .. (0,0);
    \draw (-.9,1) .. controls (-.6,.2) .. (0,0);
    \draw[thick](0,-1) -- (0,-.7) arc(0:180:.3) -- (-.6,-1);
    \draw [->] (0,-.7);
    \draw [>-] (-.6,-.7);
    }
    =
    \tikzmath{
    \begin{scope}
    \clip[rounded corners=5pt] (-1,-1) rectangle (1,1);
    \fill[\AColor] (-1,-1) rectangle (1,1);
    \end{scope}
    \draw (.9,-1) .. controls (.6,-.2) .. (0,0);
    \draw (-.9,1) .. controls (-.6,.2) .. (0,0);
    \draw[thick] (.55,0) --(.5,.2)  arc(10:150:.5);
    \draw[thick] (.4,-.4) .. controls (.2,-.7) .. (.2,-1);
    \draw[->] (.2,-.8); 
    \draw[thick] (-.6,0) ..controls (-.7,-.7) .. (-.7,-1);
    \draw[>-] (-.7,-.7);
    \draw [->] (.5,.25);
    \roundNbox{fill=white}{(.5,-.2)}{.25}{0}{0}{${\scriptstyle Q}$}
    \roundNbox{fill=white}{(-.5,.2)}{.25}{0}{0}{${\scriptstyle Q^*}$}
    }
    \quad
    \text{ and }
    \quad
    \tikzmath{
    \begin{scope}
    \clip[rounded corners=5pt] (-1,-1) rectangle (1,1);
    \fill[\AColor] (-1,-1) rectangle (1,1);
    \end{scope}
    \draw (-.9,1) .. controls (-.6,.2) .. (0,0);
    \draw (.9,-1) .. controls (.6,-.2) .. (0,0);
    \draw[thick](0,1) -- (0,.7) arc(180:360:.3) -- (.6,1);
    \draw [->] (0,.7);
    \draw [>-] (.6,.7);
    }
    =
    \tikzmath{
    \begin{scope}
    \clip[rounded corners=5pt] (-1,-1) rectangle (1,1);
    \fill[\AColor] (-1,-1) rectangle (1,1);
    \end{scope}
    \draw (-.9,1) .. controls (-.6,.2) .. (0,0);
    \draw (.9,-1) .. controls (.6,-.2) .. (0,0);
    \draw[thick] (-.55,0) --(-.5,-.2)  arc(190:330:.5);
    \draw[thick] (-.4,.4) .. controls (-.2,.7) .. (-.2,1);
    \draw[->] (-.2,.8); 
    \draw[thick] (.6,0) ..controls (.7,.7) .. (.7,1);
    \draw[>-] (.7,.7);
    \draw [->] (-.5,-.15);
    \roundNbox{fill=white}{(-.5,.2)}{.25}{0}{0}{${\scriptstyle Q}$}
    \roundNbox{fill=white}{(.5,.-.2)}{.25}{0}{0}{${\scriptstyle Q^*}$}
    }
    \end{align}    
\end{ex}

\section{Laplacians and connectivity for quantum graphs}
Similarly as in \cite{Matsuda_2024}, we can define a quantum analog notion of graph Laplacian, which in the classical case provides an algebraic tool to compute the topological properties of graphs such as the number of connected components or to detect the existence of cycles. 

\begin{defn}\label{defn: laplacian}
    Given an irreflexive, undirected (i.e. real) quantum graph $\cG = (X,\ \hat{T})$, we define its {\bf quantum Laplacian operator} as 
    $$
        L(\cG):= (\iota^\dag\otimes \id_X)\circ(\hat{T}\otimes \id)\circ m^\dag - \hat{T}.
    $$
    Diagrammatically, we represent this as 
   $$ L(\hat{T}) = 
    \tikzmath{
        \begin{scope}
        \clip[rounded corners=5pt] (-.8,-1.5) rectangle (1.3,1);
        \fill[\AColor] (-.8,-2.5) rectangle (1.3,1);
        \end{scope}
        \draw (-.25,0.0) -- (-.25,.3);
        \filldraw (-.25,0.3) circle (.05cm);
        \draw (.75,-0.55) -- (.75,1);
        \roundNbox{fill=white}{(-.3,-.25)}{.3}{0}{0}{$\hat{T}$}
        \draw (-.25,-.55) arc (180:360:.5cm);
        \draw (.25,-1) -- (.25,-1.5);
        \filldraw (.25,-1.05) circle (.05cm);
        }
        -
        \tikzmath{
        \begin{scope}
        \clip[rounded corners=5pt] (1.3,-1.5) rectangle (2.6,1);
        \fill[\AColor] (1.3,-2.5) rectangle (2.6,1);
        \end{scope}
        \draw (2.0,-1.5) -- (2.0, 1);
        \roundNbox{fill=white}{(2.0,-.25)}{.3}{0}{0}{$\hat{T}$}
        }.
        $$
    \end{defn}

    \begin{defn}\label{defn: Degree matrix}
    We define {\bf the degree matrix of $\cG$} by
    $$
    \mathsf{Deg}(\hat{T}) = L(\hat{T}) + \hat{T} = 
    \tikzmath{
        \begin{scope}
        \clip[rounded corners=5pt] (-.8,-1.5) rectangle (1.2,1);
        \fill[\AColor] (-.8,-2.5) rectangle (2.5,1);
        \end{scope}
        \draw (-.25,0.0) -- (-.25,.3);
        \filldraw (-.25,0.3) circle (.05cm);
        \draw (.75,-0.55) -- (.75,1);
        \roundNbox{fill=white}{(-.3,-.25)}{.3}{0}{0}{$\hat{T}$}
        \draw (-.25,-.55) arc (180:360:.5cm);
        \draw (.25,-1) -- (.25,-1.5);
        \filldraw (.25,-1.05) circle (.05cm);
        }.
    $$
    \end{defn}

     \begin{defn}\label{defn: incidence matrix}
         Similarly, we define the {\bf incidence/gradient matrix} as (some positive real multiple of)
    $$
    D(\hat{T}) = \left(\hat{T} \otimes \id - \id \otimes \hat{T}\right)\circ m^\dag.
    $$ 
    Diagrammatically this corresponds to
    $$D(\hat{T}) = 
    \tikzmath{
        \begin{scope}
        \clip[rounded corners=5pt] (-.8,-1.5) rectangle (1.00,0.5);
        \fill[\AColor] (-.8,-2.5) rectangle (1.05,.5);
        \end{scope}
        \draw (-.25,0.0) -- (-.25,0.5);
        \draw (.75,-0.55) -- (.75,.5);
        \roundNbox{fill=white}{(-.3,-.25)}{.3}{0}{0}{$\hat{T}$}
        \draw (-.25,-.55) arc (180:360:.5cm);
        \draw (.25,-1) -- (.25,-1.5);
        \filldraw (.25,-1.05) circle (.05cm);
        
        \node at (1.25, -0.5) {$-$};

        \begin{scope}
        \clip[rounded corners=5pt] (1.45,-1.5) rectangle (3.22,0.5);
        \fill[\AColor] (1.48,-2.5) rectangle (3.25,.5);
        \end{scope}
        \draw (1.75,-.55) -- (1.75,0.5);
        \draw (2.75,-0.55) -- (2.75,.5);
        \roundNbox{fill=white}{(2.75,-.25)}{.3}{0}{0}{$\hat{T}$}
        \draw (1.75,-.55) arc (180:360:.5cm);
        \draw (2.25,-1) -- (2.25,-1.5);
        \filldraw (2.25,-1.05) circle (.05cm);
        }.
        $$
     \end{defn}
    
The {incidence matrix} of a \textit{classical directed graph} has an important interpretation, since it encodes the vertex-edge relationships. Each column represents an edge \( e_j \), and each row corresponds to a vertex \( v_i \). The entry \( b_{i,j} \) is \(-1\) if \( v_i \) is the \textit{tail} (source) of \( e_j \), \(+1\) if \( v_i \) is the \textit{head} (target), and \( 0 \) otherwise. This representation captures the graph's flow structure, where each column sums to zero (reflecting edge directionality: \(-1\) and \(+1\) cancel out). The matrix is fundamental in network analysis, Kirchhoff's circuit laws, and optimization problems \cite{Chung1996-vw}. 

The \textbf{Gradient} and \textbf{Divergence} operators are generalized using the incidence matrix \( D \). The \textit{discrete gradient} \(\nabla f = D^t f\) maps vertex potentials \( f \) to edge differences (analogous to directional derivatives), while the \textit{discrete divergence} \(\mathsf{div} g = D g\) computes net flow at vertices for edge flows \( g \). These adjoint operators satisfy \(\mathsf{div} (\nabla f) = L f\), where \( L = DD^t \) is the graph Laplacian. This notion extends to the quantum case as we show in the next statement. 

\begin{prop}
    Incidence matrix $D(\cG)$ is a (positive multiple of a) square root for $L(\cG)$, that is, 
    $$
    L(\cG) = C\cdot D(\cG)^\dagger\circ D(\cG)\geq 0,
    $$
    for some $C\in(0,\infty).$
\end{prop}

The proof is a particular case of that of \cite[Lemma 2.6]{Matsuda_2024}. 
\bigskip

The \textbf{graph Laplacian} in the classical case, encodes certain topological properties of a graph. Crucially, the \textit{multiplicity of the zero eigenvalue} of $L$ equals the number of \textbf{connected components}. This follows because each component contributes one dimension to the nullspace of $L$ via a constant eigenvector. The Laplacian's spectrum is thus fundamental for analyzing graph connectivity. Now, we extend these notions to the non-commutative case.

\begin{defn}\label{defn:pi0}
We define {\bf the number of connected components of an undirected quantum graph $\cG$} as 
    $$
    \#\pi_0(\cG) := \mathsf{null}(L)= \dim_\bbC(X)-\mathsf{rank}(L).
    $$
If $\cG$ is an irreflexive (ie $\hat{T}\star \id_Q=0$) quantum graph, we say {\bf $\cG$ has a cycle} if and only if  
    $$
    \Tr(L) > 2\cdot\mathsf{rank}(L).
    $$
If $$\Tr(L) = 2\cdot\mathsf{rank}(L),$$
we say {\bf $\cG$ is a forest}. 
Finally, {\bf $\cG$ is a tree} if and only if $\#\pi_0(\cG)=1$ and $\cG$ is a forest. 
\end{defn}
The previous definition is motivated by undirected simple irreflexive classical graphs, where the nullity of the Laplacian counts the number of connected components. 
For such graphs, one always has that 
$$\#E\geq \sum_{c\in \pi_0}(\#c-1) = \#V-\#\pi_0,$$
where each summand above is the number of edges in a minimal spanning tree for the connected component $c\in \pi_0.$ 
It is then obvious that a graph is a forest if and only if  equality is attained, and it has a cycle otherwise.\footnote{See \url{https://stackoverflow.com/questions/16436165/detecting-cycles-in-an-adjacency-matrix} for an interesting discussion.}
Furthermore, the left-hand-side of the inequality can be re-expressed as $\Tr(L)/2,$ corresponding to {\bf the number of edges}, since the trace is the sum over the degrees of all vertices.
That is, 
\begin{align}\label{eqn:NumQuEdges}
    \#E(\cG)=\Tr(L(\hat{T})/2.
\end{align}

\begin{remark}
    As shown in  \cite[Theorem 3.7]{Matsuda_2024} a regular quantum graph of degree $k$ is connected if and only if $k\in \mathsf{Spec}(\hat{T})$ is simple root, that is 
    $$1 = \dim(\ker (k\cdot \id_X - \hat{T}) )= \mathsf{null}(\mathsf{Deg}(\hat{T}) - \hat{T}) = \mathsf{null} (L\hat{T}) = \# \pi_0(\mathcal{G}). $$

    During the latter development stages of this manuscript, we learned about \cite[Proposition 4.4]{2025arXiv250522519C}, where a characterization of connectivity for (non-tracial) quantum graphs similar to ours is given in terms of the nullity of the Laplacian. Nevertheless, our definition for $\#\pi_0(\cG)$ is, to the best of our knowledge, new in the literature.
\end{remark}

\begin{ex}[{\bf Connectivity for trivial and complete quantum graphs}]
    Let $\ell^2(X)$ be a $\delta$-form/finite quantum set. 
    From a simple diagrammatic computation following Definitions \ref{defn: Degree matrix} and \ref{defn: laplacian}, one can show that the trivial graph $\cT$ on $Q=\ell^2(X)$ has Laplacian $L({\cT}) =0$, with $\#\pi_0(\cT) = \dim_\bbC(X) $ and $\mathsf{rank}(L(\cT))=0$ so the trivial graphs are always maximally disconnected forests. 

    Similarly, the Laplacian for the complete quantum graph $\cJ$ on $X$ is given by 
    $$
    L(\cJ) = \delta^{2}\cdot(\id_{X}-\iota\circ\iota^\dag),
    $$
    which is (a scalar multiple of) the projection onto $(\iota[\bbC])^\perp\subset \ell^2(X),$ and hence $\mathsf{rank}(L(\cJ)) = \dim_\bbC(X)-1.$ 
    Thus, $\#\pi_0(\cJ) = \mathsf{null}(L(\cJ)) = 1$, and $\cJ$ is connected. 
    Furthermore, 
    $$
    \dfrac{\Tr(L(\cJ))}{2} = \dfrac{\delta^{2}}{2}(\dim_\bbC(X)-1)\geq \dim_\bbC(X)-1 = \mathsf{rank}(L(\cJ)).
    $$
    Therefore, the complete quantum graph on any finite quantum set is always connected and it is a tree if and only if $\delta^2=2.$ 
    In particular, for classical sets $\bbC\overset{\tr}{\subset}{\bbC^N}$ this happens exactly when $N=2$, and for the tracial quantum sets $\bbC\overset{\tr}{\subset}{M_N}$ the complete graph is a tree if and only if $N=2.$ 
\end{ex}

\begin{ex}\label{ex:LaplaciansM2}
    Using the adjacency matrices for the four undirected reflexive quantum graphs on $(M_2, \tr)$ of degrees $\mathsf{deg} = 1,2,3,4$ from Matsuda's Schur idempotents and Definition \ref{defn:pi0} we obtain that 
    \begin{align*}
        L_{A_1} = 0, \quad &\sigma(L_{A_1}) = \{ 0,0,0,0 \}, \quad \mathsf{deg}(A_1) = 1,\\
        L_{A_2} = \begin{pmatrix}
            0 & 0 & 0 &0\\
            0 & 2 & 0 &0\\
            0 & 0 & 2 &0\\
            0 & 0 & 0 &0\\
        \end{pmatrix}, \quad &\sigma(L_{A_2}) = \{ 0,0,2,2 \}, \quad \mathsf{deg}(A_2) = 2,\\
        L_{A_3} = \begin{pmatrix}
            2 & 0 & 0 &-2\\
            0 & 2 & 0 & 0\\
            0 & 0 & 2 & 0\\
            -2 & 0 & 0 &2\\
        \end{pmatrix}, \quad &\sigma(L_{A_3}) = \{ 0,2,2,4 \}, \quad \mathsf{deg}(A_3) = 3,\\
        L_{A_4} = \begin{pmatrix}
            2 & 0 & 0 &-2\\
            0 & 4 & 0 & 0\\
            0 & 0 & 4 & 0\\
            -2 & 0 & 0 &2\\
        \end{pmatrix}, \quad &\sigma(L_{A_4}) = \{ 0,4,4,4 \}, \quad \mathsf{deg}(A_4) = 4.
    \end{align*}
    Therefore, 
    \begin{align*}
        \#\pi_0(A_1)=4, \qquad \#\pi_0(A_2)=2, \qquad \#\pi_0(A_3)=1, \quad \text{ and } \quad \#\pi_0(A_4)=1,
    \end{align*} 
    And the existence of cycles condition tells us that 
    \begin{align*}
        &A_1: \qquad\qquad \Tr(L_{A_1})/2=0=\mathsf{rank}(L_{A_1})  &&\Leftrightarrow\qquad \text{is a forest}\\
        &A_2: \qquad\qquad  \Tr(L_{A_2})/2=4/2=\mathsf{rank}(L_{A_2}) &&\Leftrightarrow \qquad\text{is a forest}\\
        &A_3: \qquad\qquad  \Tr(L_{A_3})/2=8/2>3=\mathsf{rank}(L_{A_3}) &&\Leftrightarrow \qquad\text{connected \& contains cycle}\\
        &A_4: \qquad\qquad  \Tr(L_{A_4})/2=12/2>3=\mathsf{rank}(L_{A_4}) &&\Leftrightarrow \qquad {\text{connected \& contains cycle.}}
    \end{align*}
    Recall \cite[\S4.3]{MR4481115} that each of these quantum graphs is quantum isomorphic to a classical graphs with consistent features.

One can similarly check that the complete graph $\cJ$ on $(M_3, \tr)$ is connected and has cycles. As usual, denoting it's adjacency operator as $\hat{J}$, we have 
 \begin{align*}
\hat{J} =
\begin{pmatrix}
3 & 0 & 0 & 0 & 3 & 0 & 0 & 0 & 3 \\
0 & 0 & 0 & 0 & 0 & 0 & 0 & 0 & 0 \\
0 & 0 & 0 & 0 & 0 & 0 & 0 & 0 & 0 \\
0 & 0 & 0 & 0 & 0 & 0 & 0 & 0 & 0 \\
3 & 0 & 0 & 0 & 3 & 0 & 0 & 0 & 3 \\
0 & 0 & 0 & 0 & 0 & 0 & 0 & 0 & 0 \\
0 & 0 & 0 & 0 & 0 & 0 & 0 & 0 & 0 \\
0 & 0 & 0 & 0 & 0 & 0 & 0 & 0 & 0 \\
3 & 0 & 0 & 0 & 3 & 0 & 0 & 0 & 3
\end{pmatrix}, 
\qquad     
L(\hat{J}) =
\begin{pmatrix}
6 & 0 & 0 & 0 & -3 & 0 & 0 & 0 & -3 \\
0 & 9 & 0 & 0 & 0 & 0 & 0 & 0 & 0 \\
0 & 0 & 9 & 0 & 0 & 0 & 0 & 0 & 0 \\
0 & 0 & 0 & 9 & 0 & 0 & 0 & 0 & 0 \\
-3 & 0 & 0 & 0 & 6 & 0 & 0 & 0 & -3 \\
0 & 0 & 0 & 0 & 0 & 9 & 0 & 0 & 0 \\
0 & 0 & 0 & 0 & 0 & 0 & 9 & 0 & 0 \\
0 & 0 & 0 & 0 & 0 & 0 & 0 & 9 & 0 \\
-3 & 0 & 0 & 0 & -3 & 0 & 0 & 0 & 6
\end{pmatrix},
\end{align*}
\[
\sigma(L(\widehat{\mathcal{J}})) =
\left\{ 0, 9, 9, 9, 9, 9, 9, 9, 9 \right\}, \quad
\#\pi_0(A_9) = 1, \quad
\mathsf{deg}(A_9) = 9,
\]

\[
\mathsf{rank}(L(\hat{J})) = 8, \quad
\Tr(L(\hat{J}))/2 = 36.
\]
Therefore, $ \hat{J}$ {\bf is connected and contains cycles.}
\end{ex}

\section{Higher regularity for quantum graphs}
In \cite{2019arXiv190208984E}, Edge used planar algebraic terms to describe certain regularity properties that a classical graph must possess in order to define a spin model of a singly-generated \emph{Yang-Baxter Planar algebra}. 
The main observation here that all those properties can be similarly defined for an arbitrary Schur idempotent over a $\delta$-form independently on whether it defines a classical or a quantum graph. 

\begin{defn} Let $\cG = (\ell^2(X), \hat{T})$ be an undirected and irreflexive quantum graph over the quantum set $X$.
Recall that $\hat{T}$ denotes the  quantum adjacency operator and $T$ its edge projector. We say $\cG$ is {\bf regular with parameter} $k$ if there exists $k \in \mathbb{C}$ such that Equation (\ref{eq:regular}) is satisfied:
\begin{equation}\label{eq:regular}
    \hat{T}\circ \hat{J} = \delta^2\cdot k\cdot\hat{J}
\end{equation}
\end{defn}

\begin{defn}\label{dfn:2PtRegular}
    Let $\cG =(X, \hat{T})$ be an undirected and irreflexive quantum graph.  We say {\bf $\cG$ is 2-point regular/strongly regular} if there exists $k,\lambda, \mu \in \mathbb{C}$ such that Equation (\ref{eq:twoRegular}) holds
    \begin{align}\label{eq:twoRegular}
        \hat{T}^2 = \lambda\hat{T}+ \mu\hat{T}^c+ k\id_X.
    \end{align}
\end{defn}

\begin{defn}\label{defn:WhiteTriangle}
    Let $X$ be a finite quantum set and $\hat{R}, \hat{S}, \hat{T}\in \End(\ell^2(X))$ be arbitrary operators. 
    We then define the following operator:
    \begin{align}\label{eqn:WhiteTriangle}
        {}_{\hat{R}}\triangle_{\hat{S}}^{{\hat{T}}}
        :=(\id_X\otimes m)\circ(\id_X\otimes\hat{T}\otimes\id_X)\circ(m^\dag\otimes \id_X) \circ (\hat{R}\otimes \hat{S})\circ m^\dag=\  
        \tikzmath{
        \begin{scope}
        \clip[rounded corners=5pt] (-1.3,-1.5) rectangle (1.5,2.8);
        \fill[\AColor] (-1.3,-1.5) rectangle (1.5,2.8);
        \end{scope}
        \draw (0,1.8) arc (180:0:.5cm) -- (1,0);
        \draw (.5,2.3) -- (0.5,2.8);
        \filldraw (.5,2.3) circle (.05cm);
        \roundNbox{fill=white}{(0, 1.5)}{.3}{0}{0}{$\hat{T}$}
        \draw (-1,2.8) -- (-1,1.2) arc (180:360:.5);
        \draw (-.5,.25) -- (-.5,.7);
        \filldraw (-.5,.7) circle (.05cm);
        \roundNbox{fill=white}{(-.5,0)}{.3}{0}{0}{$\hat{R}$}
        \roundNbox{fill=white}{(1,0)}{.3}{0}{0}{$\hat{S}$}
        \draw (-.5,-.3) arc (180:360:.75cm);
        \draw (.25,-1) -- (.25,-1.5);
        \filldraw (.25,-1.05) circle (.05cm);
        }.
    \end{align}
\end{defn}
\noindent We will typically only be concerned with the cases in which $\hat{R}, \hat{S}, \hat{T} \in \{ \hat{T}, \hat{J}, \id\}.$

\begin{defn}\label{defn:Triangle-Free}[c.f. \cite[Lemma 3.1.5]{2019arXiv190208984E}]
    We say that a quantum graph $\cG =(X, \hat{T})$ is {\bf triangle-free} if and only if the following relation holds:
\begin{equation}\label{eqn:Triangles}
{}_{\hat{T}}\triangle_{\hat{T}}^{{\hat{T}}}=0.
\end{equation}
\end{defn}

Definitions \ref{defn:WhiteTriangle}  and \ref{defn:Triangle-Free} suggests a way to define the {\bf girth of a quantum graph}. Indeed, as in Diagram (\ref{eqn:WhiteTriangle}), one could connect the vertices of a regular $n$-gon by strings to the outer disk and place an adjacency operator $\hat{T}$ on its edges. (Equivalently one could form an $n$-gon using its edge projector $T$ forming each edge). Then, the girth of $\hat{T}$ would be the least $n$ such that the polygon diagram vanishes if such $n$ exists. 

We justify Definition \ref{defn:Triangle-Free} by the following lemma:
\begin{lem}\label{lem:IndicatorTriangles}
    If $\cG= (V,\hat{T})$ is a classical graph, then the expression ${}_{\hat{T}}\triangle_{\hat{T}}^{{\hat{T}}}$ in Equation (\ref{eqn:Triangles}) indicates if any three given vertices form a triangle. 
    Moreover, if $\cG$ is 2-point regular, $\cG$ is triangle-free if and only if $\lambda = 0.$
\end{lem}
\begin{proof}
    Let $\alpha, \beta, \gamma\in V$ be arbitrary vertices. 
    A direct computation reads
    \begin{align}
        {}_{\hat{T}}\triangle_{\hat{T}}^{{\hat{T}}}: \alpha\mapsto \sum_{\substack{x, y\in V\\ x\sim\alpha\sim y\\ x\sim y}}x\otimes y.
    \end{align}
    Therefore, in the canonical basis for $V$ we get
    \begin{align}
        \left[{}_{\hat{T}}\triangle_{\hat{T}}^{{\hat{T}}}\right]_{\alpha}^{\beta\gamma} = \delta_{\alpha\sim\beta\sim\gamma\sim \alpha}.
    \end{align}
    The rest follows by noting $\lambda$ is the parameter counting how many triangles are formed in $\cG$ from any given pair of adjacent vertices $\alpha\sim \beta.$
\end{proof}

In order to define a \emph{spin model} (c.f. \cite{2019arXiv190208984E} and references therein, and Definition \ref{defn: QuSpinMod}), a graph needs to satisfy a stronger regularity condition, called \emph{3-point regularity}. For a classical graph $\cG$, we say $\cG $ is 3-point regular if the number of common neighbors to any 3 vertices depends only on their configuration inside $\cG$  (c.f. \cite[Definition 2.1.9]{2019arXiv190208984E}). 
    Up to rotations, the possible configurations are:
    \begin{align}\label{eqn:configurations}
        \tikzmath{
        \filldraw (0,1) circle (.05cm);
        \draw (0,0) -- (0,1);
        \draw (0,0) -- (-.86,-.5);
        \draw (0,0) -- (.86,-.5);
        \filldraw (-.86,-.5) circle (.05cm);
        \filldraw (.86,-.5) circle (.05cm);
        \filldraw (0,0) circle (.05cm);
        \node at (.2, .2) {$x$};
        \node at (0,1.2) {$ \alpha$};
        \node at (-1.1,-.5) {$ \beta$};
        \node at (1,.1-.5) {$ \gamma$};
        },\qquad
        \tikzmath{
        \filldraw (0,1) circle (.05cm);
        \draw (0,0) -- (0,1);
        \draw (0,0) -- (-.86,-.5);
        \draw (0,0) -- (.86,-.5);
        \draw(0,1) -- (.86, -.5);
        \filldraw (-.86,-.5) circle (.05cm);
        \filldraw (.86,-.5) circle (.05cm);
        \filldraw (0,0) circle (.05cm);
        \node at (.2, .2) {$x$};
        \node at (0,1.2) {$ \alpha$};
        \node at (-1.1,-.5) {$ \beta$};
        \node at (1,.1-.5) {$ \gamma$};
        }, \qquad
        \tikzmath{
        \filldraw (0,1) circle (.05cm);
        \draw (0,0) -- (0,1);
        \draw (0,0) -- (-.86,-.5);
        \draw (0,0) -- (.86,-.5);
        \draw(0,1) -- (.86, -.5);
        \draw(0,1) -- (-.86, -.5);
        \filldraw (-.86,-.5) circle (.05cm);
        \filldraw (.86,-.5) circle (.05cm);
        \filldraw (0,0) circle (.05cm);
        \node at (.2, .2) {$x$};
        \node at (0,1.2) {$ \alpha$};
        \node at (-1.1,-.5) {$ \beta$};
        \node at (1,.1-.5) {$ \gamma$};
        },\qquad \text{and }
        \tikzmath{
        \filldraw (0,1) circle (.05cm);
        \draw (0,0) -- (0,1);
        \draw (0,0) -- (-.86,-.5);
        \draw (0,0) -- (.86,-.5);
        \draw(0,1) -- (.86, -.5);
        \draw(0,1) -- (-.86, -.5);
        \draw (-.86, -.5) -- (.86, -.5);
        \filldraw (-.86,-.5) circle (.05cm);
        \filldraw (.86,-.5) circle (.05cm);
        \filldraw (0,0) circle (.05cm);
        \node at (.2, .2) {$x$};
        \node at (0,1.2) {$ \alpha$};
        \node at (-1.1,-.5) {$ \beta$};
        \node at (1,.1-.5) {$ \gamma$};
        }.  
    \end{align} 

We now extend this notion to a quantum graph $\cG$ by making use of the observation that substituting instances of $\hat{T}$ by its complement graph $\hat{T}^c$ or the trivial graph $\id_X$ in the left hand side of Equation (\ref{eqn:Triangles}) indicates the configurations described above in the spirit of Lemma \ref{lem:IndicatorTriangles}. 

We shall now construct another operator associated to any quantum graph:
\begin{defn}\label{defn:CommonNeighbors}
    Let $\cG = (\ell^2(X), \hat{T})$ be a quantum graph. 
    We define the operator 
    \begin{align}
        \blacktriangle:= (\hat{T}\otimes\hat{T})\circ m^\dag\circ\hat{T} 
        =
        \tikzmath{
        \begin{scope}
        \clip[rounded corners=5pt] (-.8,-2.5) rectangle (1.2,.5);
        \fill[\AColor] (-.8,-2.5) rectangle (1.2,.5);
        \end{scope}
        \draw (-.25,0) -- (-.25,.5);
        \draw (.75,0) -- (.75,.5);
        \roundNbox{fill=white}{(-.3,-.25)}{.3}{0}{0}{$\hat{T}$}
        \roundNbox{fill=white}{(.7,-.25)}{.3}{0}{0}{$\hat{T}$}
        \draw (-.25,-.55) arc (180:360:.5cm);
        \draw (.25,-1) -- (.25,-2.5);
        \filldraw (.25,-1.05) circle (.05cm);
        \roundNbox{fill=white}{(.25,-1.8)}{-.3}{0}{0}{$\hat{T}$}
        }\in \Hom(\ell^2(X)\to \ell^2(X)\otimes \ell^2(X)).
        \end{align}
\end{defn}

In the classical setting, the map $\blacktriangle$ has a precise combinatorial meaning.
\begin{prop}\label{prop:CommonNeighbors}
    If $\cG=(V, \hat{T})$ is a  classical graph, and $\alpha,\beta, \gamma\in V$ are arbitrary,     
    then
    the entry 
    $$
    [\blacktriangle]_{\alpha}^{\beta\gamma} = \#\{v\in V|\ \hat{T}_{v,\alpha} = \hat{T}_{v,\beta} = \hat{T}_{v,\gamma} = 1\},
    $$
    counts the number of common neighbors $x\in V$ to all $\alpha, \beta$ and $\gamma.$
\end{prop}
\begin{proof}
    Given vertices $\alpha, \beta, \gamma\in V,$
    $$
    \blacktriangle: \alpha \mapsto \sum_{\substack{x,y,z\in V\\ x\sim \alpha\\ y\sim x\sim z}} y\otimes z.
    $$
    Therefore, 
    $$
    \left[\blacktriangle\right]_{\alpha}^{\beta\gamma} = \#\{x\in V|\ x\sim\alpha, x\sim \beta, x\sim \gamma\}.
    $$
\end{proof}

\begin{defn}\label{defn:3ptReg}
    We say that a quantum graph $\cG =(X, \hat{T})$ is {\bf 3-point-regular}  with parameters $(q_0, q_1, q_2, q_3, \lambda, \mu, k)$ if the following identity holds 
    \begin{align}
        \delta^2\cdot\blacktriangle&= q_3\cdot \left[{}_{\hat{T}}\triangle_{\hat{T}}^{{\hat{T}}}\right] + q_2 \left[{}_{\hat{T}^c}\triangle_{\hat{T}}^{{\hat{T}}} + {}_{\hat{T}}\triangle_{\hat{T}^c}^{{\hat{T}}} + {}_{\hat{T}}\triangle_{\hat{T}}^{{\hat{T}^c}}\right] + q_1\cdot\left[{}_{\hat{T}^c}\triangle_{\hat{T}^c}^{{\hat{T}}} + {}_{\hat{T}^c}\triangle_{\hat{T}}^{{\hat{T}^c}} + {}_{\hat{T}}\triangle_{\hat{T}^c}^{{\hat{T}^c}}\right]+ q_0\cdot\left[{}_{\hat{T}^c}\triangle_{\hat{T}^c}^{{\hat{T}^c}}\right]\nonumber\\
        & + \lambda\cdot\left[{}_{\hat{T}}\triangle_{\hat{T}}^{{\id}} + {}_{\hat{T}}\triangle_{\id}^{{\hat{T}}} + {}_{\id}\triangle_{\hat{T}}^{{\hat{T}}}\right] + \mu\cdot\left[{}_{\hat{T}^c}\triangle_{\hat{T}^c}^{{\id}} + {}_{\hat{T}^c}\triangle_{\id}^{{\hat{T}^c}} + {}_{\id}\triangle_{{\hat{T}^c}}^{{\hat{T}^c}}\right]\nonumber\\ 
        &+ k\cdot\left[{}_{\id}\triangle_{\id}^{{\id}}\right]\!.
    \end{align}
\end{defn}
\noindent Notice that our expressions make no explicit reference to state sums not checkerboard shadings in the planar algebraic style of \cite[Lemma 3.1.4]{2019arXiv190208984E}, which we believe helps clarify the combinatorial significance of individual string diagrams. 

\begin{remark}
    If $\cG = (V, \hat{T})$ is a classical 3-point regular classical graph, the parameters $q_0, q_1, q_2,$ and $ q_3$ respectively give the number of common neighbors to arbitrary vertices $\alpha, \beta, \gamma$ as a function of the configuration of these three points (cf Equation (\ref{eqn:configurations})).
    Notice $\cG$ will automatically be 2-point regular (i.e. strongly regular) with parameters $(k, \lambda, \mu)$. 
\end{remark}

We shall now provide examples of 2-point/3-point regular quantum graphs over various spaces.
\begin{ex}\label{ex:M2TripleRegularity}
    We consider the quantum space $X=(M_2, \tr)$. By the classification of quantum graphs on \cite{MR4481115}, we will only consider the quantum graphs $A_3$ and $A_4$ with the notation from Example \ref{ex:LaplaciansM2}, since these are all possible graphs up to complementation. 
    One can compute directly that both these graphs are 3-point regular with the following parameters (considered as irreflexive quantum graphs):
    \begin{table}[h!]
    \centering
    \begin{tabular}{l c c c c c c c}
    {} & $k$ & $\lambda$ & $\mu$ & $q_3$ & $q_2$ & $q_1$ & $q_0$ \\
    \midrule
    $A_4:$        & 3  & 2 & 0 & 1 & 0 & 0 & 0 \\
    $A_3:$        & 2  & 0 & 2 & 0 & 0 & 0 & 0   
    \end{tabular}
    \end{table}

    Since these graphs are quantum isomorphic to classical graphs which are on the nose 3-point regular with the corresponding parameters as above, by Theorem \ref{thm:QuInvariants}, we can reach the same conclusions using the machinery of Section \ref{sec:bubble}.
\end{ex}

\begin{ex}\label{ex:M3TripleRegularity}
  {  We now consider quantum graphs defined over $X=(M_3, \tr)$ and we shall give examples of 3-point regular graphs. 

    The graphs $G_7$ and $G_3$ in the table below are each other's complements, and were obtained by applying the quantum Fourier transform to the diagonal projection $\mathsf{diag}(1,1,1,0,0,0,0,0,0)$.
    
    The graphs $G_6$ and $G_4$ below are complementary and were constructed as follows: 
    We begin by constructing a projection onto a 4-dimensional subspace \( V \) of \( M_3(\mathbb{C}) \). This subspace is defined as the linear span of the following matrices:
\[
V = \mathbb{C} \begin{pmatrix}
        1 & 0 & 0\\
        0 & 1 & 0\\
        0 & 0 & 1\\
    \end{pmatrix} \oplus \mathbb{C} \begin{pmatrix}
        0 & 1 & 0\\
        1 & 0 & 0\\
        0 & 0 & 0\\
    \end{pmatrix} \oplus \mathbb{C} \begin{pmatrix}
        0 & 0 & 1\\
        0 & 0 & 0\\
        1 & 0 & 0\\
    \end{pmatrix} \oplus \mathbb{C} \begin{pmatrix}
        0 & 0 & 0\\
        0 & 0 & 1\\
        0 & 1 & 0\\
    \end{pmatrix}.
\]

The projection matrix \( P \) onto \( V \) is constructed by taking its columns to be the vectorized forms of the basis matrices \( M_i \) of \( V \). The vectorization operation, \(\text{vec}: M_n(\mathbb{C}) \to \mathbb{C}^{n^2}\), is defined as:
\[
\text{vec}(A) = (a_{11}, \ldots, a_{1n}, a_{21}, \ldots, a_{2n}, \ldots, a_{n1}, \ldots, a_{nn})^T,
\]
for a matrix \( A = (a_{ij}) \). The columns of \( P \) are arranged to ensure the matrix is symmetric.

Subsequently, each non-zero entry (specifically, each entry of 1) in the initial matrix \( P \) is replaced by a free parameter \( p_i \). A brute-force numerical search was then conducted over the parameter space, scanning values within the interval \([0.5, 2]\), to identify a parameter set for which the matrix becomes a Schur idempotent. The resulting matrix, denoted \( G_4^{r} \), displays the specific parameter values that satisfy this condition. It can be verified that, $G_4^r$ is a \textbf{reflexive} quantum graph. In the following, we call $G_4:= G_4^{r} - Id$, the irreflexive version of $G_4^r$.

    The following table contains {\bf irreflexive real undirected/bidirected $3$-point regular quantum graphs on $M_3$}: }
    \begin{table}[h!]
    \centering
    \begin{tabular}{l c c c c c c c c c c c}
    {} & $k$ & $\lambda$ & $\mu$ & $q_3$ & $q_2$ & $q_1$ & $q_0$ & $\Tr(L)$  & $\mathsf{rank}(L)$ & $\#\pi_0$  & $\exists$\text{ cycle}\\
    \midrule
    $G_3:$        & 2  & 1    & 0    & 0      & 0     & 0     & 0    & 18    & 6     & 3    & Y\\    
    $G_4:$        & 3  & 6/8  & 6/8  & 3/8    & -3/8  & 3/8   & -3/8 & 27    & 8     & 1    & Y\\
    $9P_q:$       & 4  & 1    & 2    & 0      & 0     & 1     & 0    & 36    & 8    & 1    & Y\\   
    $G_6:$        & 5  & 14/8 & 30/8 & -3/8   & 7/8   & 21/8  & 15/8 & 45    & 8     & 1    & Y\\
    $G_7:$        & 6  & 3    & 6    & 0      & 3     & 0     & 6    & 54    & 8     & 1    & Y\\   
    $\cJ:$    & 8  & 7    & 0    & 6      & 0     & 0     & 0    & 72    & 8     & 1    & Y
    \end{tabular}
    \end{table}\\
    Here, $9P_q$ is the quantum graph analogous to the classical 9-Paley graph whose adjacency matrix is given in Construction \ref{const:q9Paley}, and is self-complimentary.  
    The adjacency matrices of remaining quantum graphs are expressed below in the orthonormal basis $\{u_{a,b}=\sqrt{3}e_{a,b}\}_{a,b=1}^3$.
$$
G_3 =
\begin{pmatrix}
2 & 0 & 0 & 0 & 0 & 0 & 0 & 0 & 0 \\
0 &-1 & 0 & 0 & 0 & 0 & 0 & 0 & 0 \\
0 & 0 &-1 & 0 & 0 & 0 & 0 & 0 & 0 \\
0 & 0 & 0 &-1 & 0 & 0 & 0 & 0 & 0 \\
0 & 0 & 0 & 0 &2 & 0 & 0 & 0 & 0 \\
0 & 0 & 0 & 0 & 0 &-1 & 0 & 0 & 0 \\
0 & 0 & 0 & 0 & 0 & 0 &-1 & 0 & 0 \\
0 & 0 & 0 & 0 & 0 & 0 & 0 &-1 & 0 \\
0 & 0 & 0 & 0 & 0 & 0 & 0 & 0 &2
\end{pmatrix}, 
\
G_4 =
\begin{pmatrix}
0 & 0 & 0 & 0 & \frac{3}{2} & 0 & 0 & 0 & \frac{3}{2} \\
0 & 0 & 0 & \frac{3}{2} & 0 & 0 & 0 & 0 & 0 \\
0 & 0 & 0 & 0 & 0 & 0 & \frac{3}{2} & 0 & 0 \\
0 & \frac{3}{2} & 0 & 0 & 0 & 0 & 0 & 0 & 0 \\
\frac{3}{2} & 0 & 0 & 0 & 0 & 0 & 0 & 0 & \frac{3}{2} \\
0 & 0 & 0 & 0 & 0 & 0 & 0 & \frac{3}{2} & 0 \\
0 & 0 & \frac{3}{2} & 0 & 0 & 0 & 0 & 0 & 0 \\
0 & 0 & 0 & 0 & 0 & \frac{3}{2} & 0 & 0 & 0 \\
\frac{3}{2} & 0 & 0 & 0 & \frac{3}{2} & 0 & 0 & 0 & 0
\end{pmatrix},
$$

$$\mathsf{Spec}(G_3) = \{2,-1,-1,-1,2,-1,-1,-1,2\},$$
$$\mathsf{Spec}(G_4) = \{-3/2,  3,   3/2, -3/2,  3/2, -3/2, -3/2,  3/2, -3/2\}$$

In fact, as is mentioned in \cite[Remark in p 39]{MR1188082} $G_4$ satisfies $\mu = (-3/2)^2 + (-3/2)$, which is the \emph{pseudo Latin square type graph} condition, and is not a conference graph. 

$$
G_6 =
\begin{pmatrix}
2 & 0 & 0 & 0 & \frac{3}{2} & 0 & 0 & 0 & \frac{3}{2} \\
0 &-1 & 0 & -\frac{3}{2} & 0 & 0 & 0 & 0 & 0 \\
0 & 0 &-1 & 0 & 0 & 0 & -\frac{3}{2} & 0 & 0 \\
0 & -\frac{3}{2} & 0 &-1 & 0 & 0 & 0 & 0 & 0 \\
\frac{3}{2} & 0 & 0 & 0 & 2 & 0 & 0 & 0 & \frac{3}{2} \\
0 & 0 & 0 & 0 & 0 & -1 & 0 & -\frac{3}{2} & 0 \\
0 & 0 & -\frac{3}{2} & 0 & 0 & 0 & -1 & 0 & 0 \\
0 & 0 & 0 & 0 & 0 & -\frac{3}{2} & 0 & -1 & 0 \\
\frac{3}{2} & 0 & 0 & 0 & \frac{3}{2} & 0 & 0 & 0 & 2
\end{pmatrix},\ 
G_7 =
\begin{pmatrix}
0 & 0 & 0 & 0 & 3 & 0 & 0 & 0 & 3 \\
0 & 0 & 0 & 0 & 0 & 0 & 0 & 0 & 0 \\
0 & 0 & 0 & 0 & 0 & 0 & 0 & 0 & 0 \\
0 & 0 & 0 & 0 & 0 & 0 & 0 & 0 & 0 \\
3 & 0 & 0 & 0 & 0 & 0 & 0 & 0 & 3 \\
0 & 0 & 0 & 0 & 0 & 0 & 0 & 0 & 0 \\
0 & 0 & 0 & 0 & 0 & 0 & 0 & 0 & 0 \\
0 & 0 & 0 & 0 & 0 & 0 & 0 & 0 & 0 \\
3 & 0 & 0 & 0 & 3 & 0 & 0 & 0 & 0
\end{pmatrix}\!.
$$
Their spectra is given by 
$$\mathsf{Spec}(G_6) = \{1/2, 5, -5/2, -5/2, 1/2, 1/2, 1/2, 1/2, -5/2\},$$ 
$$\mathsf{Spec}(G_7) = \{-3,  6, -3,  0,  0,  0,  0,  0,  0\}.$$
We observe that the parameter transformations between the complementary pairs of graphs above satisfy \cite[Proposition 3.23]{2024arXiv240406157G}, but were obtained directly aided by computer software.
Notice, moreover that $G_3$ is disconnected by our Laplacian analysis, and this conclusion is consistent with Gromada's \cite[Proposition 3.27]{2024arXiv240406157G} for strongly regular graphs with $\mu=0.$ 
A similar conclusion follows for the irreflexive $\cJ.$  
Our construction of $G_4$ gives a new example of a strongly regular quantum graph with $\mu\neq 0.$ 
The computational verification of this claim in our Software's Section \href{https://colab.research.google.com/drive/11rctqnjD3bxMejfEKWnJfr2qUy6oUYA3?usp=sharing}{4.2}. 
\end{ex}

Recall that when $\lambda =0 $ the graph is triangle-free, and so $q_0$ can take any complex value. 
Similarly, when $\mu =0$, the parameters $q_1$ and $q_2$ are free. 
\begin{remark}
    
    The graph $G_7$ from Example \ref{ex:M3TripleRegularity} shares all its parameters with the classical {\bf complete tripartite graph} $\cK_{3,3,3}$, which can readily be seen to be connected and $3$-point regular. 
    This suggests that $\cK_{3,3,3}$ can be transported to $M_3$ via a quantum isomorphism $M_3\cong_q\bbC^9.$
\end{remark}

\begin{ex}\label{ex:params}
    In  Constructions \ref{const:q16Clebsch} and \ref{const:qHS} we shall later obtain 3-point regular quantum graphs over the tracial spaces $M_4$ and a noncommutative C*-algebra of dimension $100$ which are quantum isomorphic to the classical $16$-Clebsch and the Higman-Sims graphs, respectively, whose regularity parameters are given by (c.f. \cite[Table II]{MR1247144}, \cite{2019arXiv190208984E}):
    \begin{table}[h!]
    \centering
    \begin{tabular}{l c c c c c c c c c c c c }
    {} & $\mathsf{dim}$ & $k$ & $\lambda$ & $\mu$ & $q_3$ & $q_2$ & $q_1$ & $q_0$ & $\Tr(L)$  & $\mathsf{rank}(L)$ & $\#\pi_0$  & $\exists$\text{ cycle}\\
    \midrule
        $16Cl_q:$        & 16  & 5  & 0   & 2    & 1& 0 & 0 & 0 & 80   & 15 &1 & Y\\   
        $HS_q:$          &100  & 22 & 0   & 6    & 2     & 0 & 0 & 0 & 2200 & 99 &1  & Y\\ 
    \end{tabular}
    \end{table}
\end{ex}
Thee values for $16Cl_q$ were also verified by direct computation aided by computer software. 
However, all values can be obtained from the classical $16$-Clebsch and Higman-Sims graph thanks to Theorem  \ref{thm:QuInvariants}.

\section{Deforming quantum graphs by bubbling}\label{sec:bubble}
We now summarize a construction from \cite[\S3]{MR3907958} used to deform a given (quantum) graph into another quantum isomorphic (quantum) graph. 
\begin{construction}\label{const:Bubbling}[{\bf Bubbling quantum graphs}]
    Let $\cG=(\ell^2(X), \hat{T})$ be a real undirected quantum graph, and let
    \begin{align*}
    \bbX: (H\otimes \overline{H})\otimes &\ell^2(X)\to \ell^2(X)\otimes (H\otimes \overline{H})\\
    &\tikzmath
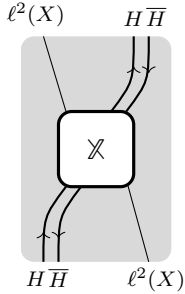
{
        \begin{scope}
        \clip[rounded corners=5pt] (-1,-1.5) rectangle (1,1.5);
        \fill[\AColor] (-1,-1.5) rectangle (1,1.5);
        \end{scope}
         \draw[thick] (-.7, -1.5) .. controls (-.7,-.9) .. (-.4, -.5);
         \draw[thick] (-.5, -1.5) .. controls (-.5,-.9) .. (-.2, -.5);
         \draw[->] (-.7,-1.1);
         \draw[>-] (-.5,-1.1);
         \draw (.4, -.4) -- (.7, -1.5);
         \draw[thick] (.2,.5) .. controls (.5,.9) .. (.5,1.5);
         \draw[thick] (.4,.5) .. controls (.7,.9) .. (.7,1.5);      \draw[->] (.5,1.1);   
         \draw[>-] (.7,1.1);
         \draw (-.4,.5) -- (-.7, 1.5);
        \roundNbox{fill=white}{(0,0)}{.5}{0}{0}{$\mathbb{X}$}
        \draw (-.8,-1.5) node[below]{${\scriptstyle H}$};
        \draw (-.5,-1.45) node[below]{${\scriptstyle \overline{H}}$};
        \draw (.8,1.5) node[above]{${\scriptstyle \overline{H}}$};
        \draw (.5,1.5) node[above]{${\scriptstyle H}$};
        \draw (-.8,1.5) node[above]{${\scriptstyle \ell^2(X)}$};
        \draw (.8,-1.45) node[below]{${\scriptstyle \ell^2(X)}$};
    }
    \end{align*}
    be a quantum graph isomorphism (i.e. $\mathbb{X}$ satisfies Equations  (\ref{eqn:QuantumFunction}), (\ref{eqn:QuantumBijection}) and (\ref{eqn:QuantumGraphHom})) which is assumed to be a Frobenius C*-algebra in the C*-tensor category $\mathsf{QAut}(\cG)\cong \Rep(C(\Aut^+(\cG)))$ (c.f. Equivalence (\ref{eqn:QAutRep})). Our goal is to use a \emph{dagger splitting} of $\mathbb{X}$ by some quantum bijection $P$: 
    \begin{align}\label{eqn:Splitting}
    \mathbb{X} = P\bullet P^*,    
    \end{align}
    to construct a new quantum graph $\cG_q = (\ell^2(X_q), \hat{T_q})$ from $P$ that is quantum isomorphic to $\cG.$

    Consider the map
    \begin{align}\label{eqn:pi}
    \pi: \overline{H}\otimes& \ell^2(X)\otimes H\to \overline{H}\otimes \ell^2(X)\otimes H\\
    &\dfrac{1}{\mathsf{dim}(H)}\cdot\tikzmath
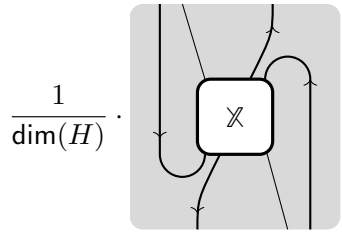
{
        \begin{scope}
        \clip[rounded corners=5pt] (-1.4,-1.5) rectangle (1.4,1.5);
        \fill[\AColor] (-1.4,-1.5) rectangle (1.4,1.5);
        \end{scope}
         \draw[thick] (-.4,-.5) arc (360:180:.3)  -- (-1, 1.5);
         \draw[->](1,.5);
         \draw[>-](-.5,-1.2);
         \draw[thick] (-.5, -1.5) .. controls(-.5,-1.1) .. (-.2, -.5);
         \draw(.4, -.4)-- (.7, -1.5);
         \draw[->](.5,1.2);
         \draw[>-](-1,-.2);
         \draw[thick] (.2,.5)  ..controls (.5,1.1) .. (.5,1.5);
         \draw[thick] (.4,.5) arc (180:0:.3)  -- (1, -1.5);
         \draw (-.4,.5) -- (-.7, 1.5);
        \roundNbox{fill=white}{(0,0)}{.5}{0}{0}{$\mathbb{X}$}
    },\nonumber
    \end{align}
    which is readily seen to be an orthogonal projection $\pi = \pi^2 = \pi^\dag.$
    Let's now consider the Hilbert space $\ell^2(X_q)$ defined as  the range of $\pi,$ and consider the splitting isometry
    $$
    \iota: \ell^2(X_q)\to \overline{H}\otimes \ell^2(X)\otimes H
    $$
    satisfying \cite[Equation (62)]{MR3907958}
    \begin{align}\label{eqn:SplittingPi}
        &\iota\circ\iota^\dag = \pi\qquad\qquad\qquad \text{ and} &&\iota^\dag\circ\iota  = \id_{\ell^2(X_q)}.
    \end{align}
    Finally, consider the following map
    \begin{align}\label{eqn:P}
        P: H\otimes&\ell^2(X_q)\to \ell^2(X)\otimes H,\\
        &P:=\sqrt{\dim(H)}\cdot\tikzmath{
        \begin{scope}
        \clip[rounded corners=5pt] (-1,-1) rectangle (1,1);
        \fill[\AColor] (-1,-1) rectangle (1,1);
        \end{scope}
         \draw(0,-1) --(0,-.1);
         \draw(0,.1) -- (0,1);
         \draw [->] (.5,.7);
         \draw [->] (-.5,-.7);
         \draw[thick] (-.5,-1) .. controls (-.5,-.4) .. (-.1,-.1);
         \draw[thick] (.5,1) .. controls (.5,.4) .. (.1,.1);
        \roundNbox{fill=white}{(0,0)}{.15}{0}{0}{}
        }
        :=
        \tikzmath{
        \begin{scope}
        \clip[rounded corners=5pt] (-1.2,-1) rectangle (1,1);
        \fill[\AColor] (-1.2,-1) rectangle (1,1);
        \end{scope}
         \draw(0,-1) --(0,-.3);
         \draw(0,.3) -- (0,1);
         \draw[thick] (.2,.3) ..controls (.5,.6).. (.5,1);
         \draw[->] (.5,.8);
         \draw[thick] (-.2,.3) arc (0:180:.2);
         \draw[thick] (-.6, .3) .. controls(-.7,-.6).. (-.7, -1);
         \draw[ ->] (-.7,-.6);
        \roundNbox{fill=white}{(0,0)}{.3}{0}{0}{$\iota$}
        }.\nonumber 
    \end{align}\
    satisfying $P\bullet P^* =\mathbb{X}$. 
    For simplicity, we will use the following diagram to denote the conjugate of $P$ as no ambiguity may arise:
    \begin{align}
     P^{*}: \overline{H}\otimes&\ell^2(X)\to \ell^2(X_q)\otimes \overline{H},\\
        &P^*:=\tikzmath{
        \begin{scope}
        \clip[rounded corners=5pt] (-1,-1) rectangle (1,1);
        \fill[\AColor] (-1,-1) rectangle (1,1);
        \end{scope}
         \draw(0,-1) --(0,-.1);
         \draw(0,.1) -- (0,1);
         \draw [>-] (.5,.7);
         \draw [>-] (-.5,-.7);
         \draw[thick] (-.5,-1) .. controls (-.5,-.4) .. (-.1,-.1);
         \draw[thick] (.5,1) .. controls (.5,.4) .. (.1,.1);
        \roundNbox{fill=white}{(0,0)}{.15}{0}{0}{}
        }.\nonumber    
    \end{align}

    We now equip the Hilbert space $\ell^2(X_q)$ with the structure of a quantum space: 
    \begin{align}\label{eqn:Xq}
        m_q = \tikzmath{
        \begin{scope}
        \clip[rounded corners=5pt] (-2,-1.5) rectangle (2,2);
        \fill[\AColor] (-2,-1.5) rectangle (2, 2);
        \end{scope}
        \draw (0,2) node[above]{$\scriptstyle \ell^2(X_q)$};
        \draw (-.5,-1.5) -- (-.5, -.2) -- (-.5,0) arc (180:0:.5cm) --(.5,-1.5);
        \draw (0,.5) -- (0,2);
        \filldraw (0,.5) circle (.05cm);
        \draw[thick] (-.4,-.8) .. controls (-.1,-.7) .. (.0,-.4);
        \draw[thick] (0,-.4)  .. controls (.1,-.1) .. (.5,-0);
        \draw[thick] (.5,0) -- (.7,.2) .. controls (.9,.3) .. (.9,1.2) arc (0:180:.41);
        \draw[thick] (-.6,-.8) arc (360:180:.3) -- (-1.2, 0) .. controls (-1.2,.2).. (0,1.2) ;
        \draw [->] (.9,1);
        \draw [>-] (-1.2,-.6);
        \draw [->] (0,-.35);
        \draw (-.6,-1.5) node[below]{$\scriptstyle \ell^2(X_q)$};
        \draw (0.6,-1.5) node[below]{$\scriptstyle \ell^2(X_q)$};
        \roundNbox{fill=white}{(-.5,-.8)}{.15}{0}{0}{}
        \roundNbox{fill=white}{(.5,0)}{.15}{0}{0}{}
        \roundNbox{fill=white}{(0,1.2)}{.15}{0}{0}{}
    }
    \qquad \text{ and }\qquad 
    i_q=
    \tikzmath{
        \begin{scope}
        \clip[rounded corners=5pt] (-1,-1) rectangle (1.2,2);
        \fill[\AColor] (-1,-1) rectangle (1.2,2);
        \end{scope}
        \draw (0,2) node[above]{${\scriptstyle \ell^2(X_q)}$};
        \draw (0,.5) -- (0,2);
        \draw[thick] (0,1) .. controls (.4,1.5) ..  (.6,1) ;
        \draw[thick] (.6,1) .. controls (.7,.7)  ..  (.3,.1);
        \draw[thick] (.3,.1) .. controls (0,-.3) .. (-.4,.1);
        \draw[thick] (-.4,.1) .. controls (-.5,.3) ..  (0,1);
        \draw [>-] (.65,.8);
        \draw[->](-.45,.3);
        \filldraw (0,.5) circle (.05cm);
        \roundNbox{fill=white}{(0,1)}{.15}{0}{0}{}
    }.
    \end{align}
    From the defining properties of quantum isomorphisms, it is easy to see that $X_q$ is again a $\delta$-form.
    With this structure, it is readily seen that $P$ is a quantum bijection between the quantum spaces $\ell^2(X)$ and $\ell^2(X_q)$.

    We are now ready to {\bf bubble $\hat{T}$} to obtain a $\hat{T}_q\in \End(\ell^2(X_q))$ as follows:
    \begin{align}\label{eqn:BuubleTea}
    \hat{T}_q  &:=\dfrac{1}{\dim(H)}\cdot (\id_X\otimes \ev_H)\circ (P^*\otimes \id_H) \circ (\id_{\bar{H}}\otimes \hat{T}\otimes \id_H)
     \circ (\id_{\bar{H}}\otimes \id_H\otimes P) \circ 
     (\ev_H^\dag\otimes \id_X)\nonumber\\
    &=
    \dfrac{1}{\dim(H)}\cdot
    \tikzmath{
        \begin{scope}
        \clip[rounded corners=5pt] (-1.2,-2) rectangle (1.2,2);
        \fill[\AColor] (-1.2,-2) rectangle (1.2,2);
        \end{scope}
         \draw(.3,-1) .. controls (.3,-.8) .. (0,-.3);
         \draw(.5,-1) .. controls (.8,-1.3) .. (.8,-2);
         \draw(-.3,1) .. controls (-.3,.8) .. (0,.3);
         \draw(-.5,1) .. controls (-.8,1.3) .. (-.8,2);
         \draw[thick] (-.3,1.1) -- (-.2,1.2) arc (160:0:.5) -- (.7,-.5) .. controls (.7,-.9) .. (.45,-1);
         \draw[thick] (.3,-1.1) -- (.2,-1.2) arc (340:180:.5) -- (-.7,.5) .. controls (-.7,.9) .. (-.45,1);
         \draw [->] (.75,.4);
         \draw[>-] (-.75,-.2);
        \roundNbox{fill=white}{(0,0)}{.3}{0}{0}{$\hat{T}$}
        \roundNbox{fill=white}{(.4,-1)}{.15}{0}{0}{}
        \roundNbox{fill=white}{(-.4,1)}{.15}{0}{0}{}
    }
    \end{align}
    Which can readily be seen to define a $\star$-idempotent. 
    Notice that $\cG_q$ is real/un-directed/reflexive if and only if $\cG$ is such, and that any arbitrary operator on $\ell^2(X)$ and its tensor powers can be bubbled by similar means whether or not is a $\star$-idempotent.
\end{construction}

We now show that {\bf bubbling $\hat{T}$ preserves its spectrum along with multiplicities}. This is non-trivial since this process involves an amplification followed by conjugation by an isometry rather than merely a unitary conjugation. 
The first step in this direction is to show that bubbling commutes with polynomials in $\hat{T}.$

\begin{lem}\label{lem:BubbleProperties}
    With the notation and assumptions of Construction \ref{const:Bubbling}, if $f(z)\in \bbC[z]$ is a polynomial, then 
    $$
        f(\hat{T})_q =  f(\hat{T}_q).
    $$
    Consequently, as sets, 
    $$
    \Spec\left(\hat{T}\right)  = \Spec\left(\hat{T}_q\right).
    $$
    Additionally, if $\{Q_\lambda\}_{\lambda\in \Spec(\hat{T})}$ is the set of spectral projections for $\hat{T}$ and $\{R_\lambda\}_{\lambda\in \Spec(\hat{T}_q)}$ is that of $\hat{T}_q$, then for all $\lambda\in\Spec(\hat{T})$
    $$
    \bbX\bullet Q_\lambda = Q_\lambda\bullet\bbX\quad\text{ and }\quad (Q_\lambda)_q = R_\lambda. 
    $$
\end{lem}
\begin{proof}
    Since $\hat{T}$ and $\bbX$ commute, for all $k\in\bbN$ it holds that $(\hat{T}^k)_q = (\hat{T}_q)^k.$ 
    The first claim then follows by linearity.  

    Let $m_{\hat{T}}(z)$ be the minimal polynomial for $\hat{T}$ and similarly $m_{\hat{T}_q}(z)$ for $\hat{T}_q.$ 
    From the claim above, we have that 
    $$
    m_{\hat{T}}(\hat{T}) = 0 = m_{\hat{T}}(\hat{T}_q).
    $$
    Thus,  
    $$
    m_{\hat{T}}(z) = m_{\hat{T}_q}(z).
    $$
    Since minimal polynomials uniquely factor over $\bbC$ as 
    $$
    m_{\hat{T}}(z) = \Pi_{\lambda\in\Spec(\hat{T})}(z-\lambda\id),
    $$
    the equality of sets $\Spec(\hat{T})= \Spec(\hat{T})_q$ follows.

    Since $\hat{T}$ and $\hat{T}_q$ are self-adjoint, they are diagonalizable, with spectral decompositions 
     $$
     \hat{T} = \sum_{\lambda\in \Spec(\hat{T})}\lambda Q_\lambda\quad \text{ and }\quad \hat{T}_q = \sum_{\lambda\in \Spec(\hat{T})}\lambda R_\lambda.
     $$
     For any $\lambda\in \Spec({\hat{T}}),$ the indicator function $1_{\{\lambda\}}(z) = \prod_{\mu\neq \lambda}\dfrac{z-\mu}{\lambda - \mu}$ is a polynomial. 
    Thus, the spectral projections are polynomials in their respective adjacency operator.
    Indeed, 
    $$
    Q_\lambda = 1_{\{\lambda\}}(\hat{T}) \quad \text{ and }\quad R_\lambda = 1_{\{\lambda\}}(\hat{T}_q).
    $$
    It follows at once that 
    $$
    \bbX\bullet Q_\lambda = Q_\lambda \bullet \bbX,
    $$
    since 
    $$
    \bbX\bullet \hat{T} = \hat{T} \bullet \bbX
    $$
    holds. 
    Similarly, recalling that $\bbX = P\bullet P^*,$ since 
    $$
    P\bullet \hat{T}_q = \hat{T} \bullet P,
    $$
    it follows that 
    $$
    P\bullet R_\lambda = Q_\lambda \bullet P.
    $$
    This implies the last claim
    $$
    (Q_\lambda)_q = R_\lambda.
    $$
\end{proof}

\begin{remark}
    That bubbling preserves the set spectra can be shown directly: 
    For $\lambda\in \mathbb{C}$, set $S:=S(\lambda):= \hat{T}-\lambda\id,$ and we denote its bubbled version by $S_q.$
    Similarly, bubbling in the reversed direction is denoted $S_{\overline{q}}.$

    We shall now show that $\mathsf{spec}\left(\hat{T}\right) = \mathsf{spec}\left(\hat{T}_q\right)$ considered as sets. 
    That is, $\hat{T}$ and $\hat{T}_q$ have the same eigenvalues. 
    Let $\lambda\notin\mathsf{spec}(\hat{T})$ so that there exists $R$ with $RS = \id = SR.$ Then a short diagrammatic computation yields $R_qS_q = (RS)_q = \id = S_qR_q,$ since $S$ commutes with $\bbX$. 
    This shows $\lambda\notin \mathsf{spec}(\hat{T}_q)$.
    
    Conversely, if $\lambda\notin\mathsf{spec}(\hat{T}_q)$, there exists $R$ such that $S_qR= \id = RS_q.$ 
    A short diagrammatic computation yields $SR_{\overline{q}} = (S_qR)_q = \id_{\overline{q}} = q = R_{\overline{q}}S.$
    Therefore, $\Spec(\hat{T}) = \Spec(\hat{T}_q)$ as sets. 
\end{remark}

\begin{prop}\label{prop:BubleSpectrum}
    With the notation and assumptions of Construction \ref{const:Bubbling}, we have that 
    for every $S\in \End(\ell^2X),$ 
    $$
    \Tr\left(S\right) = \Tr\left(S_q\right).
    $$
    Furthermore, for $k\in \bbN,$ we have 
    $$
    \Tr(\hat{T}^k) = \Tr(\hat{T}_q^k).
    $$
    Consequently, 
    $$
    \mathsf{spec}\left(\hat{T}\right) = \mathsf{spec}\left(\hat{T}_q\right)
    $$ 
    counted with multiplicities. 
    
\end{prop}
\begin{proof}
    That bubbling preserves the trace of $S$ follows from a single diagrammatic computation using the Duality Equations (\ref{eqn:DiagDualityQAut}). 
    
    That $\Tr(\hat{T}^k) = \Tr(\hat{T}_q^k)$ follows from Lemma  \ref{lem:BubbleProperties} and the previous claim. 
     
     Finally, we establish that $\hat{T}$ and $\hat{T}_q$ have the same spectra counting algebraic multiplicities. 
     We shall maintain the notation from the proof of Lemma \ref{lem:BubbleProperties}.  
     Notice that since we are considering self-adjoint $\hat{T}$, algebraic and geometric multiplicities match for $\hat{T}$ and $\hat{T}_q$, respectively. 
     Algebraic multiplicities are given by 
     $$
     \mathsf{AlgMult}_{\hat{T}}(\lambda) = \Tr(Q_\lambda)\quad \text{ and }\quad \mathsf{AlgMult}_{\hat{T}_q}(\lambda) = \Tr(R_\lambda)
     $$
     for all $\lambda\in \Spec(\hat{T}),$ respectively. 
     Since $\Tr(\hat{T}^k) = \Tr(\hat{T}_q^k)$ for all $k\in \bbN,$ for any polynomial $f(z) = \sum_{\ell=0}^L c_\ell z^\ell$ we have that 
     $$
     \Tr(f(\hat{T})) = \sum_{\ell=0}^Lc_\ell\Tr(\hat{T}^\ell) = \sum_{\ell=0}^L c_\ell\Tr(\hat{T}_q^\ell) = \Tr(f(\hat{T}_q)). 
     $$
    Finally, we obtain for every $\lambda$
    $$
    \mathsf{AlgMult}_{\hat{T}}(\lambda) =  \Tr(Q_\lambda) = \Tr(1_{\{\lambda\}}(\hat{T})) = \Tr(1_{\{\lambda\}}(\hat{T}_q)) = \Tr(R_\lambda) = \mathsf{AlgMult}_{\hat{T}_q}(\lambda).
    $$
     This completes the proof. 
\end{proof}

Upon review of the proofs of the results in Lemma \ref{lem:BubbleProperties} and Proposition \ref{prop:BubleSpectrum}, we notice that the fundamental identity used was that 
$$
    \hat{T}\bullet \bbX = \bbX\bullet \hat{T}.
$$
However, one can directly show by a quick diagrammatic check that 
\begin{align}\label{eqn:Commutators}
    L(\hat{T})\bullet \bbX = \bbX \bullet L(\hat{T}). 
\end{align}
Therefore, similar statements to those in Lemma \ref{lem:BubbleProperties} and Proposition \ref{prop:BubleSpectrum} hold for $L(\hat{T})$, proven by similar arguments. 
We record this observation in the following corollary.
\begin{cor}\label{cor:BubblingEigenspaces}
    Let $\cG = (X,\hat{T})$ be a real undirected quantum graph and $P\bullet P^* = \bbX \in \QAut(\cG)$ as in Construction \ref{const:Bubbling}.  
    For $S\in \{\hat{T}, L(\hat{T})\}$ we have that, if $f(z)\in \bbC[z]$ is a polynomial, then 
    $$
        f(S)_q =  f(S_q),
    $$
    and 
    $$
    \mathsf{Spec}(S) = \mathsf{Spec}(S_q),
    $$
    accounting for multiplicities.
    Moreover, for every $\lambda\in \mathsf{Spec}(S)$ we have  
    $$
        \mathsf{Eigenspace}_S(\lambda)  \cong \mathsf{Eigenspace}_{S_q}(\lambda).
    $$
    That is, bubbling $\hat{T}$ or $L(\hat{T})$ preserves their respective spectra and eigenspaces. 
\end{cor}
\begin{proof}
It suffices to check the claims involving the Laplacian only. 

The claim that $f(L(\hat{S}))_q =  f(L(\hat{S}_q))$ follows by checking directly that bubbling commutes with words in $\hat{T}$ and $\mathsf{Deg}(\hat{T}).$
The remainder of the proof follows similarly as those of Lemma \ref{lem:BubbleProperties} and Proposition \ref{prop:BubleSpectrum}, and hence we omit the remaining details. 
\end{proof}

We shall now prove that bubbling commutes with with the inverse Quantum Fourier Transform; that is, bubbling is compatible with taking the edge projector.
\begin{prop}\label{prop:BubbleEdgeProjector}
    Let $\cG = (X,\ \hat{T})$ be a quantum graph with edge projector $T.$
    Then 
    $$
    T_q = (\hat{T}_q)^\vee.
    $$
    Here, $\hat{T}^\vee := T$ denotes the inverse of $T\mapsto \hat{T}.$
\end{prop}
\begin{proof}
    The edge projector of $\hat{T}_q$ is given by 
    \begin{align}
    \dim(H)^{-3}\cdot
    \tikzmath{
    \begin{scope}
    \clip[rounded corners=5pt] (-1.1,-2) rectangle (1.1,2);
    \fill[\AColor] (-1.1,-2) rectangle (1.1,2);
    \end{scope}
    \draw (-.7, -2) --(-.7,1) arc (180:0:.3cm) -- (-.1,0);
    \draw(-.4,1.3) -- (-.4,2);
    \filldraw (-.4,1.3) circle (.05cm);
    \draw(.3,-1) -- (.3,-2);
    \filldraw (.3,-1) circle (.05cm);
    \draw (0,0) -- (0,-.7) arc (180:360:.3cm) -- (.6, 2);
    \roundNbox{fill=white}{(0,.1)}{.2}{0}{0}{${\scriptstyle\hat{ T}}$}
    \draw[thick] (-.4, 1.7) arc (90:450:.4);
    \draw [->] (0,1.4);
    \roundNbox{fill=white}{(-.4,1.7)}{.14}{0}{0}{}
    \roundNbox{fill=white}{(-.7,1)}{.14}{0}{0}{}
    \roundNbox{fill=white}{(-.1,1)}{.14}{0}{0}{}
    \draw[thick] (.3, -.7) arc (90:450:.4);
    \draw [->] (.7,-1.1);
    \roundNbox{fill=white}{(.3,-1.5)}{.14}{0}{0}{}
    \roundNbox{fill=white}{(.6,-.8)}{.14}{0}{0}{}
    \roundNbox{fill=white}{(0,-.8)}{.14}{0}{0}{}
    \draw[thick] (0,.55) arc (90:450:.45);
    \draw [->] (.45,.1);
    \roundNbox{fill=white}{(0,-.4)}{.14}{0}{0}{}
    \roundNbox{fill=white}{(-.1,.55)}{.14}{0}{0}{}
    } 
    &=
    \dim(H)^{-2}\cdot
    \tikzmath{
    \begin{scope}
    \clip[rounded corners=5pt] (-1.1,-2) rectangle (1.7,2);
    \fill[\AColor] (-1.1,-2) rectangle (1.7,2);
    \end{scope}
    \draw (-.7, -2) --(-.7,1) arc (180:0:.3cm) -- (-.1,0);
    \draw(-.4,1.3) -- (-.4,2);
    \filldraw (-.4,1.3) circle (.05cm);
    \draw(.3,-.9) -- (.3,-2);
    \filldraw (.3,-.9) circle (.05cm);
    \draw (0,0) -- (0,-.6) arc (180:360:.3cm) -- (.6, 2);
    \roundNbox{fill=white}{(0,-.1)}{.2}{0}{0}{${\scriptstyle\hat{T}}$}
    \draw[thick] (.4,.8) arc (90:450:1);
    \draw [->] (1.4,.1);
    \roundNbox{fill=white}{(.3,-1.2)}{.14}{0}{0}{}
    \roundNbox{fill=white}{(.6,.8)}{.14}{0}{0}{}
    \draw[thick] (-.4, 1.5) arc (90:450:.5);
    \draw [->] (0.1,1.2);
    \roundNbox{fill=white}{(-.4,1.58)}{.14}{0}{0}{}
    \roundNbox{fill=white}{(-.7,.6)}{.14}{0}{0}{}
    \roundNbox{fill=white}{(-.15,.65)}{.25}{0}{0}{${\scriptstyle \bbX}$}
    }\nonumber \\ 
    &=
    \dim(H)^{-2}\cdot
    \tikzmath{
    \begin{scope}
    \clip[rounded corners=5pt] (-2,-2) rectangle (1.7,2);
    \fill[\AColor] (-2,-2) rectangle (1.7,2);
    \end{scope}
    \draw (-1.5, -2) --(-1.5,.5) arc (180:0:.5cm) arc (200:360:.3)  --(-.1,2);
    \filldraw (-.1,.3) circle (.05cm);
    \draw(-.1,.3) -- (0,.1);
    \draw(.3,-.9) -- (.3,-2);
    \filldraw (.3,-.9) circle (.05cm);
    \draw (0,0) -- (0,-.6) arc (180:360:.3cm) -- (.6, 2);
    \roundNbox{fill=white}{(0,-.1)}{.2}{0}{0}{${\scriptstyle\hat{T}}$}
    \draw[thick] (.4,.8) arc (90:450:1);
    \draw [->] (1.4,.1);
    \roundNbox{fill=white}{(.3,-1.2)}{.14}{0}{0}{}
    \roundNbox{fill=white}{(.6,.8)}{.14}{0}{0}{}
    \draw[thick] (-.9,1) -- (0,1.7) ..controls (0.3,1.4) .. (-.55,.75) .. controls (-1,.6) .. (-.9,1);
    \draw[->](.2, 1.5);
     \roundNbox{fill=white}{(-.5,.35)}{.14}{0}{0}{}
     \roundNbox{fill=white}{(-.55,.75)}{.14}{0}{0}{}
     \roundNbox{fill=white}{(-.9,1)}{.14}{0}{0}{}
     \roundNbox{fill=white}{(.1,.7)}{.14}{0}{0}{}
     \roundNbox{fill=white}{(0,1.2)}{.14}{0}{0}{}
     \roundNbox{fill=white}{(-.05,1.7)}{.14}{0}{0}{}
    }
    =T_q,
    \end{align}
    is the outcome of bubbling $T.$
    The second equality above follows from the compatibility of $\bbX = P\bullet P^*$ with comultiplication. The last equality is obtained  from the Duality Equations (\ref{eqn:DualityQAutG}) for $P$ and $P^*.$ This completes the proof. 
\end{proof}

Interestingly, the diagrams $\blacktriangle$ from Definition \ref{defn:CommonNeighbors} and $\triangle$ from Definition \ref{defn:WhiteTriangle} involved in Definition \ref{defn:3ptReg} transform nicely under bubbling.
\begin{prop}\label{prop:BubblingTriangles}
    Let $X$ be a finite quantum set and and let $\hat{R}, \hat{S}, \hat{T}$ be undirected real quantum graphs over $X$. Let $\bbX\in \QAut(\hat{R})\cap\QAut(\hat{S})\cap\QAut(\hat{T})$ be a quantum graph automorphism with splitting $P\bullet P^* = \bbX$ as in Construction \ref{const:Bubbling}. (We are mainly interested in the case where $\{\hat{R},\hat{S},\hat{T}\} \subseteq \{\hat{T},\hat{T}^c,\hat{J},\id_X\}$.)
    We then have that 
    \begin{enumerate}
        \item $(\blacktriangle_{\hat{T}})_q = \blacktriangle_{\hat{T}_q}$,
        \item $\left({}_{\hat{R}}\triangle^{\hat{T}}_{\hat{S}}\right)_q = {}_{\hat{R}_q}\triangle^{\hat{T}_q}_{\hat{S}_q}$.
    \end{enumerate}
\end{prop}
\begin{proof}
    1.- follows easily from a diagrammatic computation and the splitting of the quantum graph automorphism $\bbX = P\bullet P^*.$

    2.- We have that
    \begin{align}
        {}_{\hat{R}_q}\triangle^{\hat{T}_q}_{\hat{S}_q}
        =
        \dim(H)^{-6}\cdot
        \tikzmath{
        \begin{scope}
        \clip[rounded corners=5pt] (-1.6,-2.5) rectangle (1.7,3.8);
        \fill[\AColor] (-1.6,-2.5) rectangle (1.7,3.8);
        \end{scope}
        \draw[thick] (.5,3.3) arc (90:450:.4);
        \draw[thick] (-.4,1.7) -- (-.4,2) arc (180:0:.4) -- (.4,1.7) arc (360:180:.4); 
        \draw [->] (.4,1.7);
        \draw[thick] (.6,.5) -- (.6,1) arc (180:0:.4) -- (1.4,.5) arc (360:180:.4); 
        \draw [->] (1.4, .5);
        \draw[thick] (-.9,-.6) -- (-.9,-.5) arc (180:0:.4) -- (-.1,-.6) arc (360:180:.4); 
        \draw (0, 1.3) -- (0,2.4) arc (180:0:.5cm) -- (1,-.8);
        \draw (.5,2.9) -- (0.5,3.8);
        \filldraw (.5,2.9) circle (.05cm);
        \roundNbox{fill=white}{(0, 1.85)}{.25}{0}{0}{${\scriptstyle\hat{T}}$}
        \draw (-1,3.8) -- (-1,1.35) arc (180:360:.5);
        \draw (-.5, 0) -- (-.5,.85);
        \filldraw (-.5,.85) circle (.05cm);
        \draw (-.5,-.5) -- (-.5,0);
        \roundNbox{fill=white}{(-.5,-.5)}{.2}{0}{0}{${\scriptstyle \hat{R}}$}
        \roundNbox{fill=white}{(1,.7)}{.2}{0}{0}{${\scriptstyle\hat{S}}$}
        \draw (-.5,-.7) arc (180:360:.75cm);
        \draw (.25,-1.4) -- (.25,-2.5);
        \filldraw (.25,-1.45) circle (.05cm);
        \roundNbox{fill=white}{(0, 2.4)}{.15}{0}{0}{${}$}
        \roundNbox{fill=white}{(0, 1.3)}{.15}{0}{0}{${}$}
        \roundNbox{fill=white}{(.5, 3.3)}{.15}{0}{0}{${}$}
        \roundNbox{fill=white}{(1, .1)}{.15}{0}{0}{${}$}
        \roundNbox{fill=white}{(1, 1.3)}{.15}{0}{0}{${}$}
        \draw[->](-.1, -.45);
        \roundNbox{fill=white}{(-.5, -1)}{.15}{0}{0}{${}$}
        \roundNbox{fill=white}{(-.5, -.05)}{.15}{0}{0}{${}$}
        \roundNbox{fill=white}{(.8, 2.7)}{.15}{0}{0}{${}$}
        \draw[->] (-.15,.65) ;
        \roundNbox{fill=white}{(.2, 2.7)}{.15}{0}{0}{${}$}
        \draw[<-](.9,3);
        \draw[thick] (.25,-1.1) arc (90:450:.4);
        \roundNbox{fill=white}{(.25, -1.9)}{.15}{0}{0}{${}$}
        \roundNbox{fill=white}{(.6, -1.4)}{.15}{0}{0}{${}$}
        \roundNbox{fill=white}{(-.1, -1.4)}{.15}{0}{0}{${}$}
        \draw[<-](.6,-1.65);
        \draw[thick] (-.5,1.2) arc (90:450:.4); 
        \roundNbox{fill=white}{(-.85, .9)}{.15}{0}{0}{${}$}
        \roundNbox{fill=white}{(-.15, .9)}{.15}{0}{0}{${}$}
        \roundNbox{fill=white}{(-.5, .4)}{.15}{0}{0}{${}$}
        }
        &\overset{\clubsuit}{=}
        \dim(H)^{-2}\cdot
        \tikzmath{
        \begin{scope}
        \clip[rounded corners=5pt] (-1.2,-2.3) rectangle (1.2,2);
        \fill[\AColor] (-1.2,-2.3) rectangle (1.2,2);
        \end{scope}
        \draw[thick] (-1, 0) .. controls (-1,1) .. (-.8,1.2) ..controls(-.5,1.6) .. (-.1,1.7) --(0,1.7);
        \draw[thick] (-1,0) .. controls (-.9,-.5)..(0,-.5);
        \draw[->](1,0);
        \draw[thick] (1, 0) .. controls (1,1) .. (.8,1.2) ..controls(.5,1.6) .. (.1,1.7) --(0,1.7);
        \draw[thick] (1,0) .. controls (.9,-.5)..(0,-.5);
        \draw(-.5,-.5) -- (-.5,2);
        \draw(.5,-.5) -- (.5,2);
        \draw (-.5,-.5)-- (-.5,-1) arc (180:360:.5) --(.5,-.5);
        \draw (0,-2.3)--(0,-1.5);
        \draw[thick] (-.5,-1.5) .. controls(-.4,-1.1) .. (0,-1);
        \draw[thick] (.5,-1.5) .. controls(.4,-1.1) .. (0,-1);
        \draw[thick] (-.5,-1.5) .. controls(-.5,-1.7) .. (0,-1.8);
        \draw[thick] (.5,-1.5) .. controls(.5,-1.7) .. (0,-1.8);
        \draw[->] (.5,-1.5);
        \roundNbox{fill=white}{(-.5, 0)}{.25}{0}{0}{${\scriptstyle\hat{R}}$}
        \roundNbox{fill=white}{(.5, 0)}{.25}{0}{0}{${\scriptstyle\hat{S}}$}
        \roundNbox{fill=white}{(0, 1)}{.3}{.25}{.25}{$T$}
        \filldraw (0,-1.5) circle (.05cm);
        \roundNbox{fill=white}{(-.4, -1.2)}{.15}{0}{0}{${}$}
        \roundNbox{fill=white}{(.4, -1.2)}{.15}{0}{0}{${}$}
        \roundNbox{fill=white}{(-0, -1.8)}{.15}{0}{0}{${}$}
        \roundNbox{fill=white}{(-.5, 1.6)}{.15}{0}{0}{${}$}
        \roundNbox{fill=white}{(-.5, -.53)}{.15}{0}{0}{${}$}
        \roundNbox{fill=white}{(.5, 1.6)}{.15}{0}{0}{${}$}
        \roundNbox{fill=white}{(.5, -.53)}{.15}{0}{0}{${}$}
        }\\
        &\overset{\heartsuit}{=}
        \dim(H)^{-1}\cdot
        \tikzmath{
        \begin{scope}
        \clip[rounded corners=5pt] (-1.2,-2.3) rectangle (1.2,2);
        \fill[\AColor] (-1.2,-2.3) rectangle (1.2,2);
        \end{scope}
        \draw[thick] (-1, 0) .. controls (-1,1) .. (-.8,1.2) ..controls(-.5,1.6) .. (-.1,1.7) --(0,1.7);
        \draw[thick] (-1,0) .. controls (-.9,-1.7)..(0,-1.8);
        \draw[->](1,0);
        \draw[thick] (1, 0) .. controls (1,1) .. (.8,1.2) ..controls(.5,1.6) .. (.1,1.7) --(0,1.7);
        \draw[thick] (1,0) .. controls (.9,-1.7)..(0,-1.8);
        \draw(-.5,-.5) -- (-.5,2);
        \draw(.5,-.5) -- (.5,2);
        \draw (-.5,-.5)-- (-.5,-1) arc (180:360:.5) --(.5,-.5);
        \draw (0,-2.3)--(0,-1.5);
        \roundNbox{fill=white}{(-.5, 0)}{.25}{0}{0}{${\scriptstyle\hat{R}}$}
        \roundNbox{fill=white}{(.5, 0)}{.25}{0}{0}{${\scriptstyle\hat{S}}$}
        \roundNbox{fill=white}{(0, 1)}{.3}{.25}{.25}{$T$}
        \filldraw (0,-1.5) circle (.05cm);
        \roundNbox{fill=white}{(-0, -1.8)}{.15}{0}{0}{${}$}
        \roundNbox{fill=white}{(-.5, 1.6)}{.15}{0}{0}{${}$}
        \roundNbox{fill=white}{(.5, 1.6)}{.15}{0}{0}{${}$}
        }
        \overset{\spadesuit}{=}
        \left({}_{\hat{R}}\triangle^{\hat{T}}_{\hat{S}}\right)_q.\nonumber
    \end{align}
    Here, equality $\clubsuit$ is an application of Proposition \ref{prop:BubbleEdgeProjector} as well and commuting $\hat{R}$ and $\hat{S}$ past $\bbX = P\bullet P^*$. Equality $\heartsuit$ follows from Equation (\ref{eqn:QuantumFunction}) and finally equality $\spadesuit$ is Definition \ref{defn:WhiteTriangle}.
\end{proof}

We are now in position to state and prove our first main result, where we establish that many other relevant features are preserved by bubbling. 
\begin{thm}\label{thm:QuInvariants}{\bf [Theorem \ref{thmalpha:QuInvariants}]}
    Let $\cG = (\ell^2(X), \hat{T})$ be a quantum graph. 
    Consider a C*-Frobenius algebra quantum graph isomorphism $\mathbb{X}\in\mathsf{QAut}(\cG)$ of the form 
    $$
    \bbX: (H\otimes \overline{H})\otimes \ell^2(X)\to \ell^2(X)\otimes (H\otimes \overline{H}),   
    $$
    where $H$ is a finite-dimensional Hilbert space.
    Assume 
    $$
    P: H\otimes \ell^2(X_q) \to \ell^2(X)\otimes H
    $$
    is a quantum isomorphism between the quantum spaces $\ell^2(X)$ and $\ell^2(X_q)$ satisfying $P\bullet P^* \cong \mathbb{X}$ as in Construction \ref{const:Bubbling}.  
    Then, bubbling $\cG$ into the quantum graph $\cG_q := (\ell^2(X_q), \hat{T_q})$ satisfies:
    \begin{enumerate}[label={\bf (B\arabic*)}]
        \item\label{item:AlgRels} If $f(z)=\sum_{k=0}^na_kz^k\in\mathbb{C}[z]$ is a polynomial, then 
        $$
        f\left(\hat{T}_q\right)= \left(f\left(\hat{T}\right)\right)_q.
        $$
        That is, polynomials in $\hat{T}$ commute with bubbling. In particular $\hat{T}_q$ satisfies any polynomial equation satisfied by $\hat{T}$ and vice versa.  Similarly for $L(\hat{T})$. 
        \item\label{item:TopSpecFeatures} The Laplacian (Definition \ref{defn: laplacian}) and Degree Matrix (Definition \ref{defn: Degree matrix}) transform as 
        \begin{align*}
            L\left(\hat{T}_q\right) = L\left(\hat{T}\right)_q\qquad \text{ and }\qquad\mathsf{Deg}\left(\hat{T}_q\right) = \mathsf{Deg}\left(\hat{T}\right)_q\!.
        \end{align*}
        Moreover, for $S\in \{\hat{T}, L(\hat{T})\}$
        \begin{align*}
            \mathsf{Spec}(S) = \mathsf{Spec}(S_q)
        \end{align*}
        counting multiplicities, and 
        \begin{align*}
            \mathsf{Eigenspace}(S)_\lambda \cong \mathsf{Eigenspace}(S_q)_\lambda\qquad \forall \lambda\in \mathsf{Spec}(S) = \mathsf{Spec}(S_q).
        \end{align*}
        Consequently, the number of edges (Equation (\ref{eqn:NumQuEdges})), number of connected components, existence of cycles (Definition \ref{defn:pi0}) are preserved by bubbling. 
        \item\label{item:HighReg} For $k\in\{1,2,3\}$, $\hat{T}$ is $t$-point regular (Definition \ref{defn:3ptReg}) if and only if $\hat{T}_q$ is $t$-point regular, moreover  with the same parameters. Moreover, girth is preserved by bubbling and consequently being triangle-free is also preserved.
        \item\label {item:BubblingQAUT} The quantum automorphism group is preserved up to monoidal equivalence. 
        That is
        $$
        \Aut^+(\cG) \overset{\otimes}{\cong}\Aut^+(\cG_q).
        $$
    \end{enumerate}
\end{thm}
\begin{proof}
    
    \ref{item:AlgRels}: Is easy to check on monomials since $\hat{T}$ commutes with $\mathbb{X}$. The rest follows by linearity.

    \smallskip

    \ref{item:TopSpecFeatures}:  
    That bubbling transforms the Laplacian and Degree operators as stated is an easy check. 
    In Corollary \ref{cor:BubblingEigenspaces} we established that bubbling commutes with taking spectra of $S$.
    In particular, the nullity, rank and trace of the Laplacian are preserved and hence the number of connected components, the existence of cycles and the number of edges are invariants. 

    \smallskip

    \ref{item:HighReg}:    If $t=1,2$, then $t$-point regularity is defined by polynomial relations and so in these cases \ref{item:HighReg} follows from \ref{item:AlgRels}. 
    
    To show that $3$-point regularity and its parameters are preserved by bubbling follows from repeated application of Proposition \ref{prop:BubblingTriangles} to the cases where $\hat{R}, \hat{S}, \hat{T}\in \{\hat{T}, \hat{T}^c, \hat{J}, \id\}$, where it is readily seen that bubbling preserves de defining equation for 3-point regularity in Definition \ref{defn:3ptReg}. 

    \smallskip
    
    \ref{item:BubblingQAUT}:  
    The quantum graph isomorphism $\cG\cong_q \cG_q$ yields a finite-dimensional $*$-representation of the \textit{Linking Algebra} (see \cite{MR4091496} for definitions) between these two graphs. 
    This can be verified directly, using the unitarity of $P$ along with properties of quantum functions (c.f. \cite[Remarks 2.31 and 2.39]{MR4481115})
    Therefore, the Linking Algebra is non-zero, and by  \cite[Theorem 4.7]{MR4091496}, we conclude the existence of the desired monoidal equivalence. 
\end{proof}

\begin{remark}[{\bf Property (T) discrete quantum groups}]\label{remark:Prop(T)}
    We now comment on Items \ref{item:BubblingQAUT} and \ref{item:(T)} regarding {\bf approximation properties} of (the duals of) quantum automorphism groups of certain quantum graphs. 
    As explained in the proof of Theorem \ref{thm:QuInvariants} Item \ref{item:BubblingQAUT}, we have an equivalence of unitary tensor categories 
    $ \Rep(\Aut^+(\cG))\cong \Rep(\Aut^+(\cG_q)).$
    
    In the case of the Higman-Sims graph $HS$, it is known that $\Aut^+(HS)\overset{\otimes}{\cong} SO_q(5),$ with $q$ the square of the golden ratio \cite{MR1469634}, which is a compact quantum group of Kac-type, whose dual is an infinite discrete group with Property (T) \cite{MR3849625}. 
    Consequently, the dual of $\Aut^+(HS_q)$ is an infinite discrete quantum group with Property (T). 
\end{remark}

Before closing this section, we shall give some examples of quantum isomorphisms that will be important in forthcoming sections in the context of finite groups, Cayley graphs, and \emph{nice} unitary error bases. 
A {\bf unitary error basis} (UEB) in $\bbC^N$ consists of a collection of unitary matrices $\{V_i\}$ which moreover forms an orthonormal basis for $\End(\bbC^N)$ with respect to the Hilbert-Schmidt inner product given by the unnormalized trace.

\begin{defn}\label{defn:CentralType}\cite[Definition 7.12.21]{MR3242743}
    Let $\Gamma$ be a finite group. 
    Let $L\leq \Gamma$ and $\psi:L \times L \to U(1)$ be a $2$-cocycle. 
    We sat that $\psi$ is \emph{non-degenerate} if the twisted group algebra satisfies 
    $$
    \bbC^\psi[L]\cong \End(\cH),
    $$
    for some Hilbert space $\cH.$
    That is, if $\bbC^\psi[L]$ is simple. 
    Recall the $*$-algebra  structure on $\bbC^\psi[L]$ (c.f. \cite[Equation (65)]{MR3907958}) is given by
    \begin{align*}
    \ell\star_\psi\ell' := \dfrac{1}{\sqrt{|L|}}\psi(\ell, \ell')\ell\ell'\qquad \text{and}\qquad e_\psi:= \sqrt{|L|}\overline{\psi}(e,e)e.
    \end{align*}

    We say that $L$ is a {\bf subgroup of central type} if it admits some non-degenerate cocycle. 
\end{defn}

\begin{defn}\label{defn:niceUEB}
A {\bf nice unitary error basis} (c.f. \cite[Definition 4.5]{MR3907958}) for a group of central type $L$ is a collection of unitaries $\{U_\ell\}_{\ell\in L}\subset \End(\cH),$ where $\cH$ is a Hilbert space with $\dim(\cH)^2 = |L|$, satisfying for all $\ell, \ell':$
\begin{enumerate}
    \item $\{U_\ell\}_{\ell\in L}\subset \cU(\cH)$ is an orthonormal basis; i.e.  $\Tr(U_\ell^\dag U_{\ell'})= \dim(\cH)\delta_{\ell=\ell'}$,
    \item $U_{\ell'}U_\ell = \psi(\ell',\ell)U_{\ell' \ell}$.
\end{enumerate}
Here, $\psi$ is a non-degenerate $2$-cocycle on $L$. 
Without loss of generality we shall assume that $\psi(\ell, e) = 1 = \psi(e, \ell)$ for all $\ell\in L$, and thus $U_e= \id_\cH.$
\end{defn}

\begin{ex}\label{ex:niceUEB}
There always exist nice UEBs for all $N$, for example $\{V_{r,s}\} = \{X^rZ^s\}_{r,s=1}^N,$ where $X$ and $Z$ are generalized Pauli matrices: letting $\omega = \exp(2\pi \sqrt{-1}/N)$, $X|k\rangle = \omega^k|k\rangle$ and $Z|k \rangle= |k+1\rangle$ $(\mathsf{mod}\ N).$ 

Any UEB $\{V_{rs}\}_{r,s=1}^N$ defines a quantum isomorphism  
\begin{align}\label{ex:UEB}
    U': \mathbb{C}^{N^2}&\to M_N\otimes \End(\mathbb{C}^N)\\
        \delta_{\alpha\beta}&\mapsto \sum_{r,s} e_{rs} \otimes V_{\alpha\beta}e_{rs}V^*_{\alpha\beta},\nonumber
\end{align}
where $\{e_{\alpha\beta}\}$ is a system of matrix units of $M_N.$
\end{ex}

\section{Applications}\label{sec:Applications}
\subsection{Deforming 3-point regular classical graphs}\label{sec:DeformSpin}
We shall use the bubbling procedure described in Construction \ref{const:Bubbling} to deform the classical 9-Paley, 16-Clebsch and Higman-Sims graphs into 3-point regular connected quantum graphs. 
Similarly, we also obtain non-commutative analogues of the Shl\"afli and the Shrikhande graphs. 
A key feature that allows bubbling in these cases is that these are Cayley graphs of generally non-abelian groups that contain a copy of some $\bbZ_N\times \bbZ_N$, and so we can make use of nice unitary error bases (c.f. Example \ref{ex:niceUEB}) to construct a quantum automorphism with the structure of a C*-Frobenius algebra that we can split. 

In the setting of Cayley graphs over a finite abelian group $\Gamma\cong \bbZ_n\times \bbZ_n$, the bubbling procedure specializes to the twisted Cayley construction introduced by Gromada on \cite[\S6]{MR4514486}, which is a special case of \cite[\S3]{MR3907958}. 
Indeed, let $\cG=\mathrm{Cay}(\Gamma,S)$ and let $\mathbb{X}$ be the quantum automorphism arising from the left regular representation of $\Gamma$ equipped with a non-degenerate $2$-cocycle $\sigma$, such that $\bbC[\Gamma]^\sigma\cong \End(H)$ for some finite-dimensional Hilbert space. 
Therefore, the associated quantum isomorphism $\bbX$ (Equation (\ref{ex:UEB})) carries the structure of a simple $C^*$-Frobenius algebra in $\mathrm{QAut}(\cG)$. 
Applying the bubbling construction amounts to forming the orthogonal projection (Equation (\ref{eqn:pi}))
\[
\pi=\frac{1}{|\Gamma|}\sum_{\alpha\in\Gamma}\overline{U_\alpha}\otimes (\alpha\rhd-) \otimes U_\alpha
\]
on $H^*\otimes \ell^2(\Gamma)\otimes H$, where $(\alpha\rhd-)$ denotes translation by $\alpha$. Passing to the Fourier basis of $\ell^2(\Gamma)$, the operator $\pi$ diagonalizes and its range $\ell^2(X_q)=\mathrm{Ran}(\pi)$ is canonically identified with the twisted group algebra $\bbC[\Gamma]^\sigma$.
Under this identification, the multiplication induced by the splitting $\mathbb{X} = P\bullet P^*$ coincides with the twisted convolution product, while the bubbled adjacency operator $\widehat T_q$ acts diagonally with the same eigenvalues as the original Cayley adjacency operator. Consequently, bubbling along $(\Gamma, \sigma)$ reproduces exactly the twisted Cayley graphs of Gromada. 

\begin{construction}[{\bf A quantum 9-Paley graph}]\label{const:q9Paley}
It is a well-known fact that the {$9$-Paley} graph is realized as a Cayley graph over $\Gamma = \bbZ_3\times \bbZ_3$ with respect to the canonical generators $\varepsilon_1, \varepsilon_2$. We now closely follow the cocycle twist procedure on $\Gamma$, with $n = 3$, in order to find a suitable quantum group for constructing a non-classical analogue of the 9-Paley, denoted $9P_q$. We shall reproduce some parts of \cite{MR4514486} here for the reader's convenience. 

Let $\tau_1, \tau_2$ be the corresponding generators of $C(\Gamma) \simeq \mathbb{C}\Gamma$, where
$$
C(\Gamma) = C^*(\tau_1, \tau_2 \mid \tau_1^n = \tau_2^n = \tau_1 \tau_1^* = \tau_1^* \tau_1 = \tau_2 \tau_2^* = \tau_2^* \tau_2 = 1,\ \tau_1 \tau_2 = \tau_2 \tau_1).
$$

Then, we take as unitary bicharacter $\sigma: \Gamma\times \Gamma \to \mathbb{C}$ given by 
$$\sigma(\tau_1, \tau_2) = \sigma(\tau_1, \tau_1) = \sigma(\tau_2, \tau_2) = 1,$$
$$\sigma(\tau_2, \tau_1) = \omega,$$
where $\omega = e^{2\pi i/n}$ is a primitive $n$-th root of unity. 

Thus, applying the corresponding cocycle twist, as defined on \cite[Definition 5.6]{MR4514486}, we obtain
$$
C(\check{\Gamma}) = \check{C}^*(\check{\tau}_1, \check{\tau}_2 \mid \check{\tau}_1^n = \check{\tau}_2^n = \check{\tau}_1 \check{\tau}_1^* = \check{\tau}_1^* \check{\tau}_1 = \check{\tau}_2 \check{\tau}_2^* = \check{\tau}_2^* \check{\tau}_2 = 1,\ \check{\tau}_1 \check{\tau}_2 = \omega \check{\tau}_2 \check{\tau}_1).
$$
We have to verify that the space obtained after the twisting process is effectively $M_3(\mathbb{C})$; to that end, consider the map $\varphi: C(\check{\Gamma})\to M_3(\mathbb{C})$, defined as
$$\varphi(\check{\tau_1}) = \begin{pmatrix}
    1 & 0 & 0\\
    0 & \omega & 0\\
    0 & 0 & \omega^2\\
\end{pmatrix}, \qquad \varphi(\check{\tau_2}) = \begin{pmatrix}
    0 & 1 & 0\\
    0 & 0 & 1\\
    1 & 0 & 0\\
\end{pmatrix},$$
with its corresponding inverse 
$$\varphi^{-1}(e_{ij}) = \frac{1}{3} \sum_k \omega^{-ik} \tau_1^k \tau_2^{i-j}. $$
An straightforward computation shows that this maps are each other's inverses, and the above matrices satisfy the relations on $C(\check{\Gamma})$. 

We now focus on the Cayley graph $\mathrm{Cay}(\Gamma, S)$ with $S = \{\varepsilon_1, \varepsilon_2\}$ in order to a quantum analogue of this Cayley graph, over the quantum space $\check{\Gamma} = {M}_n(\mathbb{C})$. 
We shall  use of the following result which provides a formula for the eigenvalues of the adjacency matrix of a Cayley graph. 
\begin{thm}[\cite{Lovsz1975}]\label{thm: Formula eig cayley}
    Let $\Gamma$ be a finite abelian group and $S \subseteq \Gamma$. Denote by $\hat{T}$ the adjacency matrix associated to the Cayley graph $\mathrm{Cay}(\Gamma, S)$. 
    Then $\mathrm{Irr}\,\Gamma$ forms the eigenbasis of $\hat{T}$. 
    Given $\mu \in \Gamma$, the eigenvalue corresponding to $\tau_\mu$ is given by
$$
\lambda_\mu = \sum_{\vartheta \in S} \tau_\mu(-\vartheta).
$$
Moreover, the change of basis operator is given by the Fourier transform $\mathcal{F}: C(\Gamma) \to \mathbb{C}\Gamma$. 
\end{thm}

Using Theorem \ref{thm: Formula eig cayley}, we can compute the spectrum of the adjacency matrix of $\mathrm{Cay}(\check{\Gamma}, \check{S})$, and thus, get the diagonalization of the adjacency operator:
$$
\tilde{T}_{ab}^{cd} = \delta_{ac} \delta_{bd} \lambda_{ab},\qquad \lambda_{ab} = \sum_{i=1}^{n-1} \omega^{ia} + \sum_{j=1}^{n-1} \omega^{jb} = 3\delta_{a,0} + 3\delta_{b,0} - 2.
$$

By twisting from $\mathbb{C}^9$ into $M_3(\mathbb{C})$, the adjacency matrix does not change, so $\check{T} = \hat{T}$. Finally, after a change of basis from the Fourier basis into $(e_{ij})$, we get:
$$
\hat{T}_{ij}^{kl} = \sum_{a,b} \varphi_{ij}^{ab} [\varphi^{-1}]_{ab}^{kl} \lambda_{ab} = \frac{1}{3} \sum_{a,b} \delta_{b, i-j} \delta_{b, k-l} \omega^{a(i-k)} (3\delta_{a,0} + 3\delta_{b,0} - 2)
= \delta_{i-j, k-l} + 3 \delta_{ijkl} - 2 \delta_{ik} \delta_{jl}.
$$

We can then explicitly calculate the adjacency matrix 
$$\hat{T} = \begin{pmatrix}
 2 & 0 & 0 & 0 & 1 & 0 & 0 & 0 & 1 \\
 0 & -1 & 0 & 0 & 0 & 1 & 1 & 0 & 0 \\
 0 & 0 & -1 & 1 & 0 & 0 & 0 & 1 & 0 \\
 0 & 0 & 1 & -1 & 0 & 0 & 0 & 1 & 0 \\
 1 & 0 & 0 & 0 & 2 & 0 & 0 & 0 & 1 \\
 0 & 1 & 0 & 0 & 0 & -1 & 1 & 0 & 0 \\
 0 & 1 & 0 & 0 & 0 & 1 & -1 & 0 & 0 \\
 0 & 0 & 1 & 1 & 0 & 0 & 0 & -1 & 0 \\
 1 & 0 & 0 & 0 & 1 & 0 & 0 & 0 & 2
\end{pmatrix},
$$
with spectrum 
$$\mathsf{Spec}(\hat{T}) = \{ 1,  4, -2, -2,  1,  1,  1, -2, -2 \}, \quad \mathsf{Spec}(L(\hat{T})) = \{  3, 0,  6,  3,  6,  3,  3,  6,  6\}.$$
It can be verified, either abstractly by by Theorem \ref{thm:QuInvariants} or by hand, using our Software, Section \href{https://colab.research.google.com/drive/11rctqnjD3bxMejfEKWnJfr2qUy6oUYA3?usp=sharing}{4.1}, that  $9P_q$ inherits the properties of its classical analogue $9P$. That is, it is an irreflexive,  connected, 3-point regular graph with the parameters from Example \ref{ex:M3TripleRegularity}, containing cycles. 

Notice that, our notation $\hat{T}$ represents the quantum adjacency operator on $M_3(\mathbb{C})$ after the cocycle twist in the Fourier basis. 
In contrast, Gromada in \cite[\S6]{MR4514486} uses $T$ for this same operator but in the standard $\{e_{ij}\}$ basis.
Similarly, Gromada uses $\hat{T}$ to denote this operator expressed in the Fourier basis, while for us it is the adjacency operator expressed in the standard $\{e_{ij}\}$ basis.
\end{construction}

\begin{construction}[{\bf A quantum 16-Clebsch  graph}]\label{const:q16Clebsch}
By repeating the cocycle twisting process from Example \ref{const:q9Paley}, we shall obtain a quantum analogue of the $16$-Clebsch graph over $M_4(\mathbb C)$, denoted $16CL_q$. 
We first realize the classical 16-Clebsh Graph as the Cayley Graph $\mathrm{Cay}(\Gamma, S)$, where $\Gamma = \mathbb Z_2^4$ and $S$ is the set of canonical generators $\{\varepsilon_i\}_{i=1}^4$ joint with the element $\sum_i \varepsilon_i$. 

Let $\tau_1,..., \tau_n$ the generators of $C(\Gamma)$, where 
$$C(\Gamma) = C^*(\tau_1, \ldots, \tau_n \mid \tau_i = \tau_i^*,\ \tau_i^2 = 1,\ \tau_i \tau_j = \tau_j \tau_i).
$$
Then, we define a cocycle $\sigma:\Gamma\times \Gamma\to \mathbb C$ as 
$$\sigma_{ij} := \sigma(\varepsilon_i, \varepsilon_j) = 
\begin{cases}
-1 & \text{if } i > j, \\
+1 & \text{if } i \leq j,
\end{cases}.$$

As a result, the twisted algebra $C(\check{\Gamma})$, obtain by using \cite[Definition 5.6]{MR4514486}, is exactly the Clifford algebra with 4 generators, defined as follows:
$$
\mathrm{Cl}_4 = C^*(\check{\tau}_1, \ldots, \check{\tau}_4 \mid \check{\tau}_i = \check{\tau}_i^*,\ \check{\tau}_i^2 = 1,\ \check{\tau}_i \check{\tau}_j = -\check{\tau}_j \check{\tau}_i) \simeq C(\check{\Gamma}).
$$
A standard result on classification of complex Clifford algebras gives  $\mathrm{Cl}_4 \simeq M_4(\mathbb C)$. Then the adjacency operator $\hat{T} \colon \mathrm{Cl}_4 \to \mathrm{Cl}_4$, defined by
$$
\hat{T} \check{\tau}_\mu = \lambda_\mu \check{\tau}_\mu \quad \text{with} \quad \lambda_\mu = 4-2\cdot \mathrm{deg}(\mu),
$$
where $\check{\tau}_\mu = \check{\tau}_{\mu_1}...\check{\tau}_{\mu_k}$ for $\mu = e_{\mu_1} + ... e_{\mu_k}\in \mathbb{Z}_2^n$, defines a quantum graph $\operatorname{Cay}(\check{\Gamma}, \check{S})$ over the quantum space $M_4(\mathbb C)$, which is quantum isomorphic to the classical Cayley graph $\operatorname{Cay}(\mathbb{Z}_2^4, S)$. 

It then remains to make the change of basis from the Fourier basis to the standard basis of $\mathrm{Cl}_4$. Take the basis for $\mathrm{Cl}_4$ as the standard representation, given by the tensor product of the $2\times 2$ Pauli matrices \cite{cartier1996algebraic}:
$$\check{\tau}_1 = \sigma_1 \otimes I_2, \quad \check{\tau}_2 = \sigma_3 \otimes \sigma_1, \quad \check{\tau}_3 = \sigma_3 \otimes \sigma_2\quad \check{\tau}_4 = \sigma_3 \otimes \sigma_3, $$
where
$$\sigma_1 = \begin{pmatrix}
    0 & 1\\
    1 & 0
\end{pmatrix}, \quad \sigma_2 = \begin{pmatrix}
    0 & -i\\
    i & 0
\end{pmatrix} \quad \sigma_3 = \begin{pmatrix}
    1 & 0\\
    0 & -1
\end{pmatrix}. $$
The set $\{ \prod_{i\in S}\check{\tau}_i: S\subseteq\{1,...,4\} \}$ is a basis for $\mathrm{Cl}_4$, and by applying the corresponding matrix change $P$ to the diagonal matrix $\hat{T}$, we get
{\small $$\hat{T} = P^{-1} \tilde{T} P = \left(\begin{array}{cccccccccccccccc}
 2 & 0 & 0 & 0 & 0 & 1 & 0 & 0 & 0 & 0 & 1 & 0 & 0 & 0 & 0 & 1 \\
 0 & 0 & 0 & 0 & 1 & 0 & 0 & 0 & 0 & 0 & 0 & 1 & 0 & 0 & -1 & 0 \\
 0 & 0 & 0 & 0 & 0 & 0 & 0 & 1 & 1 & 0 & 0 & 0 & 0 & -1 & 0 & 0 \\
 0 & 0 & 0 & -2 & 0 & 0 & 1 & 0 & 0 & 1 & 0 & 0 & 1 & 0 & 0 & 0 \\
 0 & 1 & 0 & 0 & 0 & 0 & 0 & 0 & 0 & 0 & 0 & -1 & 0 & 0 & 1 & 0 \\
 1 & 0 & 0 & 0 & 0 & 2 & 0 & 0 & 0 & 0 & 1 & 0 & 0 & 0 & 0 & 1 \\
 0 & 0 & 0 & 1 & 0 & 0 & -2 & 0 & 0 & 1 & 0 & 0 & 1 & 0 & 0 & 0 \\
 0 & 0 & 1 & 0 & 0 & 0 & 0 & 0 & -1 & 0 & 0 & 0 & 0 & 1 & 0 & 0 \\
 0 & 0 & 1 & 0 & 0 & 0 & 0 & -1 & 0 & 0 & 0 & 0 & 0 & 1 & 0 & 0 \\
 0 & 0 & 0 & 1 & 0 & 0 & 1 & 0 & 0 & -2 & 0 & 0 & 1 & 0 & 0 & 0 \\
 1 & 0 & 0 & 0 & 0 & 1 & 0 & 0 & 0 & 0 & 2 & 0 & 0 & 0 & 0 & 1 \\
 0 & 1 & 0 & 0 & -1 & 0 & 0 & 0 & 0 & 0 & 0 & 0 & 0 & 0 & 1 & 0 \\
 0 & 0 & 0 & 1 & 0 & 0 & 1 & 0 & 0 & 1 & 0 & 0 & -2 & 0 & 0 & 0 \\
 0 & 0 & -1 & 0 & 0 & 0 & 0 & 1 & 1 & 0 & 0 & 0 & 0 & 0 & 0 & 0 \\
 0 & -1 & 0 & 0 & 1 & 0 & 0 & 0 & 0 & 0 & 0 & 1 & 0 & 0 & 0 & 0 \\
 1 & 0 & 0 & 0 & 0 & 1 & 0 & 0 & 0 & 0 & 1 & 0 & 0 & 0 & 0 & 2
\end{array}\right),$$}
with spectrum 
$$
\mathsf{Spec}(\hat{T})=\{ 1,  5,  1, -3,  1, -3,  1,  1, -3, -3, -3,  1,  1,  1,  1,  1\}.
$$ 
It can directly be verified that $\hat{T}$ defines an irreflexive $3$-point regular connected quantum graph with parameters given in Example \ref{ex:params}, which is quantum isomorphic to the classical $16$-Clebsch graph. 
We do this explicitly in Section \href{https://colab.research.google.com/drive/11rctqnjD3bxMejfEKWnJfr2qUy6oUYA3?usp=sharing}{4.3} of the Software
Moreover, it has the same Laplacian spectrum and the same topology as its classical analogue:
$$\sigma(L(\hat{T})) = \{ 4, 0, 4, 8, 4, 8, 4, 4, 8, 8, 8, 4, 4, 4, 4, 4 \}.$$
Independently, by Theorem \ref{thm:qDeformations} can see all of these properties are shared with its classical analogue. 
\end{construction}

\begin{construction}\label{const:qHS}
We now describe a {\bf non-classical analogue of the Higman-Sims graph}. 
The classical Higman-Sims graph, denoted $HS$, is known to be a Cayley graph over the non-abelian group\footnote{See \url{http://linear.ups.edu/eagts/section-11.html} for a construction on Sage Math. In GAP this group is indexed as (100,9).}
$$
H = \bbZ_5\times (\bbZ_5\rtimes \bbZ_4)
$$
as described in \cite{Heinze2001, MR1975767}\footnote{This construction is also attributed to  C.E. Praeger and C. Schneider by these authors.}.  
The order of $H$ is  $|H|=100$, with $\bbC[H] \cong \bbC^{\oplus 20}\oplus M_4^{\oplus 5}$, having $25$ irreducible representations of dimensions $1$ and $4$ only. 
As a group of permutations, it can be presented as a subgroup of $S_{10}$, the symmetric group in ten points as follows:
\begin{align}\label{eqn:HasPermutation}
H &\cong  \Big\langle g_1:=(1\ 2\ 3\ 4\ 5),\ g_2:=(6\ 7\ 8\ 9\ 10),\ g_3:=(2\ 3\ 5\ 4)\Big\rangle \leq S_{10}.  
\end{align}
We now look for subgroups of central type of $H$ (c.f. Definition \ref{defn:CentralType}), which must be of order a square dividing $100,$ limiting the possibilities to $2^2$ and $5^2$.
By a computational exploration on \emph{Sage Math}, all order four subgroups of $L$ are cyclic and the unique subgroup of order $25$ is  
$$
L:=\langle g_1,g_2\rangle \cong \bbZ_5\times\bbZ_5.
$$
The known classification of abelian groups of central type (c.f. \cite[Theorem 5]{MR1843319}, \cite{MR3520282}) tells us that these are exactly those of the form $\Gamma\times \Gamma$ for an abelian group $\Gamma$, and so, $L$ is the unique subgroup of $H$ of central type. 

We shall now give an explicit normalized 2-cocycle $\psi:L\times L \to U(1)$ so that the twisted group algebra $\bbC L^\psi$ is isomorphic to $\End(\cH)$ for some Hilbert space $\cH.$ 
The data of $\psi$ precisely matches that of a nice unitary error basis (c.f. Definition \ref{defn:niceUEB}). 
For $a,b,c,d\in\{0,...,4\}$ and a primitive $5$-th root of unity, let 
$$
\psi(ag_1+bg_2,\ cg_1+dg_2):= \omega^{bc},
$$
which by \cite[\S6.1]{MR4514486} is known to be a nice unitary error basis, yielding $\bbC L^\psi \cong \End(\bbC^5)\cong M_5.$ 
To describe this $*$-isomorphism explicitly we use a generalized \emph{Pauli unitary error basis} $\{X^iZ^j\}_{i,j=0}^4$, where $X|k\rangle = |k-1\rangle\ \mathsf{mod}(5)$ and $Z|k\rangle = \omega^k|k\rangle$ with $\{|k\rangle\}_0^5$ denoting the standard basis of $\cH=\bbC^5.$
The map is then given by $g_1^ig_2^j\mapsto X^iZ^j,$ and is readily seen to satisfy the properties of a nice unitary error basis. 
We have therefore witnessed that $\psi$ is a non-degenerate $2$-cocycle on $L,$ and so, $L$ is a subgroup of central type. 

We have gathered all the ingredients to construct a classical C*-Frobenius algebra in the category $\Hilb_{\Aut(HS)}\subsetneq \mathsf{QAut}(HS)$ in the sense of \cite{MR3907958}.
In fact, the purpose of realizing $HS$ as a Cayley graph over $H$ is that $L\leq H\leq \Aut(HS)$ where the action on the vertices is by left translation. 
Indeed, by \cite[Proposition 4.6]{MR3907958}, the map
\begin{align}
    \bbX_{L,\psi}:\bbC^5\otimes \overline{\bbC^5}\otimes V_{HS}\to V_{HS}\otimes \bbC^5\otimes \overline{\bbC^5}
\end{align}
given by 
\begin{align*}
    \bbX_{L,\psi}:= \dfrac{1}{5}\cdot\sum_{\ell\in L\leq\Aut(HS)}
    \tikzmath{
        \begin{scope}
        \clip[rounded corners=5pt] (-1,-2) rectangle (1,2);
        \fill[\AColor] (-1,-2) rectangle (1,2);
        \end{scope}
         \draw[thick] (-.8, -2) .. controls (-.7,.2) .. (-.3, -.6);
         \draw[thick] (.8, 2) .. controls (.7,-.2) .. (.3, .6);
         \draw[thick] (-.3,-.5) -- (-.3,-2);
         \draw[thick] (.3,.5) -- (.3,2);
         \draw[->] (-.78,-1.5);
         \draw[>-] (-.3,-1.5);
         \draw(0, -0) -- (.6, -2);
         \draw[->] (.3,1.5);   
         \draw[>-] (.78,1.5);
         \draw(-0,0) -- (-.6, 2);
        \roundNbox{fill=white}{(0,0)}{.28}{0}{0}{$\ell\rhd$}
        \draw (-.8,-1.5) node[below]{};
        \draw (-.5,-1.45) node[below]{};
        \draw (.8,1.5) node[above]{};
        \draw (.5,1.5) node[above]{};
        \draw (-.8,1.5) node[above]{};
        \draw (.8,-1.45) node[below]{};
        \roundNbox{fill=white}{(-.25,-.8)}{.3}{0}{0}{$U_\ell^*$}
        \roundNbox{fill=white}{(.25,.8)}{.3}{0}{0}{$U_\ell$}
    }.
\end{align*}
Following Equation (\ref{eqn:pi}),  rotating this map into 
$$
\pi = \dfrac{1}{25}\sum_{\ell\in L\leq\Aut(HS)} \overline{U_\ell}\otimes (\ell\rhd-)\otimes U_\ell
$$
yields a self-adjoint idempotent $\pi$ on $\overline{\bbC^5}\otimes V_{HS}\otimes \bbC^5$, whose $\dagger$-splitting yields an isometry
$$
\iota:A\to \overline{\bbC^5}\otimes V_{HS}\otimes \bbC^5$$
satisfying $\iota^\dag\circ\iota = \id_{X}$ as well as $\iota\circ\iota^\dag = \pi$. 
Here, $A$ is some Hilbert space.
Following the proof of \cite[Theorem 3.4]{MR3907958}, the linear map 
\begin{align}
    P:\bbC^5\otimes A\to V_{HS}\otimes \bbC^5
\end{align}
defined as in Equation (\ref{eqn:P}), which is a unitary by construction. 
Consequently,  the dimension of $A$ is $100.$
Using the algebra structure on $V_{HS}$ and $P$, we canonically turn $A$ into a tracial quantum space denoted $V_{HS_q}$ with structure maps as given in Equation (\ref{eqn:Xq}). 
Furthermore, we use $P$ to construct an adjacency operator $\hat{HS_q}$ on $\ell^2(V_{HS_q})$ as in Equation (\ref{eqn:BuubleTea}).
We denote the resulting {\bf quantum Higman-Sims graph} $(\ell^2(V_{HS_q}), \hat{HS_q})$ by $HS_q.$ 
Again, by \cite[Theorem 3.4]{MR3907958}, $P$ defines a quantum graph automorphism $$HS\cong_q HS_q.$$

We now verify that {\bf $HS_q$ is not classical}; i.e. that $V_{HS_q}$ is not abelian.
Notice that, if it were classical since $HS$ is strongly regular with parameters $(n=100,\ k=22,\ \lambda=0,\ \mu=6)$ (c.f. Example (\ref{ex:params})) then so would $HS_q$ by Theorem \ref{thm:QuInvariants}.
However, it is well-known that $HS$ is the unique strongly regular \emph{classical} graph with those parameters, implying $HS_q$ is classical and classically isomorphic to $HS$. 
To rule this possibility out, we use \cite[Equation (74)]{MR3907958}, computing the dimension of $Z(V_{HS_q}),$ the center of $V_{HS_q}$. 
We have that 
\begin{align}\label{eqn:DimCenter}
    \dim(Z(V_{HS_q})) = \dfrac{1}{|L|}\sum_{h\in H}\Phi^{L,\psi}_{\mathsf{Stab}_L(h)}.
\end{align}
Here, $\mathsf{Stab}_L(h) = \{\ell\in L|\ \ell\rhd h= h\}$ is the stabilizer of the vertex $h\in H$ of $HS$ inside of $L.$ 
Furthermore, for $\Lambda\leq L,$ as in \cite[Equations (69) and (72)]{MR3907958}
$$
\Phi^{L,\psi}_{\Lambda}:= \sum_{\alpha, \beta\in \Lambda,\ [\alpha, \beta]=e} \psi(\alpha,\beta)\cdot \overline{\psi}(\alpha\beta\alpha^{-1}, \alpha).
$$
Notice, however that a few simplifications appear. First, since $L$ is abelian, the commuting condition in the sum above as well as the conjugation disappear and we obtain
$$
\Phi^{L,\psi}_{\Lambda}= \sum_{\alpha, \beta\in \Lambda} \psi(\alpha,\beta)\overline{\psi}(\beta, \alpha).
$$
Furthermore, since $L\leq \Aut(HS)$ by left translation, for all $h\in HS$ we have $\mathsf{Stab}_L(h) = \{e\}$, and therefore 
$$
\Phi^{L,\psi}_{\mathsf{Stab}_L(h)} = 1\qquad \forall h\in HS.
$$
Equation (\ref{eqn:DimCenter}) then evaluates into 
\begin{align}\label{eqn:VHSqDecomp}
    \dim(Z(V_{HS_q})) = \dfrac{1}{|L|}\sum_{h\in H}\Phi^{L,\psi}_{\mathsf{Stab}_L(h)} = \dfrac{1}{25}\cdot 100=4\neq 100.
\end{align}
Therefore, $V_{HS_q}\not\cong \bbC^{100}$ is not abelian since it decomposes exactly into four matrix blocks 
\begin{align}
V_{HS_q}\cong M_{c_1}\oplus M_{c_2}\oplus M_{c_3}\oplus M_{c_4},\quad \text{ with }\quad c_1^2+c_2^2+c_3^2+c_4^2 = 100.    
\end{align}
\end{construction}

\begin{remark}
    We do not explicitly know the semisimple decomposition of $V_{KS_q}$ appearing in Equation (\ref{eqn:VHSqDecomp}), and we ignore if the Higman-Sims graph can be deformed over various C*-algebras, and perhaps over a factor.  
    Notice that by Equation (\ref{eqn:DimCenter}), if the Higman-Sims graph is a Cayley graph over, say, $\bbZ_{10}\times \bbZ_{10}$, then it could be realized over $M_{10}$. 
    However, we are not aware of any such presentation. 
\end{remark}

We summarize the conclusions of Constructions \ref{const:q9Paley}, \ref{const:q16Clebsch} and \ref{const:qHS} in the following
\begin{thm}[{\bf Theorem \ref{thmalpha:qDeformations}}]\label{thm:qDeformations}
    For each of the classical graph  $\cG\in\{9P,\ 16Cl,\ HS\}$, there exists a concrete quantum (non-classical) graph  $\cG_q\in\{9P_q,\ 16Cl_q,\ HS_q \}$ which is connected, $3$-point regular with the same parameters, Laplacian, minimal polynomial and a concrete quantum graph isomorphism 
    \begin{align}
        \cG\cong_q\cG_q.
    \end{align}
\end{thm}
\begin{proof}
    The existence of (non-classical) quantum graphs and the quantum isomorphisms was established in Constructions \ref{const:q9Paley}, \ref{const:q16Clebsch} and \ref{const:qHS}. 
\end{proof}

\begin{remark}
    There is no non-classic analogue of the pentagonal graph. \cite[Remark 4.3]{MR3907958}
\end{remark}

\begin{exs}\label{exs:OtherDeformations}
    This deformation construction by bubbling works on any graph known to be a Cayley graph over a finite group containing non-trivial subgroups of central type. 
    
    For instance, the {\bf Schl\"{a}fli Graph} is the unique classical strongly regular graph  with parameters $(27,16,10,8)$ \cite[Lemma 10.9.4]{MR1829620}, and is known to be a Cayley graph (see \cite{MR3654193} and references therein) for
    \begin{align*}
        \Gamma &= (\bbZ_3\times\bbZ_3)\rtimes\bbZ_3 \cong \left\langle a,b,c| a^3 = b^3 = c^3 = e,\ abc=ba,\ ac = ca,\ bc = cb\right\rangle,\\
        S&=\{a,\ a^2,\ b,\ b^2,\ c,\ c^2,\ cba,\ a^2b^2c^2,\ aba,\ bab\}. 
    \end{align*}
    Visibly, $\Lambda:=\left\langle a,b\right\rangle\leq \Gamma$ is a copy of $\bbZ_3\times\bbZ_3$ inside $\Gamma,$ which is a group of central type. 
    Therefore, there is quantum isomorphic a strongly regular quantum graph with the same parameters modeled over a non-abelian C*-algebra of dimension $27$ whose center has dimension $3$: either $M_3^{\oplus 3}$ or $M_5\oplus \bbC^{\oplus 2}.$

    Similarly, the {\bf Shrikhande graph} is a strongly regular classical graph with parameters $(16,6,2,2)$, which moreover is a Cayley graph \cite{MR3654193} for
    \begin{align*}
        \Gamma = \bbZ_4\times \bbZ_4\text{ with } S = \{\pm(1,0),\ \pm(0,1),\ \pm(1,-1)\}.
    \end{align*}
    The eigenvalues for the adjacency matrix of $\Gamma$ can be computed as 
    \[
    \lambda_{a,b} = 2 \, \mathrm{Re}(\omega^a) + 2 \, \mathrm{Re}(\omega^b) + 2 \, \mathrm{Re}(\omega^{a-b}), 
    \quad a,b = 0,1,\dots,n-1, \quad \omega = e^{2\pi i / n}.
    \]
    Thus, similarly as the case for the 9-Paley graph, its adjacency matrix is:
    $${ \hat{T} = \left(\begin{array}{cccccccccccccccc}
2 & 0 & 0 & 0 & 0 & 2 & 0 & 0 & 0 & 0 & 0 & 0 & 0 & 0 & 0 & 2 \\
0 & 0 & 0 & 0 & 0 & 0 & 1+i & 0 & 0 & 0 & 0 & 0 & 1-i & 0 & 0 & 0 \\
0 & 0 & -2 & 0 & 0 & 0 & 0 & 0 & 0 & 0 & 0 & 0 & 0 & 0 & 0 & 0 \\
0 & 0 & 0 & 0 & 1-i & 0 & 0 & 0 & 0 & 0 & 0 & 0 & 0 & 0 & 1+i & 0 \\
0 & 0 & 0 & 1+i & 0 & 0 & 0 & 0 & 0 & 1-i & 0 & 0 & 0 & 0 & 0 & 0 \\
2 & 0 & 0 & 0 & 0 & 2 & 0 & 0 & 0 & 0 & 2 & 0 & 0 & 0 & 0 & 0 \\
0 & 1-i & 0 & 0 & 0 & 0 & 0 & 0 & 0 & 0 & 0 & 1+i & 0 & 0 & 0 & 0 \\
0 & 0 & 0 & 0 & 0 & 0 & 0 & -2 & 0 & 0 & 0 & 0 & 0 & 0 & 0 & 0 \\
0 & 0 & 0 & 0 & 0 & 0 & 0 & 0 & -2 & 0 & 0 & 0 & 0 & 0 & 0 & 0 \\
0 & 0 & 0 & 0 & 1+i & 0 & 0 & 0 & 0 & 0 & 0 & 0 & 0 & 0 & 1-i & 0 \\
0 & 0 & 0 & 0 & 0 & 2 & 0 & 0 & 0 & 0 & 2 & 0 & 0 & 0 & 0 & 2 \\
0 & 0 & 0 & 0 & 0 & 0 & 1-i & 0 & 0 & 0 & 0 & 0 & 1+i & 0 & 0 & 0 \\
0 & 1+i & 0 & 0 & 0 & 0 & 0 & 0 & 0 & 0 & 0 & 1-i & 0 & 0 & 0 & 0 \\
0 & 0 & 0 & 0 & 0 & 0 & 0 & 0 & 0 & 0 & 0 & 0 & 0 & -2 & 0 & 0 \\
0 & 0 & 0 & 1-i & 0 & 0 & 0 & 0 & 0 & 1+i & 0 & 0 & 0 & 0 & 0 & 0 \\
2 & 0 & 0 & 0 & 0 & 0 & 0 & 0 & 0 & 0 & 2 & 0 & 0 & 0 & 0 & 2
\end{array}\right)}.
$$
Notice that $\hat{T} = -\hat{T}^*$ and $\hat{T}^\dag = \hat{T}$, is a non-real Schur-idempotent and thus defines an undirected irreflexive graph that is strongly regular with the same parameters as its classical analogue modeled over $M_4$. 
This is verified in our Software's Section \href{https://colab.research.google.com/drive/11rctqnjD3bxMejfEKWnJfr2qUy6oUYA3?usp=sharing}{4.4}.  
\end{exs}

\begin{remark}\label{rmk:non-Cayleys}
    Neither the {\bf Petersen graph} (strongly regular classical graph with parameters $(10,3,0,1)$) or the {\bf Hoffman-Singleton graph} (unique strongly regular classical graph with parameters $ (50,7,0,1)$) are Cayley graphs of groups \cite[page 10]{MR4458537}.
    Therefore, our current methods do not readily apply to them, and it is an interesting problem to determine whether they can be deformed. 
    This problem might be tractable as neither of these have quantum symmetry  \cite[Table 2]{MR4082994}. 

    Similarly, the {\bf $M22$ graph} is the unique strongly regular classical graph with parameters $(77, 16, 0, 4).$ 
    Since $77$ is square-free, even by realizing $M22$ as a Cayley graph of some group of order $77=11\times 7$, our methods would not apply as it does not contain subgroups of central type. 
\end{remark}

\subsection{Quantum spin models}\label{sec:QuSpinModels}
In this section we introduce the notion of a \emph{quantum} spin model, and give a condition so it yields a specialization of the Kauffman polynomial. 
We shall also use $3$-point regular quantum graphs to obtain explicit examples. 

\begin{defn}\label{defn: QuSpinMod}
    Let $X$ be a quantum set with value $\delta>0.$
    A {\bf symmetric quantum spin model} with {\bf the loop variable} $d\in \bbC$ with $d^2= \delta^2$, {\bf modulus} $a\in\bbC^\times$, and {\bf Boltzmann weights} matrices $W_+,W_-\in \End(\ell^2(X))$  is the $5$-tuple 
    $$
    \left(X,\ W_+,\ W_-,\ d,\ a\right)\!,
    $$
    satisfying the following relations:
    \begin{enumerate}[label={\bf (QSM\arabic*)}]
        \item\label{item:modulus} $\id_X\star W_+ = a\id_X,$ \ and \  $\id_X\star W_- = a^{-1}\id_X,$
        \item\label{item:1regular} $J\circ W_+ = d a^{-1}J,$ \  and \  $J\circ W_- = d a J,$
        \item\label{item:SchurInverses} $W_+\star W_- = J,$
        \item\label{item:CircInverses} $W_+\circ W_- = d^2 \id_X,$
        \item\label{item:StarTriangle} The  {\bf star-triangle equations}: 
        $$
        W_+ \circ m \circ (W_+\otimes W_-)=d \cdot m \circ (W_+\otimes W_-)\circ (m\otimes \id_X) \circ(\id_X\otimes W_-\otimes \id_X)\circ(\id_X\otimes m^\dag).
        $$
        That is, in diagrams 
        \begin{align*}
        d\cdot\tikzmath
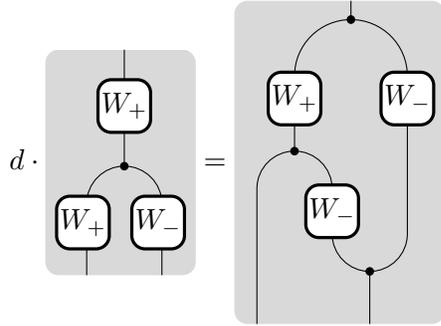
{
        \begin{scope}
        \clip[rounded corners=5pt] (-.8,-2.5) rectangle (1.2,.5);
        \fill[\AColor] (-.8,-2.5) rectangle (1.2,.5);
        \end{scope}
        \draw (-.25,-2.5) -- (-.25,-1.5);
        \draw (.75,-2.5) -- (.75,-1.5);
        \draw (-.25,-1.55) arc (180:0:.5cm);
        \draw (.25,-1) -- (.25,.5);
        \filldraw (.25,-1.05) circle (.05cm);
        \roundNbox{fill=white}{(-.3,-1.8)}{.35}{0}{0}{$W_+$}
        \roundNbox{fill=white}{(.7,-1.8)}{.35}{0}{0}{$W_-$}
        \roundNbox{fill=white}{(.25,-.25)}{-.35}{0}{0}{$W_+$}
        }
        = \tikzmath{
        \begin{scope}
        \clip[rounded corners=5pt] (-1.3,-1.5) rectangle (1.5,2.8);
        \fill[\AColor] (-1.3,-1.5) rectangle (1.5,2.8);
        \end{scope}
        \draw (0,-.3) arc (180:360:.5cm) -- (1,1.5);
        \draw (.5,-1.5) -- (0.5,-.8);
        \draw (-1, -1.5) -- (-1, .3) arc (180:0:.5);
        \draw (-.5, 1.3) -- (-.5,.8);
        \draw (-.5,1.8) arc (180:0:.75cm);
        \draw (.25,2.5) -- (.25,2.8);
        \filldraw (-.5,.8) circle (.05cm);
        \filldraw (.25,2.55) circle (.05cm);
        \filldraw (.5,-.8) circle (.05cm);
        \roundNbox{fill=white}{(-.5,1.5)}{.35}{0}{0}{$W_+$}
        \roundNbox{fill=white}{(1,1.5)}{.35}{0}{0}{$W_-$}
        \roundNbox{fill=white}{(0, 0)}{.35}{0}{0}{$W_-$}
        }.
        \end{align*}
    \end{enumerate}
    \noindent We call the quantum spin model {\bf symmetric} in case $W_{\pm}^* = W_{\pm}.$
\end{defn}

Definition \ref{defn: QuSpinMod} quantizes \cite[Proposition 2]{MR1188082}, extending it to the case where the underlying set $X$ is quantum and hence has no ``points''. 
That is, whenever $X= (\bbC^N, \tr)$ this recovers the usual notion and captures all classical examples. 
The Relations (QSM) come from invariance under the Shaded Reidemeister Moves, following the standard notation from the literature (c.f. \cite{MR1247144}).
That is, \ref{item:modulus} and \ref{item:1regular} reflect invariance under type-$\rm{I}$ moves, while \ref{item:SchurInverses} and \ref{item:CircInverses} reflects invariance under type-$\rm{II}$ moves.
Relation \ref{item:StarTriangle} comes from invariance under type-$\rm{III}$ moves. 
In fact, it can be verified by a diagrammatic computation that  \ref{item:StarTriangle} is equivalent to \cite[Equation (12)]{MR1247144} 
$$
R_1R_2R_1 = R_2R_1R_2,
$$
where $R_1:= W_+\otimes \id_X$ and $R_2:= d(\id_X\otimes m)\circ(\id_X\otimes W_-\otimes \id_X)\circ(m^\dag\otimes \id_X).$

To a quantum spin model, one can assign the C*-subalgebras of $\End(\ell^2(X))$
\begin{align}
    \langle W_+, J\ |\ \circ, \dag\rangle \quad \text{ and }\quad \langle W_-, I\ |\ \star, *\rangle,
\end{align}
for which there exists a \emph{unique} $*$-isomorphism i.e. a \emph{duality} $\Psi$ \cite[Proposition 3]{MR1188082} given by
\begin{align}\label{eqn:SelfDuality}
     W_+&\mapsto dW_-,  &&{W_-\mapsto dW_+},\\
     J&\mapsto d^2J,     &&I\mapsto J.\nonumber
\end{align}
Depending on normalization, this map is proportional to the quantum Fourier transform. 

The additional requirement for a quantum spin model to give an evaluation of the Kauffman polynomial is the  {\bf Kauffman exchange relation} (c.f. \cite[Equation (18)]{MR1188082}): 
\begin{align}\label{eqn:KauffExchange}
    W_+ + \varepsilon \cdot W_- = z (d\id_X + \varepsilon J) \qquad \exists z\in\bbC,    
\end{align}
where $\varepsilon$ is a chosen sign $\varepsilon\in\{\pm 1\}.$

We now describe how to obtain examples of quantum spin models from certain quantum graphs. 
\begin{construction}[{\bf Quantum spin models from quantum graphs}]\label{const:BoltzmannWeights}
    Let $\cG= (X, \hat{T})$ be an irreflexive undirected $3$-point regular (Definition \ref{defn:3ptReg}) quantum graph with parameters 
    $$(\delta^2,\ k,\ \lambda,\ \mu,\ q_0,\ q_1,\ q_2,\ q_3).$$ 
    Let $\{k, s,r\}  = \mathsf{Spec}(\hat{T})\subset \bbR,$ where $k\geq s,r$ and one of $r$ or $s$ is negative.
    The odered collection $\{A_0, A_1, A_2\}=\{\id_X, \hat{T}, \hat{T}^c\}$ is a {\bf cocomutative association scheme with $m=2$ classes} over $\ell^2(X)$ \cite[Definition 2.5]{2024arXiv240406157G} (c.f. \cite[\S2.3]{MR1188082}).
    (Here, we consider $\hat{T^c}$ as the irreflexive complement of $\hat{T}.$) 
    Then, the Bose-Mesner algebra 
    $$
    \cA:=\mathsf{span}_\bbC\{\id_X, \hat{T}, \hat{T}^c\}
    $$
    is a $3$-dimensional C*-algebra with respect to both the $\circ$ and $\star$ operations. 
    In the $\circ$ picture, we denote by $\{E_0, E_1, E_2\} = \{J, T, T^c\}$ the corresponding ordered basis of minimal $\circ$-idempotents of $\cA\subset \End(\ell^2(X))$.
    That is, we also have
    $$
    \cA=\mathsf{span}_\bbC\{J, T, T^c\},
    $$
    which in turn are the spectral projections of $\hat{T}.$ 

    The change of basis matrix from $\{E_i\}_{j=0}^2$ to $\{A_i\}_{i=0}^2$ is given by 
    \begin{align}
        P:=\begin{pmatrix}
            1 & k & \delta^2-k-1\\
            1 & s & -1-s\\
            1 & r & -1-r
        \end{pmatrix}\!.
    \end{align}
    That is, $P = [\Psi]_{\{E_j\}}^{\{A_i\}}$, and we say {\bf the scheme is formally self-dual} if there is an order of the eigenvalues of $\hat{T}$ and the projections $\{E_j\}$ so that 
    \begin{align}\label{eqn:SelfDualScheme}
    P^2 = \delta^2\id_X.    
    \end{align}
    
    Consider the matrices in $\cA$ given by 
    \begin{align}\label{eqn:BoltzmannWeights}
        W_+:=  t_0\id_X + t_1\hat{T} + t_2\hat{T}^c \qquad \text{and} \qquad W_-:=  t_0^{-1}\id_X + t_1^{-1}\hat{T}^c + t_2^2\hat{T}.
    \end{align}
    where $t_0, t_1, t_2\in \bbC^\times$ are to be determined. 
    In the classical setting, in \cite[Proposition 5]{MR1188082}, Jaeger establishes necessary and sufficient conditions on $(X\cong \bbC^V, a, d, W_+, W_{-})$ to determine is a spin model in the sense of Definition \ref{defn: QuSpinMod}, imposing constraints on the coefficients $t_0, t_1, t_2.$ 
    Letting $\varepsilon\in \{\pm 1\}$ be a choice of sign, the remaining coefficients are determined from solutions $t\in \bbC^{\times}$ of equation \cite[Equation (33)]{MR1188082}
    \begin{align*}
        s^2 + (r+1)^2 - \varepsilon s (r+1)(t^2 + t^{-2}) = 1.
    \end{align*}
    The remaining relevant parameters are given by 
    \begin{align*}
        a &= -s \varepsilon t + (r+1)t^{-1},\\
        t_0&=a,\\
        t_1&= \varepsilon t,\\
        t_2 &= t^{-1},\\
        d&= \varepsilon(r-s),\\
        z&= t + \varepsilon t^{-1}.
    \end{align*}
    Here, $z$ is the parameter from the Exchange Relation Equation (\ref{eqn:KauffExchange}), and $d$ satisfies $d^2 = \delta^2$.    
\end{construction}
Our definition of a quantum spin model, thus provides a \emph{matrix free approach} to spin models (c.f. \cite{MR1323744}) by tracing out images of knot diagrams, where positive crossings are replaced by $W_+$ and negative crossings by $W_-.$  

All deformations obtained in Section \ref{sec:Applications} of classical graphs known to yield evaluations of the Kauffman polynomial \cite[Theorem 3.1.10]{2019arXiv190208984E}; i.e. $9P_q, 16Cl_q$ and $HS_q$ give examples of quantum symmetric spin models for the Kauffman polynomial. 
In fact, to construct the respective Boltzmann weights $W_+$ and $W_-$, aided by Theorem \ref{thm:QuInvariants}, since all the parameters of the spin model are determined by spectral and $3$-point regularity data, we can import them from \cite[\S3.6]{MR1188082}.
Of course, to build the Boltzmann weights, one should use the adjacency operators $\hat{T}_q$, which differ from their classical analogues $\hat{T}.$

\begin{ex}[{\bf the quantum square $A_3$}]\label{ex:SpinQSquare}
From Examples on $M_2$ from Example \ref{ex:LaplaciansM2}, we know that the irreflexive quantum square $A_3$ and its complement $A_2,$ which is quantum isomorphic to the disjoint union of two cliques $\cK_2,$ yield a pair of $3$-point regular quantum graphs. 
By \cite[\S 3.6.2]{MR1188082} and an application of our Theorem \ref{thm:QuInvariants} we get 
$$
k= 2,\ s= 0,\ r = -2,\ \varepsilon \in\{\pm1\},\ d = -2\varepsilon,\ t\in \bbC^\times,\ a = -t^{-1},\ z = t+\varepsilon t^{-1}. 
$$
This yields $t_0 = a,\ t_1=\varepsilon t,\ t_2 = t^{-1}.$ 
Here, 
$$
P = \begin{pmatrix}
    1 & 2 & 1\\
    1 & 0 & -1\\
    1& -2 & 1
\end{pmatrix}.
$$
The Boltzman weights from Equation (\ref{eqn:BoltzmannWeights}) become
\begin{align*}
    W^q_+ =\begin{pmatrix}
        0              & 0         & 0         & 2\varepsilon t\\
        0              &-2t^{-1}   & 0         & 0 \\
        0              & 0         & -2t^{-1}  & 0\\
        2\varepsilon t & 0         & 0         & 0
    \end{pmatrix},
     \qquad
    W^q_- =\begin{pmatrix}
        0                   & 0         & 0         & 2\varepsilon t^{-1}\\
        0                   &-2t       & 0         & 0 \\
        0                   & 0         & -2t       & 0\\
        2\varepsilon t^{-1} & 0         & 0         & 0
    \end{pmatrix}.
\end{align*}
By construction, $W^q_+$ and $W^q_-$ give quantum spin models as in Definition \ref{defn: QuSpinMod}, and can also be directly checked by computation. 
These furthermore satisfy the Kauffman Exchange relation for Equation (\ref{eqn:KauffExchange}) with $z$ as above. 

Notice these differ from their classical analogues given by 
\begin{align*}
    W_+ =\begin{pmatrix}
        -t^{-1}          & \varepsilon t   & t^{-1}        & \varepsilon t\\
        \varepsilon t    &-t^{-1}          & \varepsilon t & t^{-1} \\
        t^{-1}           & \varepsilon t   & -t^{-1}       & \varepsilon t\\
        \varepsilon t    & t^{-1}          & \varepsilon t & -t^{-1}
    \end{pmatrix},
     \qquad
    W_- =\begin{pmatrix}
        -t                 & \varepsilon t^{-1}  & t                    & \varepsilon t^{-1}\\
        \varepsilon t^{-1} &-t                   & \varepsilon t^{-1}   & t \\
        t                  & \varepsilon t^{-1}  & -t                   & \varepsilon t^{-1}\\
        \varepsilon t^{-1} & t                   & \varepsilon t^{-1}   & -t
    \end{pmatrix}.
\end{align*}

The classical $W_+$ and $W_-$, their topological interpretation, along with their corresponding evaluations of the Kauffman polynomial are accounted for in \cite[Proposition 5]{MR1247144}. Verifications of these claims may be found in Section \href{https://colab.research.google.com/drive/11rctqnjD3bxMejfEKWnJfr2qUy6oUYA3?usp=sharing}{7.2} of our Software.   
\end{ex}

\begin{ex}[{\bf spin model from $9P_q$}]\label{ex:SpinQPaley}
In \cite[\S 3.6.4]{MR1188082}, Jaeger obtains the coefficients for the classical Paley graph $9P$ with parameters
\begin{align*}
    k =4,\ s = 1,\ r = -2,\ \varepsilon = -1,\ d = 3,\ t^2+t^{-2} =1 ,\ a = t^3,\ z = t+\varepsilon t^{-1},\ 
    P = \begin{pmatrix}
    1 & 4 & 4\\
    1 & 1 & -2\\
    1& -2 & 1
\end{pmatrix}.    
\end{align*}

Solving for $t$ we obtain $t = \pm e^{\pm i\pi/6}$). 
Therefore, $W^q_+$ for the $9P_q$ graph is given by
\begin{align*}
\begingroup
\setlength\arraycolsep{2pt}
  \begin{pmatrix}
    t^3-2\varepsilon t  &   0               &   0                 &                  0   & -\varepsilon t+2t^{-1}&   0                 &   0                 &   0                 &   -\varepsilon t+2t^{-1} \\
    0                   &t^3+\varepsilon t  &   0                 &                  0   &   0                   &-\varepsilon t-t^{-1}&-\varepsilon t-t^{-1}&  0                  &    0                     \\
    0                   &0                  &t^3+\varepsilon t    & -\varepsilon t-t^{-1}&   0                   &   0                 &   0                 &-\varepsilon t-t^{-1}&    0                     \\    
    0                   &0                  &-\varepsilon t-t^{-1}&t^3+\varepsilon t     &   0                   &   0                 &   0                 &-\varepsilon t-t^{-1} &    0                     \\
-\varepsilon t+2t^{-1}  &0                  &     0               &           0          &t^3-2\varepsilon t     &  0                  &   0                 &   0                 &    -\varepsilon t+2t^{-1}\\    
    0                   &-\varepsilon t-t^{-1}&    0              &            0         &      0                &t^3+\varepsilon t    &-\varepsilon t-t^{-1}&   0                 &    0                     \\
    0                   &-\varepsilon t-t^{-1}&     0             &             0        &      0                &-\varepsilon t-t^{-1}&t^3+\varepsilon t    &          0          &       0                  \\
    0                   &       0           &-\varepsilon t-t^{-1}&-\varepsilon t-t^{-1} &         0             &      0              &        0            &t^3+\varepsilon t    &         0                \\  
  -\varepsilon t+2t^{-1}&       0           &       0             &       0              &-\varepsilon t+2t^{-1} &      0              &      0              &      0              &t^3-2\varepsilon t        \\   
  \end{pmatrix}\!.
  \endgroup
\end{align*}
By \ref{item:CircInverses}, $W^q_- = d^2(W^q_+)^{-1}.$
By construction, these $W^q_+$ and $W^q_-$ give quantum spin models as in Definition \ref{defn: QuSpinMod}, and this can also be directly checked by computation (c.f. Software's Section \href{https://colab.research.google.com/drive/11rctqnjD3bxMejfEKWnJfr2qUy6oUYA3?usp=sharing}{7.6}). 
These furthermore satisfy the Kauffman Exchange relation for Equation (\ref{eqn:KauffExchange}) with $z$ as above. 

With the same parameters, the classical counterpart is given by
\begin{align*}
    W_+=\begin{pmatrix}
        t^3             &\varepsilon t &\varepsilon t &\varepsilon t& t^{-1}        & t^{-1}        &\varepsilon t  & t^{-1}        & t^{-1}      \\
        \varepsilon t   &t^{3}         &\varepsilon t & t^{-1}      &\varepsilon t  &t^{-1}         & t^{-1}        &\varepsilon t  & t^{-1}      \\
        \varepsilon t   &\varepsilon t &t^3           & t^{-1}      & t^{-1}        &\varepsilon t  & t^{-1}        &t^{-1}         &\varepsilon t\\
        \varepsilon t   & t^{-1}       &t^{-1}        & t^{3}       &\varepsilon t  &\varepsilon t  &\varepsilon t  & t^{-1}        & t^{-1}      \\ 
        t^{-1}          &\varepsilon t & t^{-1}       &\varepsilon t& t^3           &\varepsilon t  & t^{-1}        &\varepsilon t  & t^{-1}      \\
        t^{-1}          & t^{-1}       &\varepsilon t &\varepsilon t&\varepsilon t  & t^3           & t^{-1}        &t^{-1}         &\varepsilon t\\
        \varepsilon t   & t^{-1}       & t^{-1}       &\varepsilon t& t^{-1}        &t^{-1}         &t ^3           &\varepsilon t  &\varepsilon t\\
        t^{-1}          &\varepsilon t & t^{-1}       & t^{-1}      &\varepsilon t  & t^{-1}        &\varepsilon t  & t^3           &\varepsilon t\\
        t^{-1}          & t^{-1}       &\varepsilon t & t^{-1}      & t^{-1}        &\varepsilon t  &\varepsilon t  &\varepsilon t  &t^3
    \end{pmatrix},    
\end{align*}
with $W_- = d^2(W_+)^{-1}.$ 
The resulting classical spin model is known to yield the square of Kauffman's Bracket.
Compare with Section \href{https://colab.research.google.com/drive/11rctqnjD3bxMejfEKWnJfr2qUy6oUYA3?usp=sharing}{7.7} of our Software.
\end{ex}

\begin{ex}[{\bf spin model from $16Cl_q$}]\label{ex:SpinQClebsch}
By \cite[\S 4]{MR1188082} and \cite[pp 92, 93]{MR1247144}, and our Theorem \ref{thm:QuInvariants}, the quantum graph $Cl_q$
has  parameters 
$$
k =5,\ s = -3,\ r = 1,\ \varepsilon = \pm1,\ d = 4\varepsilon,\ t = \pm\sqrt{-\varepsilon} ,\ a = -t^{-1},\ z = 0,\ 
    P = \begin{pmatrix}
    1 & 5 & 10\\
    1 & -3 & 2\\
    1& 1 & -2
\end{pmatrix}.    
$$
The corresponding Boltzmann weights $W_\pm^\varepsilon$ as in Equation (\ref{eqn:BoltzmannWeights}) determine a quantum spin model, as can be seen by directly verifying Definition \ref{defn: QuSpinMod}.

That $z=0$ means that this evaluation of the Kauffman polynomial is topologically trivial. 
We observe, as can readily be verified by computation, that for all choices of $\varepsilon$ and $t$, $W_\pm^\varepsilon$ is a  {\bf quantum Hadamard matrix} over $M_4$ in the sense of  \cite{MR4855314}. See \href{https://colab.research.google.com/drive/11rctqnjD3bxMejfEKWnJfr2qUy6oUYA3?usp=sharing}{7.8} of our Software.
\end{ex}

\begin{ex}[{\bf spin model from $HS_q$}]\label{ex:SpinQHS}
By \cite[\S 3.6.5]{MR1188082}, and our Theorem \ref{thm:QuInvariants}, the quantum graph $HS_q$ determines a quantum spin model with the following parameters 
\begin{align*}
    k = 22,\ s = -8,\ r= 2,\ \varepsilon = -1,\ d=-10,\ t = \varphi,\ a = -\varphi^5,\ z =1,\ P = \begin{pmatrix}
    1 & 22 & 77\\
    1 & -8 & 7\\
    1& 2 & -3
\end{pmatrix}\!.
\end{align*}
Here, $\varphi = (1+\sqrt(5))/2$ is the golden ratio.
\end{ex}

We now consider the graphs $G_3$ and $G_4$ from Example \ref{ex:M3TripleRegularity} which are also $3$-point regular. 
\begin{ex}[{\bf spin model from $G_4$}]\label{ex:SpimModG4}
We already know that $G_4$ and its complement $G_6$ are both connected.
The regularity parameters of $G_4$ are $(3,\ 6/8,\ 6/8,\ 3/8,\ -3/8,\ 3/8,\ -3/8).$
Following a similar analysis to \cite[\S3.4]{MR1188082} we obtain the parameters
$$
    k=3,\ s = 3/2,\ r= -3/2,\ \varepsilon = \pm1,\ d = -3\varepsilon,\ t=\pm\sqrt{-\varepsilon},\ a = t^{-1},\ z = 0,\ P = \begin{pmatrix}
    1 & 3  & 5\\
    1 &3/2& -5/2\\
    1 &-3/2& 1/2
\end{pmatrix}\!.
$$
We see the corresponding association scheme is formally self-dual.

Notice that $\varepsilon t + t^{-1}=0$ and $\dfrac{3}{2}(\varepsilon t - t^{-1}) =3\varepsilon t$ for $\varepsilon\in\{\pm1\}.$ 
Thus, the corresponding Bolzmann weight is given by 
$$
W_+^\varepsilon = 3\cdot\begin{pmatrix}
    t^{-1}&   0  &   0   &   0   &  0  &  0   &   0   &   0   &  0\\
    0      &   0  & 0     &\varepsilon t&   0   &   0   &   0  &0   &0\\
    0      &   0  &  0     &    0   &   0& 0  & \varepsilon t &  0& 0\\
    0      & \varepsilon t& 0 & 0  & 0 &  0 &  0 &  0 &  0\\
    0      & 0 & 0 & 0 & t^{-1} & 0 & 0 & 0 & 0\\
    0      & 0 & 0 & 0 & 0 & 0 & 0 & \varepsilon t & 0 \\
    0      & 0 & \varepsilon t &0 & 0 & 0 &0 &0 &0\\
    0      & 0 & 0 & 0 & 0 &\varepsilon t & 0 & 0 & 0\\
    0      & 0 & 0 & 0 & 0 &0  &0 &0 & t^{-1}
\end{pmatrix}.
$$
It is evident that $W_+^\varepsilon\circ W_+^\varepsilon = -9\varepsilon\id$, and this too determines $W_-^\varepsilon$ by \ref{item:CircInverses}.

That $z=0$ means that this evaluation of the Kauffman polynomial is topologically trivial. 
That is, the Kauffman Exchange relation (\ref{eqn:KauffExchange}) becomes
$W_+ +\varepsilon W_- = 0.$
We observe, as can readily be verified by computation, that for all choices of $\varepsilon$ and $t$, $W_+^\varepsilon$ is a quantum Hadamard matrix over $M_3$. 
See Section \href{https://colab.research.google.com/drive/11rctqnjD3bxMejfEKWnJfr2qUy6oUYA3?usp=sharing}{7.4} of our Software for explicit computations. 
\end{ex}

\begin{remark}
Consider the graph $G_3$ and its complement $G_7$ defined over $M_3$ from Example \ref{ex:M3TripleRegularity}, whose regularity parameters mimic the classical graph pair consisting of the disjoint union of $3$ copies of $\cK_3$ and the complete tripartite graph $\cK_{3,3,3}$. 

Recall that $G_3$ is not connected, and so we shall not attempt to construct Boltzmann weights. 
However, the associated parameters are given by 
$$
    k=6,\ s = 0,\ r= -3,\ \varepsilon = \pm1,\ d = -3\varepsilon,\ P = \begin{pmatrix}
    1 & 6  & 2\\
    1 & 0 & -1\\
    1 &-3 & 2
\end{pmatrix}\!,
$$
where $P^2 = d^2\id.$
Nevertheless, \cite[Equation (33)]{MR1188082} does not determine $t$, and some numerical testing was inconclusive to give appropriate weights to build Boltzmann weights. 
The choice $r=0$, $s=-3$ does not give the necessary condition $P^2 =d^2\id.$ 
\end{remark}

\begin{remark}
    Using the Boltzmann weight matrices, one can obtain invariants for knots and links as in \cite[\S3.3]{MR990215} by presenting them as the trace closure of braids. 
    The invariant is then the (normalized) trace of the braid representation obtained from a given Boltzmann weight, which corresponds to the associated \emph{partition function}.

    All the spin models from Examples \ref{ex:SpinQSquare}, \ref{ex:SpinQPaley}, \ref{ex:SpinQClebsch} and \ref{ex:SpinQHS} with classical and quantum counterparts yield the same evaluations on knots and links. 
    This can be seen numerically directly on small examples such as the Trivial and the Trefoil knot, as well as the Hopf link, but is more cleanly visible diagrammatically from the properties of bubbling. 
\end{remark}

\bigskip

We close this manuscript with the question of whether there are other  quantum spin models (for the Kauffman polynomial) and quantum Hadamard matrices arising from quantum graphs.

\bibliographystyle{amsalpha}
\bibliography{bibliography}
\end{document}